\documentclass[11pt,reqno]{article}
\usepackage{bbding}
\usepackage{pifont}
\usepackage{amsmath}
\usepackage{mathrsfs}
\usepackage{amsfonts,bm}
\usepackage{amsthm}
\usepackage{epsfig}
\usepackage[]{harpoon}
\usepackage[active]{srcltx}
\usepackage{indentfirst,latexsym}
\usepackage{psfrag}
\usepackage[]{amssymb}
\usepackage{cite}
\usepackage{caption}
\usepackage{color}

\newcommand{\new}{\newcommand*}\new{\rnew}{\renewcommand*}
\new{\newe}{\newenvironment*}\new{\stl}{\setlength}
\stl{\textwidth}{155mm}\stl{\textheight}{22cm}\stl{\headheight}{0cm}
\stl{\topmargin}{0cm}\stl{\oddsidemargin}{0.5cm}\stl{\evensidemargin}{0cm}
\rnew{\arraystretch}{1.1}\rnew{\baselinestretch}{0.95}
\renewcommand{\thefootnote}{\ding{73}}
\newtheorem{thm}{Theorem}
\newtheorem{thmABC}{Theorem}

\newtheorem{lem}{Lemma}
\newtheorem{prop}{Proposition}
\newtheorem{defn}{Definition}
\newtheorem{rem}{Remark}

\newcommand{\eps}{\varepsilon}

\newcommand{\dps}{\displaystyle}

\newcommand{\fr}{\frac}

\newcommand{\pa}{\partial}

\numberwithin{equation}{section}
\newe{keywords}
   {\begin{quote}{\bf Keywords:}}
      {\end{quote}}
\newe{AMS}
   {\begin{quote}{\bf MSC subject classification 2020:}}
      {\end{quote}}
\newe{MSC}{\vspace*{5mm}
    {\noindent\it Mathematics Subject Classification(2020):}}{}
\new{\sect}[1]{\section{#1}\setcounter{equation}{0}
 \setcounter{thm}{0}\setcounter{lmm}{0}\setcounter{rmk}{0} }

\begin{document}

\title{Structurally stable singularities and Lipschitz stable optimal transport metrics for the compressible Euler equations}

\author{
Geng Chen\thanks{Department of Mathematics,
University of Kansas, Lawrence, KS
66045 ({\tt gengchen@ku.edu}).}
\and
Yanbo Hu\thanks{Department of Mathematics, Zhejiang University of Science and Technology, Hangzhou 310023, PR China ({\tt yanbo.hu@hotmail.com}).}
\and
Yannan Shen\thanks{Department of Mathematics,
University of Kansas, Lawrence, KS
66045 ({\tt yshen@ku.edu}).}
}

\rnew{\thefootnote}{\fnsymbol{footnote}}

\date{}

\maketitle
%%========================= abstract ==================================
\begin{abstract}
It is well known that solutions to the compressible Euler equations can develop singularities in finite time. In this paper, we carry out a detailed analysis on behaviors of
solutions up to the time of the first singularity for the one-dimensional compressible Euler equations with general smooth initial data.

Our main results consist of three parts. First, for an open dense set of $C^3$ initial data, we show that the solution of Euler equations is twice continuously differentiable except at most finitely many points when the first singularity happens, using Thom's Transversality Theorem.
Second, for any initial data in the  open dense set of $C^3$ functions given in the first result, we provide the precise asymptotic description of the solution in a semi-neighborhood in the $(x,t)$-plane of each singular point at the time of the first singularity, and verify that the solution has a cusp-type singularity with H\"older exponent $1/3$ at each singular point.
The proofs of the first two results are based on the representation of the solution in terms of a semilinear system.

Third, for smooth initial data with small BV norm, we construct two Finsler type optimal transport metrics, then under these metrics show that the solution depends Lipschitz continuously on the initial data up to the time of the first singularity, with uniformly bounded Lipschitz constants.  In particular, the $C^{1/3}$ generic singularity is stable in this sense. On the other hand, since our first two results hold for an open and dense set of initial data, any H\"older continuous cusp singularity with exponent other than $1/3$ is unstable under initial perturbations.

\end{abstract}

\begin{keywords}
Compressible Euler equations, Singularity, Generic regularity, Optimal Transport Metric, Transversality.
\end{keywords}

\begin{AMS}
35L65, 35L67, 76N15
\end{AMS}

\tableofcontents

%%$$$$$$$$$$$$$$$$$$$$$$$$$$$$$$$$ section 1 %$$$$$$$$$$$$$$$$$$$$$$$$$$$$$$$
\section{Introduction}\label{S1}

The one-dimensional compressible isentropic Euler equations in Lagrangian coordinates read that \cite{CF, Daf}
\begin{align}\label{1.1}
\left\{
\begin{array}{l}
\dps \tau_t-u_x=0, \\[6pt]
\dps u_t+p(\tau)_x=0, \\[6pt]
\dps p(\tau)=\fr{1}{\gamma}\tau^{-\gamma},
\end{array}
\right.
\end{align}
where $\tau>0$ is the specific volume, $u$ is the velocity, and $p(\tau)$ is the pressure equation of state for polytropic gases with adiabatic exponent $\gamma>1$.
We consider system \eqref{1.1} with the initial data
\begin{align}\label{1.3}
\tau(x,0)=\tau_0(x),\qquad u(x,0)=u_0(x).
\end{align}

As a fundamental model in fluid dynamics, the compressible Euler equations have a long history of research. The most prominent feature and major difficulty of this system lies in the fact that, even for sufficiently smooth and small initial data, its solutions may develop singularities in finite time due
to the quasilinear structures. The formation and evolution of shock wave (gradient blowup) starting from smooth initial data have been intensively studied.
It has become evident that a precise characterization of the solution at the time of the first singularity, also referred to as the pre-shock, is key to the investigation of this problem. In the present paper, for the compressible Euler equations \eqref{1.1} with general smooth initial data \eqref{1.3}, we study properties of the solution before and at the time when the
first singularity forms. Our result will provide the basis for the research in the next step: studying the development of shock waves.

\subsection{The background and prior results}\label{S11}
The phenomenon where solutions with even smooth initial data for the compressible Euler equations may develop gradient blow-ups in finite time originates in Bernhard Riemann's seminal paper in 1850s \cite{Rie}.
There exists an extensive mathematical literature on the compressible Euler equations, especially concerning the Riemann problem and the theory of weak solutions. For a broad overview, we refer the reader to the monographs and surveys \cite{Bre, CW, Daf, Maj} etc.
%Nevertheless, comparatively few rigorous results concern the mathematical analysis of the properties of solutions at singularity formation.
%and still fewer address the shock development and propagation problem.

The breakdown for smooth solutions to
the compressible Euler equations and the $2\times2$ reducible homogeneous hyperbolic systems with general initial data in one space dimension (1-D) were originated from Lax \cite{Lax} based on the characteristic method.
Subsequently, the more general hyperbolic systems were investigated among others in \cite{John, Liu1, Liu2, Maj}. These early results can be directly applied to the compressible Euler equations \eqref{1.1} with $\gamma\geq3$ to establish a complete dichotomy result, that is, the finite time singularity forms if and only if the initial compression exists. For the more physical relevant case when $1<\gamma<3$, an appropriate lower-bound estimate for the density is needed in order to resolve the degeneracy of the corresponding Riccati equations caused by the vanishing density. In this regard, the first author and his collaborators have carried out
a series of studies, in which they derived the optimal time-dependent lower bound for the density and obtained a complete picture of finite time singularity formation for system \eqref{1.1} with $1<\gamma<3$, see \cite{Chen1, Chen2, Chen3, Chen4} and the references therein.

The singularity formation for 3-D compressible Euler equations is much more involved. In 1980s, Sideris \cite{Sid} introduced a  virial-type integral method to replace the traditional characteristic analysis and proved that its smooth solution has a finite lifespan. It is worth pointing out that this result yields the conclusion of singularity formation, yet do not offer a precise description of the singularity structure.

Meanwhile, the constructive proofs of shock formation have received a lot attention over recent years. More precisely, these results, more related to the current paper, provide detailed descriptions on behaviors at the pre-shock shock, which are essential for analyzing the subsequent shock development.
In \cite{Leb}, Lebaud carried out pioneering constructive analysis for the 1-D isentropic Euler equations \eqref{1.1}. Under the assumptions of simple-wave initial data and a non-degenerate critical point, he proved that the solution develops a \(C^{1/3}\) cusp-type pre-shock. Moreover, he constructed weak solutions satisfying the Lax entropy condition for a short time after the blow-up time. Subsequently, Chen and Dong \cite{CD} improved this result by removing the simple-wave assumption, while Kong \cite{Kong} further expanded the theory to general \(2\times 2\) genuinely nonlinear hyperbolic systems. As the first contribution addressing the non-isentropic case, Yin \cite{Yin} demonstrated that, for smooth 3-D spherically-symmetric initial data with small perturbations, gradient breakdown occurs and \(C^{1/3}\) pre-shock structure is generated; furthermore, the solution admits a local extension near the blow-up point to a weak entropy solution.

In his seminal monograph \cite{Chr1}, Christodoulou developed a novel geometric framework for the study of the shock formation and development to the multi-dimensional (m-D) relativistic irrotational flows. Later, Christodoulou and Miao \cite{CM} applied this framework to study shock formation for the m-D compressible Euler equations with isentropic irrotational and small perturbed initial data. In \cite{CL}, Christodoulou and Lisbach employed the geometric framework to establish analogous to those in \cite{Yin} for the spherically symmetric isentropic Euler equations. The symmetry assumption was removed in \cite{Chr2}. Building upon this geometric framework,
Luk and Speck \cite{LS1, LS2} established shock formation results for the 2-D and 3-D isentropic Euler equations with non-zero vorticity. It is mentioned that the above construction results of shock formation are built upon perturbations of simple plane waves, which can give rise to many potential singularity types, yet these results make no distinction between these different scenarios.
In \cite{SV}, Shkoller and Vicol developed a new geometric framework for the description of the acoustic characteristic surfaces and then constructed the boundary structure of the maximal globally hyperbolic development for the m-D (multiple dimensional) compressible Euler equations.

Over the past several years, Buckmaster, Shkoller, Vicol and their collaborators through a series of papers \cite{BI, BDSV, BSV1, BSV2, BSV3, NRSV, NSV1, NSV2} have established a transformative constructive theory of shock formation and development for m-D compressible non-isentropic Euler equations, which has profoundly advanced the geometric and analytical understanding of finite-time gradient blowup arising from generic smooth initial data.
Initiating this systematic program, Buckmaster, Shkoller and Vicol \cite{BSV1} established shock formation for the 2-D isentropic compressible Euler equations with nontrivial vorticity by studying perturbations of purely azimuthal waves. The authors gave
a precise estimation of the blowup time and location and proved that the wave profile is of $C^{1/3}$ at the blowup time. This result provides a fundamental constructive framework for m-D shock development. In \cite{BDSV}, Buckmaster, Drivas, Shkoller and Vicol
explored the simultaneous development of shock waves and cusp singularities from smooth initial data for the 2-D compressible non-isentropic azimuthally-symmetric Euler equations.
The authors constructed a solution that initially exhibits a $C^{1/3}$ singularity, with a discontinuous shock arising immediately following the pre-shock. Moreover, it was also verified that as the shock front develops, two extra cusp-type singular characteristic surfaces associated with cusp-type singularities arise from the pre-shock. This finding resolves, in the azimuthally-symmetric setting, a conjecture put forward by Landau and Lifshitz \cite{LL}.

Furthermore, Buckmaster, Shkoller and Vicol constructed the first stable generic shocks and provided explicit estimations for the blow-up time and location for the 3-D isentropic Euler equations even with nonzero vorticity in \cite{BSV3}, with a significant extension to the non-isentropic setting given in \cite{BSV2}. On the other hand,
Buckmaster and Iyer \cite{BI} constructed unstable shocks for the 2-D isentropic compressible Euler equations in the azimuthally-symmetric setting, showing that the solution converges asymptotically to the unstable $C^{1/5}$ self-similar profile of Burgers' equation.
The main analytical tool employed throughout the works \cite{BI, BSV1, BSV2, BSV3, BDSV} for establishing the stability of an explicit blowup profile is the method of modulated self-similar analysis.

By utilizing a characteristics approach, the results of \cite{BDSV} was revisited by Neal, Shkoller, and Vicol in \cite{NSV1}. Under initial assumptions more general than those in \cite{BDSV}, the authors established the stable formation of $C^{1/3}$ cusps in finite time from smooth generic initial conditions for the 2-D compressible Euler equations in azimuthal symmetry.
In the context of the 1-D non-isentropic compressible Euler equations, Neal, Rickard, Shkoller, and Vicol \cite{NRSV} rigorously constructed a \(C^\nu\) pre-shock as a first singularity for all $\nu\in[1/2,1)$. Moreover, they also proved that this formation is
stable with respect to small perturbations of the initial data in the topology of $C^{1,\mu}$ with $\mu<(1-\nu)/\nu$. Recently, Neal,  Shkoller, and Vicol \cite{NSV2} further constructed the $C^{1/(2n+1)}$ cusp type pre-shocks for all integers $n\geq1$, and shown that the solutions are codimension-$(2n-2)$ stable in the space $W^{2n+2,\infty}$. In addition,
the finite-time vorticity blowups for smooth solutions of the 2-D and 3-D compressible Euler equations with smooth, localized, and nonvacuous initial data were verified in \cite{Chenj1, Chenj2}.

The first goal of this paper is to prove the existence of generic singularity for the initial value problem with general initial data and then following by a precise asymptotic description of the solution near the point where singularity forms. The main techniques we use include the transformation from the original system to a new semilinear system under characteristic coordinates and the Thom's transversality Theorem. A similar framework was first used in \cite{BC1,BHY} for the variational wave equation. Our result gives a new method, different from the one established by Buckmaster, Shkoller, Vicol and their collaborators, to consider the pre-shock solution, and especially it works for the general initial data.

Next, we study the Lipschitz stability of solutions until the first blowup, using the Finsler optimal transport metric. This result is based on the framework of Bressan and the first author in \cite{BC1,BC2,BHY}, also when they studied the variational wave equation.

More precisely, upon eliminating $u$, the Cauchy problem \eqref{1.1}-\eqref{1.3} can be reformulated as a Cauchy problem for the second-order quasilinear wave equation in $\tau$:
\begin{align}\label{1.1a}
\left\{
\begin{array}{l}
\dps \tau_{tt}-(c^2(\tau)\tau_x)_x=0,\\
\dps \tau(x,0)=\tau_0(x),\quad \tau_t(x,0)=u_{0}'(x),
\end{array}
\right.
\end{align}
where $c(\tau)=\sqrt{-p'(\tau)}=\tau^{-\fr{\gamma+1}{2}}$ is the speed of sound.
A large amount of research has been devoted to singularity formation and lifespan estimates for smooth solutions of the general second-order quasilinear wave equations, which will not be reviewed in the present paper. Interested readers may consult the relevant monographs \cite{Ali, Chr1, Chr3, John1, Maj, Spe} and the references therein.

As first noticed in \cite{CS} by the first and third authors, the equation \eqref{1.1a} is in fact closely related to the following quasilinear wave equation:
\begin{align}\label{1.4}
\varphi_{tt}-a(\varphi)(a(\varphi)\varphi_x)_x=0,
\end{align}
especially when we study the formation of singularity and the behavior of solution near the singularity.
The equation \eqref{1.4} is derived from a variational principle for nematic liquid crystals and is commonly referred to as the variational wave equation \cite{HS}. Here the wave speed $a$ is a smooth function of $\varphi$. The variational wave equation \eqref{1.4} with strictly positive wave speed has been intensively investigated, see partial works on the singularity formations of smooth solutions \cite{GHZ}, on the dissipative weak solutions \cite{BH, ZZ1, ZZ2}, on the conservative weak solutions \cite{BZ, BCZ, CS, HR}, and on the variational wave system \cite{CCD, ZZ3, ZZ4}. In particular, Bressan and the first author demonstrated in \cite{BC1} that, for an open dense set of $C^3$ initial data, the solution of \eqref{1.4} is piecewise smooth, whereas the gradient $\varphi_x$ blows up along finitely many characteristic curves. Their analysis relied on the energy-dependent coordinates introduced in \cite{BZ} and Thom's transversality theorem \cite{Thom}. We also refer the readers to the generic regularity result for solution of 2-D pressureless Euler equations, before, at and beyond the singularity (sticky particle) in \cite{BCH1,BCH2}.

Based on the generic regularity result in \cite{BC1}, in another paper \cite{BC2} they further constructed a Finsler type optimal transport metrics for \eqref{1.4} and proved that, under this metric, the flow is uniformly Lipschitz continuous on bounded subsets of $H^1$, the energy space for the variational wave equation. Subsequently, these results have been extended to general variational wave equations and related evolutionary equations \cite{CCS1, CCS2, CCCS}. For the open dense set of initial data constructed in \cite{BC1}, Bressan, Huang, and Yu \cite{BHY} carried out a detailed asymptotic description of the solution in a neighborhood of each singular point. Especially, they systematically analyzed the behaviors of singular points in conservative weak solutions and
derived the H\"older exponents corresponding to each class of cusp-type generic singularities.

\subsection{Results of the paper}\label{S12}

In this paper, we study the behavior of the solution at the time of the first singularity
for the one-dimensional compressible Euler equations \eqref{1.1} with
general initial data. Moreover, we construct two Finsler type optimal transport metrics for \eqref{1.1} up to the time when the first singularity occurs, and the flow is uniformly Lipschitz continuous under these metrics.

We outline our main results heuristically in what follows, while rigorous and precise statements are presented in later sections. From the complete dichotomy result in \cite{Lax,Chen3, Chen4}, we know that, for system \eqref{1.1}, singularity forms in finite time if and only if the initial data contains a compression point, that is $\varpi_{-,0}'(x)<0$ or $\varpi_{+,0}'(x)<0$ at some point $x$. Here $(\varpi_{-,0}(x), \varpi_{+,0}(x))$ are the initial data of the  Riemann invariants for \eqref{1.1} defined as $\varpi_{\pm,0}(x)=u_0(x)\pm 2(\tau_0(x))^{-(\gamma-1)/2}/(\gamma-1)$. We use $t^*$ to denote the first blow-up time for the smooth solution.

The first main result of this paper is stated in abbreviated form as follows, with its complete formulation provided in Theorem \ref{thm1} in Section \ref{S4}.
\begin{thmABC}\label{thmA}
Let $\bar{\tau}_0>0$ be a constant. There exists an open dense set of initial data
\begin{align}\label{1.5}
\mathcal{D}\subset\Big(C^3(\mathbb{R})\cap W^{1,1}(\mathbb{R})\Big)^2,
\end{align}
with
\begin{align}\label{1.6}
\begin{array}{c}
\dps\fr{1}{\overline{c}}\leq \tau_0(x)\leq \overline{c}, \qquad
\dps\min\big\{\inf \varpi_{-,0}'(x),\ \inf \varpi_{+,0}'(x)\big\}<0,
\end{array}
\end{align}
for some positive constant $\overline{c}$. For initial data $(\tau_0(x), u_0(x))$ satisfying $(\tau_0(x)-\bar{\tau}_0, u_0(x))\in \mathcal{D}$ and \eqref{1.6},
the solution $(\tau(x,t), u(x,t))$ of \eqref{1.1}-\eqref{1.3}, defined on $\mathbb{R}\times[0,t^*]$, is twice continuously differentiable in the complement of finitely many points on the line $t=t^*$.
\end{thmABC}
\begin{rem}
The requirement of initial data lying in $(W^{1,1})^2$ serves mainly to ensure that the solution is smooth in the far-field region. The last condition in \eqref{1.6} means that the initial data contain compression points.
\end{rem}
\begin{rem}
Theorem \ref{thmA} shows that, for initial data in the open dense set $\mathcal{D}$ with \eqref{1.6}, only finitely many singularities emerge on the line $t=t^*$.
\end{rem}

The second main result of this paper is given below, while its full formulation appears in Theorem \ref{thm2} of Section \ref{S5}.
\begin{thmABC}\label{thmB}
Consider initial data $(\tau_0(x), u_0(x))$ to satisfy $(\tau_0(x)-\bar{\tau}_0, u_0(x))\in \mathcal{D}$  and \eqref{1.6} and $(\tau(x,t), u(x,t))$ to be the solution to the Cauchy problem \eqref{1.1}-\eqref{1.3} on the domain $\mathbb{R}\times[0, t^*]$. Let $(x^*, t^*)$ be a singular point formed from the smooth solution $(\tau(x,t), u(x,t))$. Then $(x^*, t^*)$ is a singular point of Type-1 (defined in Section \ref{S5}). Moreover, if the gradient $\pa_x\varpi_-(x,t)$ blows up at the point $(x^*, t^*)$, then the other Riemann invariant $\varpi_+(x,t)$ is smooth up to the singular point $(x^*, t^*)$, and there exists a constant $\zeta\neq0$ such that for $t\leq t^*$
\begin{align}\label{1.7}
\begin{split}
\tau(x,t)=&\tau(x^*,t^*) +\zeta\cdot [c(\tau(x^*,t^*))(t-t^*)+(x-x^*)]^{\fr{1}{3}} \\ &+O(1)\Big(|t-t^*|+|x-x^*|\Big)^{\fr{4}{9}}, \\
u(x,t)=&u(x^*,t^*) +c(\tau(x^*, t^*))\zeta\cdot [c(\tau(x^*,t^*))(t-t^*)+(x-x^*)]^{\fr{1}{3}} \\ &+O(1)\Big(|t-t^*|+|x-x^*|\Big)^{\fr{4}{9}}.
\end{split}
\end{align}
\end{thmABC}
\begin{rem}
Theorem \ref{thmB} shows that, for initial data in the open dense set $\mathcal{D}$ with \eqref{1.6}, all singularities on the line $t=t^*$ are of a single type, that is Type-1 defined in Section \ref{S5}, representing the starting point of any generic singular curve in forward or backward direction in the $(x,t)$-plane.
More precisely, only one of the quantities $\pa_x\varpi_-(x,t)$ and $\pa_x\varpi_+(x,t)$ blows up at the generic singular point $(x^*, t^*)$.
\end{rem}
\begin{rem}
The solution $(\tau(x,t), u(x,t))$ is of $C^{1/3}$ H\"older class pre-shock at every singular point on $t=t^*$, which is consistent with the H\"older exponent obtained in the earlier works \cite{CD, Leb}. However, unlike \cite{CD, Leb} and other prior constructions of cusp-type singularities with H\"older exponent $1/3$ in relevant literatures, our result is established for generic initial data, rather than relying on constructions with special point-wise conditions imposed on the initial data.
\end{rem}
\begin{rem}
Since $\mathcal{D}$ is an open dense set in $(C^3(\mathbb{R})\cap W^{1,1}(\mathbb{R}))^2$, it follows that any cusp-type singularity starting from the general smooth initial data on the line $t=t^*$ with H\"older exponent other than $1/3$ is unstable. Moreover, the $C^{1/3}$ generic singularity is stable
in the topologies induced by the two metrics $d_{*}^{(\iota)}\ (\iota=i,ii)$ in Theorem \ref{thmC} below. Their stability follows from the smoothness of the terms $\mathcal{J}_{\ell}^\pm\ (\ell=1,\cdots,6)$ in \eqref{7.17} and the definitions of $d_{*}^{(\iota)}\ (\iota=i,ii)$ in \eqref{7.27}.
\end{rem}

The third main result of this paper is presented as follows, with its compete statement given in Theorems \ref{thm4} and \ref{thm5} in Section \ref{S7}.
\begin{thmABC}\label{thmC}
Let $\bar{\tau}_0, \bar{c}$ and $\eta$ be three positive constants and denote
\begin{align}\label{1.8}
\mathcal{D}_{\eta}=\left\{(\tau_0(x), u_0(x))\left|
\begin{array}{ll}
\dps (\tau_0-\bar{\tau}_0, u_0)\in \Big(C^1(\mathbb{R})\cap W^{1,1}(\mathbb{R})\Big)^2,\\[6pt]
\dps \fr{1}{\bar{c}}\leq \tau_0(x)\leq \bar{c},\\[10pt]
\dps \int_{\mathbb{R}}[|\varpi_{-,0}'(x)|+|\varpi_{+,0}'(x)|] {\rm d}x\leq \eta,
\end{array}
\right.
\right\}.
\end{align}
Let $(\tau(x,t), u(x,t))$ and $(\hat{\tau}(x,t), \hat{u}(x,t))$ be two solutions of the compressible Euler equations \eqref{1.1} with initial data in $\mathcal{D}_{\eta}$. Denote by $t^*$ and $\hat{t}^*$ the corresponding times at which their first singularities form. Set $t_d=\min\{t^*, \hat{t}^*\}$. Then, if $\eta$ is sufficiently small, for any time $t\leq t_d$, the two geodesic distances
$d_{*}^{(i)}(\cdot,\cdot)$ and $d_{*}^{(ii)}(\cdot,\cdot)$, defined in Definition \ref{def5} in Section \ref{S7}, satisfy
\begin{align}\label{1.9}
d_{*}^{(\iota)}((\tau, u),(\hat{\tau}, \hat{u}))\leq \widehat{C} d_{*}^{(\iota)}((\tau_0, u_{0}),(\hat{\tau}_0, \hat{u}_{0})),\quad (\iota=i,ii),
\end{align}
for some positive constant $\widehat{C}$. That is, under the metrics $d_{*}^{(\iota)}$,
the solution of \eqref{1.1} depends Lipschitz continuously on the initial data up to the formation time of the first singularity, provided that the parameter $\eta$ in the initial condition \eqref{1.8} is sufficiently small.
\end{thmABC}

%+++++++++++++++++++++++++++++++++++++++++++++++++++++++++++++++++++++
\begin{figure}[htbp]
\begin{center}
\qquad
\begin{minipage}[t]{0.45\textwidth}
\includegraphics[scale=0.45]{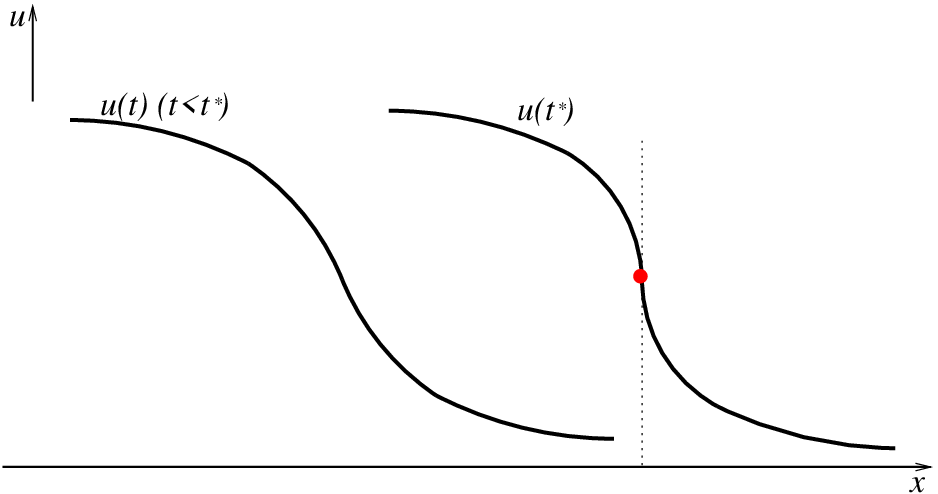}
\end{minipage}\qquad
\begin{minipage}[t]{0.45\textwidth}
\includegraphics[scale=0.42]{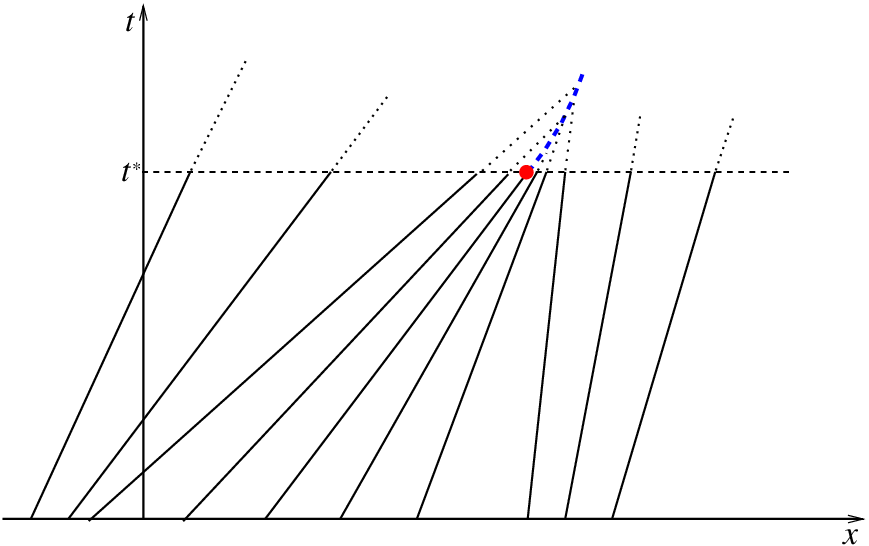}
\end{minipage}
\caption{Schematic illustration of singularity formation. Left: the profiles of $u$; Right: the compression wave. The formation of singularities lead the solution to lose Lipschitz stability with respect to the standard Sobolev metric.}
\label{fig1}
\end{center}
\end{figure}
%+++++++++++++++++++++++++++++++++++++++++++++++++++++++++++++++++++++

\begin{rem}
This conclusion requires a smallness condition on the initial data, which serves chiefly to establish the Lipschitz properties for the norms of tangent vectors of smooth solutions. See Lemma \ref{lem6} in Section \ref{S6}.
\end{rem}

\begin{rem}
It is worth pointing out that the above results (Theorems \ref{thmA}-\ref{thmC}) do not depend on the specific form of the wave speed $c$. Consequently, our framework can be applied to the general quasilinear wave equations $\psi_{tt}-(c^2(\psi)\psi_x)_x=0$ with generic conditions on $c(\cdot)$.
\end{rem}

We comment that
both metrics $d_{*}^{(i)}$ and $d_{*}^{(ii)}$ with main terms defined in \eqref{6.17a}  and \eqref{6.17}, respectively, include differences between derivatives of solutions. This makes these two metrics totally different from the $L^1$ metric, where the Lipschiz stability for the general small BV solution under the $L^1$ metric has been established in 1990s in \cite{BCo,ly, BLY, Bre}. Roughly speaking, two metrics in this paper are higher order metrics comparing to the $L^1$ metric. As we know, once a shock forms, i.e. when solution forms a discontinuity, the derivative of solution becomes a delta function (see Fig. \ref{fig1}), so one needs to add a transportation (or call it a shift) in the metric in order to obtain a Lipschitz metric. A natural choice is the optimal transport metric. See Remarks \ref{rem_12} for more details. The comparison between our metrics and the $L^1$ metric will be given in  Proposition \ref{prop3} in Section \ref{S7}. We also refer interesting readers to the H\"older stability theory under the $L^2$ metric for Euler equations in \cite{CKV,CFK1,CFK2}.
%but without an introduction since it is not very related to this paper.

\subsection{Structure and strategy of the paper}\label{S13}

We here present the structure of the paper.

\begin{itemize}

\item[] In Section \ref{S2}, we collect necessary preliminaries, including the Riemann invariants and classical dichotomy result for the compressible Euler equations \eqref{1.1}.

\item[] In Section \ref{S3}, we introduce a new set of dependent and
independent variables to transform \eqref{1.1} into an equivalent semilinear system for smooth solutions. Furthermore, we establish the existence of smooth solutions up to the boundary corresponding to the time $t=t^*$.

\item[] In Section \ref{S4}, we establish the generic structure of solutions to the semilinear
hyperbolic system at the singular points, and complete the proof of Theorem \ref{thm1}.

\item[] In Section \ref{S5}, we analyze the asymptotic behavior of the solution in a neighborhood of a singular point for the generic initial data, and complete the proof of Theorem \ref{thm2}.

\item[] In Section \ref{S6}, we construct two Finsler-type norms on tangent vectors for smooth solutions and demonstrate their Lipschitz continuity.

\item[] In Section \ref{S7}, we generalize the Lipschitz continuity previously established for smooth solutions to complete the proof of Theorems \ref{thm4} and \ref{thm5}. Finally, we compare these two metrics with the Sobolev metric.

\end{itemize}

In what follows, we give a more detailed explanation of the structure and research strategy of this paper.

In Section \ref{S2}, we first define the Riemann invariants $(\varpi_-, \varpi_+)$ for system \eqref{1.1}, together with their gradient variables $(R, S)=(\pa_x\varpi_-, \pa_x\varpi_+)$. Using these quantities, we derive a new system \eqref{2.9} in terms of $(R, S, \tau, u)$. The classical dichotomy result and the boundedness estimates for $\tau$ then follow from the specific structure of system \eqref{2.9}.

To handle the blow-up of $R$ and $S$, in Section \ref{S3}, we introduce new variables $\alpha\in[-\pi,\pi)$ and $\beta\in[-\pi,\pi)$ such that the blow-up of $R$ and $S$ correspond to $\alpha=-\pi$ and $\beta=-\pi$, respectively. Inspired by the work of Bressan and Zheng in \cite{BZ} and Chen and Shen in \cite{CS}, we introduce new dependent variables $(p, q)$ and coordinate variables $(X,Y)$.
This transforms system \eqref{2.9} into an equivalent semilinear system formulated in the coordinates $(X,Y)$ with unknowns $(\alpha, \beta, p, q, \tau, u)$, that is \eqref{3.14a} or \eqref{3.14}. Then we show this semilinear system admits a unique smooth solution up to the boundary $t(X,Y)=t^*$. In other words, the solution $(\alpha, \beta, p, q, \tau, u)(X,Y)$ to system \eqref{3.14a} remains smooth at $t(X,Y)=t^*$. Moreover, singularities occur when $\alpha(X,Y)=-\pi$ or $\beta(X,Y)=-\pi$, at which the Jacobian matrix $\pa(x,t)/\pa(X,Y)$ becomes non-invertible.

In Section \ref{S4}, we first show the local existence of smooth solutions to system \eqref{3.14a} near any given point $(X_0, Y_0)$ on the boundary $t(X,Y)=t^*$, which follows from the positive lower and upper bounds of $\tau$ at $t=t^*$. This implies that the solution $(\alpha, \beta, p, q, \tau, u)(X,Y)$ can be locally extended across the boundary $t(X,Y)=t^*$.
Then, motivated by the work of Bressan and Chen in \cite{BC1}, we prove that if $(X_0,Y_0)$ is a singular point on $t(X,Y)=t^*$ such that
$$
(\alpha, \alpha_X, \alpha_{XX})(X_0, Y_0)=(-\pi, 0, 0)\quad / \quad (\alpha, \beta, \alpha_{X})(X_0, Y_0)=(-\pi, -\pi, 0),
$$
there exist 3-parameter families of smooth solutions $(\tau^\theta, u^\theta, \alpha^\theta, \beta^\theta, p^\theta, q^\theta)$ to system \eqref{3.14a} such that $(\tau^0, u^0, \alpha^0, \beta^0, p^0, q^0)=(\tau, u, \alpha, \beta, p, q)$ and
$$
{\rm rank}\ D_\theta(\alpha^\theta, \alpha_{X}^\theta, \alpha_{XX}^{\theta})=3\quad / \quad {\rm rank}\ D_\theta(\alpha^\theta, \beta^\theta, \alpha_{X}^{\theta})=3.
$$
Based on the above conclusions and Thom's transversality theorem \cite{BC1, Thom}, we are able to show that the set consisting of all $C^2$ solutions $(\alpha, \beta, p, q, \tau, u)(X,Y)$ of system \eqref{3.14a} defined on a finite region, for which $\tau$ has positive upper and lower bounds, and which excludes the quantities given below, forms an open and dense subset of $C^2$.
\begin{align}\label{1.10}
\left\{
\begin{array}{r}
(\alpha, \alpha_X, \alpha_{XX})=(-\pi, 0, 0), \\
(\beta, \beta_Y, \beta_{YY})=(-\pi, 0, 0),
\end{array}
\right. \qquad
\left\{
\begin{array}{r}
(\alpha, \beta, \alpha_{X})=(-\pi, -\pi, 0), \\
(\alpha, \beta, \beta_{Y})=(-\pi, -\pi, 0).
\end{array}
\right.
\end{align}
The above results are presented in Lemmas \ref{lem4} and \ref{lem5}. Then Theorem \ref{thm1} (i.e. Theorem \ref{thmA}) follows directly from Lemma \ref{lem5}. In addition, we also obtain that, for the solution $(\alpha, \beta, p, q, \tau, u)(X,Y)$ of \eqref{3.14a} with the generic initial data,
none of the values in \eqref{1.10} are attainable at any point $(X,Y)$ on the boundary $t(X,Y)=t^*$.

In Section \ref{S5}, building upon the results established in Section \ref{S4}, we first classify the singular points satisfying $\alpha(X,Y)=-\pi$ on the boundary $t(X,Y)=t^*$ into three types:
\begin{itemize}

\item[] Type-1: Points where $\alpha=-\pi$, $\alpha_X= 0$, but $\beta\neq-\pi$, $\alpha_{XX}\neq0$.

\item[] Type-2: Points where $\alpha=-\pi$, but $\alpha_X\neq 0$, $\beta\neq-\pi$.

\item[] Type-3: Points where $\alpha=-\pi$, $\beta= -\pi$, but $\alpha_X\neq 0$, $\beta_Y\neq0$.

\end{itemize}
The singular points with $\beta(X,Y)=-\pi$ on the boundary $t(X,Y)=t^*$ can be classified analogously. Considering the aforementioned three types of singularities, we apply the governing system \eqref{3.14a} to perform a rigorous Taylor expansion analysis for $\tau(X,Y), u(X,Y), t(X,Y)$ and $x(X,Y)$ around the singular point $(X_0, Y_0)$. It is verified that, as
$t^*$ is the first occurrence time of singularities, only Type-1 singularities can emerge
on the boundary $t(X,Y)=t^*$, while Type-2 and Type-3 singularities are excluded.

In fact, three types of generic singularities are for the semi-linear system, corresponding to, the initial or ending point of the singularity curve (Type-1), the point on the singular curve (Type-2), the intersection point of two singular curves in two different families (Type-3).
Since the inverse transformation from the semilinear system
back to the original system only holds until the time of first blowup, when Type-2 generic
singularity won't happen.
And after a slight perturbation on the initial data, one can also avoid the case that two singularities in two directions both form exactly at the same point at the onset of the first singularity. This rules out the Type-3 singularity for solutions satisfying generic initial data (an open dense set $\mathcal{D}$). To study the development of shock wave, one can apply the Rankine-Hugoniot condition in the inverse transformation. We will not address this issue in the current paper.

For a Type-1 singular point $(X_0, Y_0)$, we utilize the Taylor expansions of $t(X,Y)$ and $x(X,Y)$ around $(X_0, Y_0)$ to express $X$ and $Y$ in terms of $t$ and $x$. Substituting the resulting expressions into the Taylor expansions of $\tau(X,Y)$ and $u(X,Y)$, we then obtain the expansions of $\tau$ and $u$ near the singular point $(x^*, t^*)$. With these explicit expressions, the proof of Theorem \ref{thm2} (i.e. Theorem \ref{thmB}) is thus completed.

\begin{figure}[htbp]
\begin{center}
\qquad
\begin{minipage}[t]{0.45\textwidth}
\includegraphics[scale=0.45]{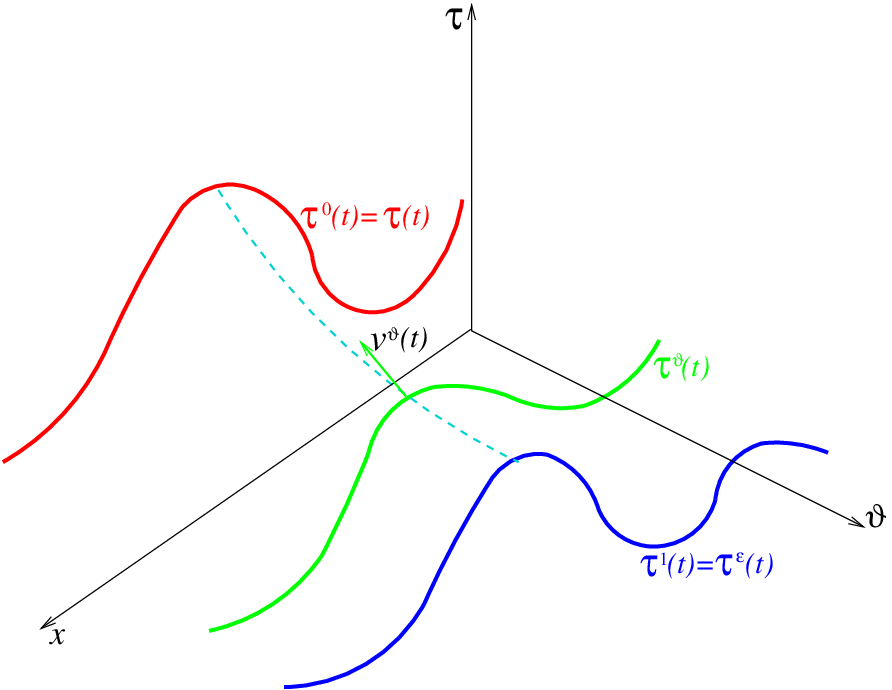}
\end{minipage}\qquad
\begin{minipage}[t]{0.45\textwidth}
\includegraphics[scale=0.45]{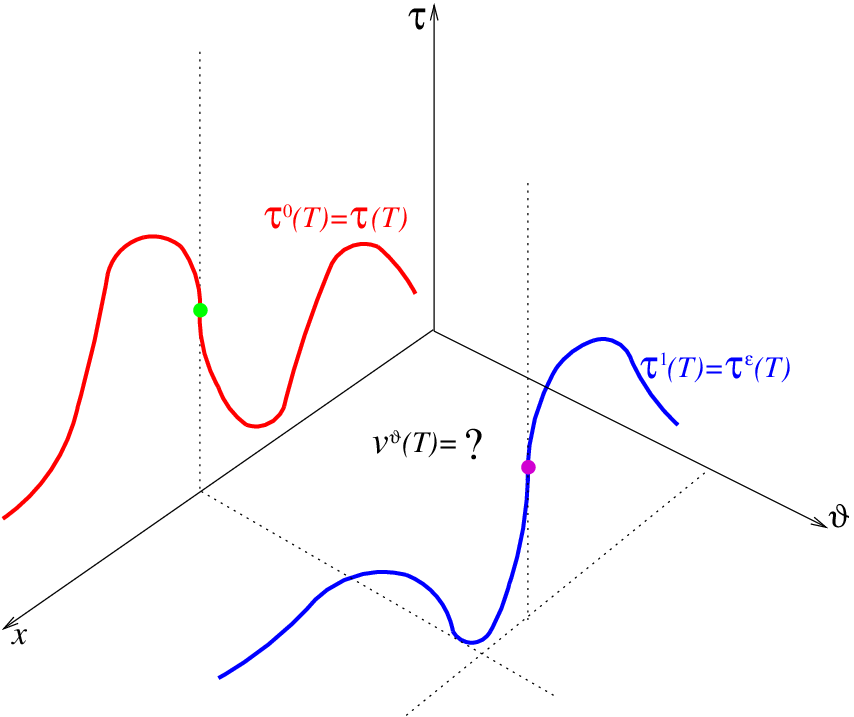}
\end{minipage}
\caption{Compare two solutions $\tau(x)$ and $\tau^\eps(x)$ at a given time $t$. Left: method of homotopy; Right: loss of regularity. Here ${\rm v}^\vartheta(t)={\rm d}\Upsilon^t/{\rm d}\vartheta$.}
\label{fig2}
\end{center}
\end{figure}

To account for the transportation cost, motivated by the work of Bressan and Chen in \cite{BC2} concerning the variational wave equation \eqref{1.4}, we are led to first construct the geodesic distance for smooth solutions in Section \ref{S6}. That is, for two given smooth solutions $(\tau(x,t), u(x,t))$ and $(\tau^\eps(x,t), u^\eps(x,t))$ for \eqref{1.1}-\eqref{1.3} defined on $\mathbb{R}\times [0,T]$, we consider all possible smooth paths $\Upsilon^t: \vartheta\mapsto (\tau^\vartheta, u^\vartheta)\ (\vartheta\in[0,1])$ satisfying $\Upsilon^t(0)=(\tau, u)$ and $\Upsilon^t(1)=(\tau^\eps, u^\eps)$. The length of each such path is measured by integrating the norm of the tangent vector ${\rm d}\Upsilon^t/{\rm d}\vartheta$, see the left figure in Fig. \ref{fig2}. Then we define the distance between $(\tau(x,t), u(x,t))$ and $(\tau^\eps(x,t), u^\eps(x,t))$ to be the infimum taken over the lengths of all such admissible paths as follows
\begin{align}\label{1.11}
d_*((\tau,u),(\tau^\eps, u^\eps))=\inf_{\Upsilon^t}\|\Upsilon^t\|:=\inf_{\Upsilon^t} \int_{0}^1\bigg\|\fr{{\rm d}\Upsilon^t}{{\rm d}\vartheta}\bigg\|_{(\tau^\vartheta, u^\vartheta)(t)}\ {\rm d}\vartheta,
\end{align}
where the subscript $(\tau^\vartheta, u^\vartheta)(t)$ is used to indicate that the the norm depends on the flow $(\tau, u)$. The crucial ingredient is how to define the Finsler norm $\|\cdot\|_{(\tau^\vartheta, u^\vartheta)(t)}$ so as to capture the behavior of the compressible Euler equations such that
\begin{align}\label{1.12}
\|\Upsilon^t\|\leq \widehat{C}\|\Upsilon^0\|,\qquad \forall\ t\in[0,T],
\end{align}
for some positive constant $\widehat{C}$. The inequality \eqref{1.12} leads to the Lipschitz continuity of the metric $d_*$ on the initial data.

%+++++++++++++++++++++++++++++++++++++++++++++++++++++++++++++++++++++
%+++++++++++++++++++++++++++++++++++++++++++++++++++++++++++++++++++++

\begin{figure}[htbp]
\begin{center}
\qquad
\begin{minipage}[t]{0.45\textwidth}
\includegraphics[scale=0.45]{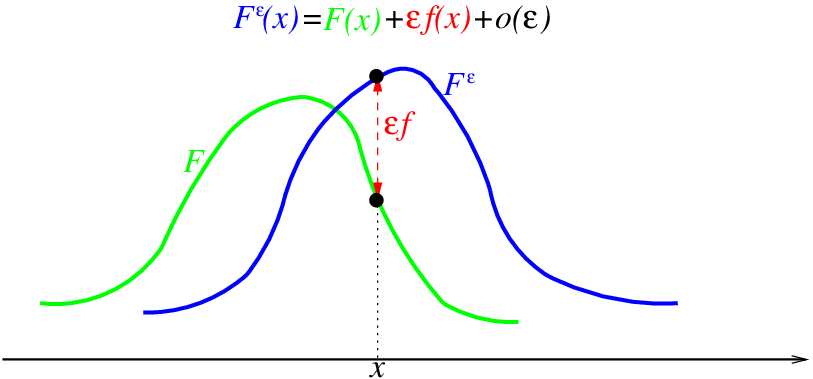}
\end{minipage}\quad
\begin{minipage}[t]{0.45\textwidth}
\includegraphics[scale=0.45]{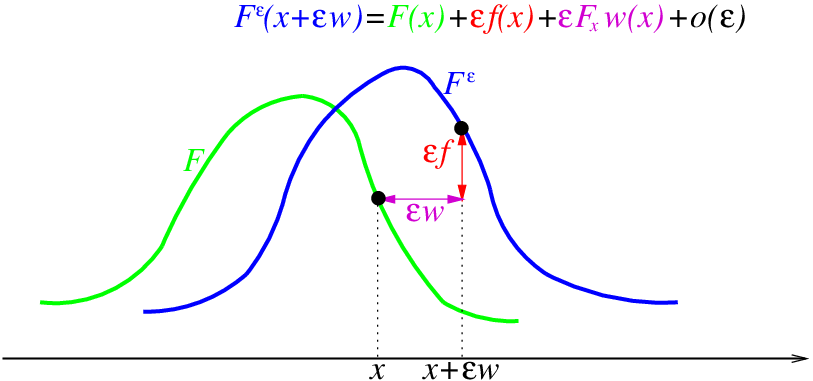}
\end{minipage}
\caption{A sketch of how to deform from the function $F$ to $F^\eps$. Left: a vertical shift $\eps f$; Right: a horizontal shift $\eps F_x w$ followed by a vertical displacement $\eps f$. The total shift as $\eps \tilde{f}=\eps (f+F_xw)$.}
\label{fig3}
\end{center}
\end{figure}

The Finsler norm $\|\cdot\|_{(\tau^\vartheta, u^\vartheta)(t)}$ is defined by measuring the cost of shifting from one solution to the other, in both forward and backward directions, with the densities $|R|$ and $|S|$, respectively. By the conservation laws for $(R, S)$, we can obtain that
\begin{align}\label{1.13}
\int_{\mathbb{R}}(|R|+|S|)\ {\rm d}x\leq \int_{\mathbb{R}}(|R_0(x)|+|S_0(x)|)\ {\rm d}x\leq \eta.
\end{align}
To circumvent the difficulties arising from absolute values, we replace the density $|R|$ with either $\sqrt{1+R^2}$ or $R^2/\sqrt{1+R^2}$, and the density $|S|$ with
either $\sqrt{1+S^2}$ or $S^2/\sqrt{1+S^2}$. When measuring the cost of transporting a function $F$ to $F^\eps$, we observe that the corresponding tangent flow $f$ only accounts for the vertical displacement between $F$ and $F^\eps$. To characterize the transportation cost more precisely, it is necessary to introduce an additional quantity, denoted by $w$, that measures the horizontal shift in the $x$-direction. It is therefore natural to take both vertical and horizontal shifts into account for estimating the transportation cost. More precisely, we measure the cost of transporting $F$ to $F^\eps$ by
\begin{align}\label{1.14}
\begin{split}
[{\rm change\ in}\ F]=&o(\eps)\ {\rm order\ of}\ \Big(F^\eps(x^\eps, t)-F(x, t)\Big) \\
=&\tilde{f}:=f(x,t)+F_x(x,t)w(x,t),
\end{split}
\end{align}
for any time $t\in[0, T]$. See Fig. \ref{fig3} for the illustration.

%+++++++++++++++++++++++++++++++++++++++++++++++++++++++++++++++++++++
%+++++++++++++++++++++++++++++++++++++++++++++++++++++++++++++++++++++

Then the basic structure of the Finsler norm takes the form
\begin{align}\label{1.15}
{\rm Terms\ in\  backward\ direction}\ \ +\ \ {\rm Terms\ in\ forward\ direction},
\end{align}
with
\begin{align}\label{1.16}
\begin{split}
&{\rm Terms\ in\  backward\ direction}  \\
=&\int_{\mathbb{R}}\bigg\{[{\rm change\ in}\ x] +[{\rm change\ in}\ \tau] +[{\rm change\ in}\ u] \\
&\qquad +[{\rm change\ in}\ \arctan R]\bigg\}\sqrt{1+R^2}\ {\rm d}x \\
&+\int_{\mathbb{R}}[{\rm change\ of\ base\ measure\ with\ density} \ \sqrt{1+R^2}]\ {\rm d}x.
\end{split}
\end{align}
The `Terms in forward direction' are symmetric. The structure in \eqref{1.16} merely gives a crude picture of the actual norm. Three additional modifications are required so that the metric satisfies the desired uniform Lipschitz property.

\noindent (i) Interaction potentials. From the governing system for the density $(\sqrt{1+R^2}, \sqrt{1+S^2})$ (i.e. \eqref{6.11}), we note that the forward or backward energy
might increase in the wave interaction due to the nonhomogeneous terms on the right-hand side of the system, even though the total energy is controlled. To overcome this difficulty, we introduce a pair of interaction potentials $\mathcal{W}^-$/$\mathcal{W}^+$ for backward/forward directions defined as
\begin{align}\label{1.17}
\begin{split}
\mathcal{W}^-=\eta+\int_{-\infty}^x\fr{S^2}{\sqrt{1+S^2}}(y,t)\ {\rm d}y,\quad \mathcal{W}^+=\eta+\int_{x}^\infty \fr{R^2}{\sqrt{1+R^2}}(y,t)\ {\rm d}y.
\end{split}
\end{align}
Serving as weight functions, they produce time-decay effects and thereby control the growth of energy.

\noindent (ii) Relative shifts. When a backward wave is horizontally shifted by $w$
and a forward wave by $z$, their relative shift equals $(w-z)$. This relative shift will induce new wave interactions. Hence, the cost functionals require appropriate adjustment to account for the relative-shift effect. Incorporating appropriate relative-shift term is an extremely delicate task. Minor modifications can substantially affect the construction of the metric and even destroy its Lipschitz property. Through careful computations, we add the following relative-shift contributions to [change in R] and [change in S] in \eqref{1.14}:
\begin{align}\label{1.18}
\begin{array}{l}
\dps [{\rm change\ in}\ R]=\tilde{r}:=r+R_x w-\fr{c'}{4c^2}RS(w-z),\\[6pt]
\dps [{\rm change\ in}\ S]=\tilde{s}:=s+S_x z-\fr{c'}{4c^2}RS(w-z).
\end{array}
\end{align}

\noindent (iii) Conserved quantities. Our calculations show that the foregoing metric still cannot achieve the desired Lipschitz property. To close the chain of estimates, we add terms generated by two conserved quantities $(r+R_x w+Rw_x)$ and $(s+S_x z+Sz_x)$ to the metric. Their origin lies in the conservation laws obeyed by $R$ and $S$.

The final two constructed `Terms in backward direction' are then given by
\begin{align}\label{1.19}
\begin{split}
&{\rm Terms\ in\  backward\ direction\ (i)}  \\
=&\int_{\mathbb{R}}\bigg\{|w| +\Big|\upsilon +\fr{Rw}{2c}-\fr{Sz}{2c}\Big| +\Big|\mu +\fr{Rw+Sz}{2}\Big|\bigg\}\mathcal{W}^-\sqrt{1+R^2}\ {\rm d}x \\
&+\int_{\mathbb{R}}|r+R_xw+Rw_x|\mathcal{W}^-\ {\rm d}x,
\end{split}
\end{align}
and
\begin{align}\label{1.20}
\begin{split}
&{\rm Terms\ in\  backward\ direction\ (ii)}  \\
=&\int_{\mathbb{R}}\bigg\{|w| +\Big|\upsilon +\fr{Rw}{2c}-\fr{Sz}{2c}\Big| +\Big|\mu +\fr{Rw+Sz}{2}\Big|\bigg\}\mathcal{W}^-\sqrt{1+R^2}\ {\rm d}x \\
&+\int_{\mathbb{R}}|r+R_xw+Rw_x|\mathcal{W}^-\ {\rm d}x  +\int_{\mathbb{R}}\fr{|\tilde{r}|}{\sqrt{1+R^2}}\mathcal{W}^-\ {\rm d}x \\
&+\int_{\mathbb{R}}\Big|\fr{R\tilde{r}}{\sqrt{1+R^2}} +\sqrt{1+R^2}\Big(w_x+\fr{c'}{4c^2}(w-z)S\Big)\Big|\mathcal{W}^-\ {\rm d}x.
\end{split}
\end{align}
Their explicit forms of the two metrics are presented in \eqref{6.17a} and \eqref{6.17} below. Under the smallness assumption on $\eta$, we establish their Lipschitz properties.

\begin{rem}
The two metrics constructed in this paper differ substantially from the metric established by Bressan and the first author in \cite{BC2} for the variational wave equation \eqref{1.4}. The main distinctions are summarized in the following three aspects.

First, our construction of the metrics draws heavily upon the conservation laws satisfied by the governing system, which captures the special structure of the compressible Euler equations.

Second, the density $(\sqrt{1+R^2}, \sqrt{1+S^2})$ adopted here does not admit an energy equation analogous to the one available for the variational wave equation.

Third, the relative-shift term demands delicate tuning, and its choice is entirely different from the corresponding quantity employed for the variational wave equation.
\end{rem}

To generalize the metrics constructed for smooth solutions in Section \ref{S6}, we first need to resolve the definition of the tangent vector ${\rm d}\Upsilon^t/{\rm d}\vartheta$ in \eqref{1.11}, which may not be well-defined due to the presence of singular points on the line $t=T$; see the right figure in Fig. \ref{fig2} for illustration. This potential obstacle is resolved by Theorem \ref{thm1} (i.e. Theorem \ref{thmA}), namely, we show that these metrics are well-posed for generic solutions by expressing them in terms of the $(X,Y)$ coordinates. Moreover, it is verified that the metrics are continuous with respect to $t$. Consequently, the existence of the generic singular points on the line $t=T$ does not affect their Lipschitz properties. These results are presented in Theorem \ref{thm4} in Section \ref{S7}. Finally, to define the geodesic distance
$d_{*}^{(\iota)}((\tau, u),(\hat{\tau}, \hat{u}))$ of the two solutions $(\tau, u)$ and $(\hat{\tau}, \hat{u})$ for $t\leq t_d$, we take the limit of the sequence $d_{*}^{(\iota)}((\tau^{(n)}, u^{(n)}),(\hat{\tau}^{(n)}, \hat{u}^{(n)}))$, where $(\tau^{(n)}, u^{(n)})$ and $(\hat{\tau}^{(n)}, \hat{u}^{(n)})$ are two sequences of generic solutions to \eqref{1.1}-\eqref{1.3} satisfying certain initial conditions. This completes the proof of Theorem \ref{thm5} (i.e. Theorem \ref{thmC}).

\section{Preliminaries}\label{S2}

For smooth solutions, system \eqref{1.1} can be rewritten as
\begin{align}\label{2.1}
\left(
\begin{array}{l}
\tau \\
u
\end{array}
\right)_t +
\left(
\begin{array}{cc}
0  &  -1 \\
-c^2(\tau)  &  0
\end{array}
\right)
\left(
\begin{array}{l}
\tau \\
u
\end{array}
\right)_x=0.
\end{align}
The two eigenvalues of \eqref{2.1} are
\begin{align}\label{2.2}
\lambda_-=-c(\tau),\qquad \lambda_+=c(\tau).
\end{align}
We introduce the Riemann invariants
\begin{align}\label{2.3}
\varpi_+=u+\fr{2}{\gamma-1}\tau^{-\fr{\gamma-1}{2}},\qquad \varpi_-=u-\fr{2}{\gamma-1}\tau^{-\fr{\gamma-1}{2}},
\end{align}
so that
\begin{align}\label{2.3a}
u=\fr{1}{2}(\varpi_++\varpi_-),\qquad \tau^{-\fr{\gamma-1}{2}}=\fr{\gamma-1}{4}(\varpi_+-\varpi_-).
\end{align}
Moreover, we can obtain the characteristic form of \eqref{2.1}
\begin{align}\label{2.4}
\left\{
\begin{array}{l}
\pa_t\varpi_--c(\tau)\pa_x\varpi_-=0, \\
\pa_t\varpi_++c(\tau)\pa_x\varpi_+=0.
\end{array}
\right.
\end{align}
It follows by \eqref{2.4} that
\begin{align}\label{2.5}
\varpi_-(x,t)=\varpi_-(x_-(0;x,t),0),\qquad \varpi_+(x,t)=\varpi_+(x_+(0;x,t),0),
\end{align}
where $x_\pm(s;x,t)$ are the $\lambda_\pm$-characteristic curves passing through point $(x,t)$ defined as follows
\begin{align}\label{2.6}
\left\{
\begin{array}{l}
\dps\fr{{\rm d}x_\pm(s;x,t)}{{\rm d}s}=\pm c(\tau(x_\pm(s;x,t),s)),\quad s\in[0,t], \\[5pt]
x_\pm(t;x,t)=x.
\end{array}
\right.
\end{align}
From \eqref{2.5}, we have
\begin{align}\label{2.7}
\inf \varpi_{\pm,0}(x)\leq \varpi_\pm(x,t) \leq \sup \varpi_{\pm, 0}(x),
\end{align}
where
$$
\varpi_{\pm,0}(x)=u_0(x)\pm \fr{2}{\gamma-1}(\tau_0(x))^{-\fr{\gamma-1}{2}}.
$$

We further introduce the gradient variables
\begin{align}\label{2.8}
R=\pa_x \varpi_-,\qquad S=\pa_x \varpi_+,
\end{align}
and then by \eqref{2.3a}
\begin{align}\label{2.8a}
\tau_t=u_x=\fr{R+S}{2},\qquad \tau_x=\fr{R-S}{2c(\tau)}.
\end{align}
From \eqref{2.4} and \eqref{2.8a}, the governing system of $(\tau, u, R, S)$ can be written as
\begin{align}\label{2.9}
\left\{
\begin{array}{l}
\dps R_t-cR_x=\fr{c'}{2c}(R-S)R, \\[8pt]
\dps S_t+cS_x=\fr{c'}{2c}(S-R)S, \\[8pt]
\dps \tau_t=\fr{R+S}{2}, \\[8pt]
\dps u_t=\fr{c}{2}(R-S),
\end{array}
\right.
\end{align}
subject to the initial data
\begin{align}\label{2.9a}
(\tau, u, R, S)|_{t=0}=(\tau_0(x), u_0(x), \varpi_{-,0}'(x), \varpi_{+,0}'(x)).
\end{align}
We also note the relation
\begin{align}\label{2.10}
\fr{c'(\tau)}{2c(\tau)}=-\fr{\gamma+1}{4\tau}.
\end{align}

Combining \eqref{2.4} and \eqref{2.9}, one can establish the following classical results \cite{Chen3, Chen4}.
\begin{prop}[\cite{Chen3, Chen4}] \label{prop-c}
Let $\gamma>1$. Assume that $(\varpi_{+,0}(x), \varpi_{-,0}(x))$ satisfy the conditions
\begin{align}\label{2.11}
\begin{split}
\inf(\varpi_{+,0}-\varpi_{-,0})(x)>0, \\
(\varpi_{+,0}, \varpi_{-,0})\in C^1,\quad \|(\varpi_{+,0}, \varpi_{-,0})\|_{C^1}<\infty.
\end{split}
\end{align}
Denote $R_0(x)=\varpi_{-,0}'(x)$ and $S_0(x)=\varpi_{+,0}'(x)$. Then the Cauchy problem \eqref{1.1}-\eqref{1.3} admits a unique global-in-time $C^1$ solution if and only if
\begin{align}\label{2.12}
R_0(x)\geq0,\quad S_0(x)\geq0, \quad \forall\ x\in \mathbb{R}.
\end{align}
Equivalently, the smooth solution of problem \eqref{1.1}-\eqref{1.3} develops singularities in finite time if and only if
\begin{align}\label{2.13}
R_0(x)<0,\quad {\rm or}\quad S_0(x)<0, \quad {\rm at\ some}\ x\in \mathbb{R}.
\end{align}
Moreover, let $t^*$ stand for the first time when one of $R$ and $S$ develops a singularity. Then the variables $(R, S, \tau)$ satisfy
\begin{align}\label{2.14}
\begin{split}
R(x,t)<&M,\qquad S(x,t)<M, \\
\inf_{\mathbb{R}\times[0,t^*]}\tau(x,t)>&\fr{1}{2}\Big(\fr{\gamma-1}{4}(\sup \varpi_{+,0}(x)-\inf \varpi_{-,0}(x))\Big)^{-\fr{2}{\gamma-1}},\\
\sup_{\mathbb{R}\times[0,t^*]}\tau(x,t)<&2(\sup \tau_0(x) +Mt^*).
\end{split}
\end{align}
where $M=1+\|R_0\|_{L^\infty} +\|S_0\|_{L^\infty}$.
\end{prop}

\begin{rem}\label{r1}
It is known by \eqref{2.14} and the definition of $c$ that there exists a uniform constant $\bar{c}>0$ such that
\begin{align}\label{2.22}
\fr{1}{\bar{c}}\leq c(\tau), -c'(\tau), c''(\tau)\leq \bar{c},
\end{align}
on the domain $\mathbb{R}\times[0,t^*]$.
\end{rem}

In this paper, we are focused on the behavior of the solution near the singular points, and accordingly assume that \eqref{2.13} holds. Throughout the paper, $t^*$ denotes the first time at which the smooth solution develops a singularity. Thus, the solution exists and is smooth for all $t<t^*$.
Furthermore, we will use $\overline{C}$ generically to denote a positive constant depending only on the constant $\overline{c}$, which may take different values at different places.

\section{The solution in a new coordinate plane}\label{S3}

In order to analyze the behaviors near singular points, we introduce a new set of dependent and independent variables to reformulate
the problem \eqref{2.9}-\eqref{2.9a} on the domain $\mathbb{R}\times[0,t^*]$ in this section.

\subsection{A semilinear hyperbolic system}\label{S31}

To deal with the unbounded values of $(R, S)$, we introduce a new pair of dependent
variables:
\begin{align}\label{3.1}
\alpha=2\arctan R,\qquad \beta=2\arctan S,
\end{align}
so that
\begin{align}\label{3.2}
\fr{1}{1+R^2}=\cos^2\fr{\alpha}{2},\quad \fr{R}{1+R^2}=\fr{1}{2}\sin \alpha,\quad \fr{1}{1+S^2}=\cos^2\fr{\beta}{2},\quad \fr{S}{1+S^2}=\fr{1}{2}\sin \beta.
\end{align}
Here and below, we restrict $\alpha$ and $\beta$ to lie in the interval $[-\pi,\pi)$, which follows from the facts from Proposition \ref{prop-c} that $R$ and $S$ have positive upper bounds and that  blow-up for $R$ and $S$ correspond, respectively, to $\alpha=-\pi$ and $\beta=-\pi$.
It suggests by \eqref{2.9} and \eqref{3.1} that
\begin{align}\label{3.3}
\begin{split}
\alpha_t-c\alpha_x=&\fr{c'}{c}\fr{R(R-S)}{1+R^2}=\fr{c'}{c}\fr{\sin\fr{\alpha}{2}}{\cos\fr{\beta}{2}} \sin\fr{\alpha-\beta}{2},\\
\beta_t+c\beta_x=&\fr{c'}{c}\fr{S(S-R)}{1+S^2}=\fr{c'}{c}\fr{\sin\fr{\beta}{2}}{\cos\fr{\alpha}{2}} \sin\fr{\beta-\alpha}{2}.
\end{split}
\end{align}

We further define a coordinate transformation $(x,t)\rightarrow(X,Y)$ as follows:
\begin{align}\label{3.4}
X=x_-(0;x,t),\qquad Y=-x_+(0;x,t).
\end{align}
Obviously, there hold
\begin{align}\label{3.5}
X_t-c(\tau)X_x=0,\qquad Y_t+c(\tau)Y_x=0.
\end{align}
Moreover, for any smooth function $f$, one finds by \eqref{3.5} that
\begin{align}\label{3.6}
\begin{split}
f_t-cf_x=&(f_XX_t+f_YY_t)-c(f_XX_x+f_YY_x)=(Y_t-cY_x)f_Y=-2cY_xf_Y,\\
f_t+cf_x=&(f_XX_t+f_YY_t)+c(f_XX_x+f_YY_x)=(X_t+cX_x)f_X=2cX_xf_X.
\end{split}
\end{align}
Now we introduce the new variables
\begin{align}\label{3.7}
p=\fr{\sqrt{1+R^2}}{X_x},\qquad q=\fr{\sqrt{1+S^2}}{-Y_x}
\end{align}
Note from \eqref{3.2} and \eqref{3.7} that
\begin{align}\label{3.8}
\fr{1}{X_x}=p\cos\fr{\alpha}{2},\qquad \fr{1}{-Y_x}=q\cos\fr{\beta}{2}.
\end{align}
Combining \eqref{3.3}, \eqref{3.6} and \eqref{3.8} yield
\begin{align}\label{3.9}
\begin{split}
\alpha_Y=&\fr{1}{-2cY_x}(\alpha_t-c\alpha_x)=\fr{c'}{2c^2}q\sin\fr{\alpha}{2} \sin\fr{\alpha-\beta}{2}, \\
\beta_X=&\fr{1}{2cX_x}(\beta_t+c\beta_x)=\fr{c'}{2c^2}p\sin\fr{\beta}{2}\sin\fr{\beta-\alpha}{2},
\end{split}
\end{align}
Furthermore, we compute by \eqref{3.3}, \eqref{3.5}, \eqref{3.6} and \eqref{3.8}
\begin{align}\label{3.10}
p_Y=&\fr{1}{-2cY_x}(p_t-cp_x) \notag \\ =&\fr{1}{2c}q\cos\fr{\beta}{2}\bigg\{\fr{\sin\fr{\alpha}{2}}{2X_x\cos^2\fr{\alpha}{2}} (\alpha_t-c\alpha_x) -\fr{1}{X_{x}^2\cos\fr{\alpha}{2}}(X_{xt}-cX_{xx})\bigg\}
\notag \\ =&\fr{1}{2c}q\cos\fr{\beta}{2} \bigg\{\fr{p\sin\fr{\alpha}{2}}{2\cos\fr{\alpha}{2}} \cdot\fr{c'}{c}\fr{\sin\fr{\alpha}{2}}{\cos\fr{\beta}{2}}\sin\fr{\alpha-\beta}{2} -p\cdot c'\fr{\tan\fr{\alpha}{2}-\tan\fr{\beta}{2}}{2c}\bigg\} \notag \\
=&\fr{c'}{4c^2}pq\cos\fr{\alpha}{2}\sin\fr{\beta-\alpha}{2},
\end{align}
and
\begin{align}\label{3.11}
q_X=&\fr{1}{2cX_x}(q_t+cq_x) \notag \\ =&\fr{1}{2c}p\cos\fr{\alpha}{2}\bigg\{\fr{\sin\fr{\beta}{2}}{-2Y_x\cos^2\fr{\beta}{2}} (\beta_t+c\beta_x) +\fr{1}{Y_{x}^2\cos\fr{\beta}{2}}(Y_{xt}+cY_{xx})\bigg\}
\notag \\ =&\fr{1}{2c}p\cos\fr{\alpha}{2} \bigg\{\fr{q\sin\fr{\beta}{2}}{2\cos\fr{\beta}{2}} \cdot\fr{c'}{c}\fr{\sin\fr{\beta}{2}}{\cos\fr{\alpha}{2}}\sin\fr{\beta-\alpha}{2} +q\cdot c'\fr{\tan\fr{\alpha}{2}-\tan\fr{\beta}{2}}{2c}\bigg\} \notag \\
=&\fr{c'}{4c^2}pq\cos\fr{\beta}{2}\sin\fr{\alpha-\beta}{2}.
\end{align}
Here we used the fact by \eqref{2.8a} and \eqref{3.1} that
$$
\pa_xc=c'\tau_x=\fr{c'}{2c}(R-S)=\fr{c'}{2c}\Big(\tan\fr{\alpha}{2}-\tan\fr{\beta}{2}\Big).
$$
In addition, utilizing \eqref{2.8a}, \eqref{3.6} and \eqref{3.8} again lead to
\begin{align}\label{3.12}
\begin{split}
\tau_X=&\fr{1}{2cX_x}(\tau_t+c\tau_x)=\fr{1}{2c}p\cos\fr{\alpha}{2} \Big(\fr{R+S}{2}+c\cdot\fr{R-S}{2c}\Big)  \\
=&\fr{1}{2c}p\cos\fr{\alpha}{2}R=\fr{1}{2c}p\sin\fr{\alpha}{2}, \\
\tau_Y=&\fr{1}{-2cY_x}(\tau_t-c\tau_x)=\fr{1}{2c}q\cos\fr{\beta}{2} \Big(\fr{R+S}{2}-c\cdot\fr{R-S}{2c}\Big)   \\
=&\fr{1}{2c}q\cos\fr{\beta}{2}S=\fr{1}{2c}q\sin\fr{\beta}{2}.
\end{split}
\end{align}
and
\begin{align}\label{3.13}
\begin{split}
u_X=&\fr{1}{2cX_x}\cdot cR=\fr{1}{2}p\sin\fr{\alpha}{2},\\
u_Y=&\fr{1}{-2cY_x}\cdot (-cS)=-\fr{1}{2}q\sin\fr{\beta}{2}.
\end{split}
\end{align}

We combine \eqref{3.9}-\eqref{3.13} to obtain a semilinear hyperbolic system for the variables $(\alpha, \beta, p, q, \tau)$ in the coordinates $(X,Y)$
\begin{align}\label{3.14a}
\left\{
\begin{array}{l}
\dps \alpha_Y=\fr{c'}{2c^2}q\sin\fr{\alpha}{2}\sin\fr{\alpha-\beta}{2}, \\[8pt]
\dps \beta_X=\fr{c'}{2c^2}p\sin\fr{\beta}{2}\sin\fr{\beta-\alpha}{2},\\[8pt]
\dps p_Y=\fr{c'}{4c^2}pq\cos\fr{\alpha}{2}\sin\fr{\beta-\alpha}{2},\\[8pt]
\dps q_X=\fr{c'}{4c^2}pq\cos\fr{\beta}{2}\sin\fr{\alpha-\beta}{2},\\[8pt]
\dps \tau_X=\fr{1}{2c}p\sin\fr{\alpha}{2}, \quad \Big({\rm or}\ \tau_Y=\fr{1}{2c}q\sin\fr{\beta}{2}\Big), \\[8pt]
\dps u_X=\fr{1}{2}p\sin\fr{\alpha}{2}, \quad \Big({\rm or}\ u_Y=-\fr{1}{2}q\sin\fr{\beta}{2}\Big),
\end{array}
\right.
\end{align}
or equivalently, the following system by substituting $h=\tau^{\fr{\gamma-1}{2}}$ for $\tau$
\begin{align}\label{3.14}
\left\{
\begin{array}{l}
\dps \alpha_Y=-\fr{\gamma+1}{4}hq\sin\fr{\alpha}{2}\sin\fr{\alpha-\beta}{2}, \\[8pt]
\dps \beta_X=-\fr{\gamma+1}{4}hp\sin\fr{\beta}{2}\sin\fr{\beta-\alpha}{2},\\[8pt]
\dps p_Y=-\fr{\gamma+1}{8}hpq\cos\fr{\alpha}{2}\sin\fr{\beta-\alpha}{2},\\[8pt]
\dps q_X=-\fr{\gamma+1}{8}hpq\cos\fr{\beta}{2}\sin\fr{\alpha-\beta}{2},\\[8pt]
\dps h_X=\fr{\gamma-1}{4}h^2p\sin\fr{\alpha}{2}, \quad \Big({\rm or}\ h_Y=\fr{\gamma-1}{4}h^2q\sin\fr{\beta}{2}\Big), \\[8pt]
\dps u_X=\fr{1}{2}p\sin\fr{\alpha}{2}, \quad \Big({\rm or}\ u_Y=-\fr{1}{2}q\sin\fr{\beta}{2}\Big).
\end{array}
\right.
\end{align}
We comment that one may use either $h_X$ or $h_Y$ and $u_X$ or $u_Y$ in \eqref{3.14} because there hold $h_{XY}=h_{YX}$ and $u_{XY}=u_{YX}$ by
\begin{align}\label{3.15}
h_{XY}=&\Big(\fr{\gamma-1}{4}h^2p\sin\fr{\alpha}{2}\Big)_Y \notag \\
=& \fr{(\gamma-1)^2}{8}h^3pq\sin\fr{\alpha}{2}\sin\fr{\beta}{2} =\Big(\fr{\gamma-1}{4}h^2q\sin\fr{\beta}{2}\Big)_X=h_{YX},
\end{align}
and
\begin{align}\label{3.15a}
u_{XY}=\Big(\fr{1}{2}p\sin\frac{\alpha}{2}\Big)_Y=0,\qquad u_{YX}=\Big(-\fr{1}{2}q\sin\frac{\beta}{2}\Big)_X=0.
\end{align}
Moreover, we see by \eqref{3.6} and \eqref{3.8} that
\begin{align}\label{3.16}
\left\{
\begin{array}{l}
\dps t_X=\fr{1}{2}h^{\fr{\gamma+1}{\gamma-1}}p\cos\fr{\alpha}{2}, \\[8pt]
\dps t_Y=\fr{1}{2}h^{\fr{\gamma+1}{\gamma-1}}q\cos\fr{\beta}{2},
\end{array}
\right.
\qquad
\left\{
\begin{array}{l}
\dps x_X=\fr{1}{2}p\cos\fr{\alpha}{2}, \\[8pt]
\dps x_Y=-\fr{1}{2}q\cos\fr{\beta}{2},
\end{array}
\right.
\end{align}
and subsequently
\begin{align}\label{3.17}
\fr{\pa(x,t)}{\pa(X,Y)}=\fr{1}{2}h^{\fr{\gamma+1}{\gamma-1}}pq\cos\fr{\alpha}{2}\cos\fr{\beta}{2}.
\end{align}

Note, system \eqref{3.14} for Euler equations was first found in \cite{CS}. In fact, a general system with a parameter $\lambda\in[0,1]$ was considered in \cite{CS}, where $\lambda=1$ and $\lambda=1/2$ are corresponding to equations \ref{1.1a} and \ref{1.4}, respectively.

\subsection{Smooth solutions for the new system}\label{S32}

We first consider the boundary conditions of system \eqref{3.14} in the
coordinates $(X,Y)$, corresponding to \eqref{2.9a} in the original coordinates $(x,t)$.
It is noted that
\begin{align}\label{3.18}
\begin{split}
R(x,0)=&R_0(x)=u_{0}'(x)+(\tau_0(x))^{-\fr{\gamma+1}{2}}\tau_{0}'(x), \\
S(x,0)=&S_0(x)=u_{0}'(x)-(\tau_0(x))^{-\fr{\gamma+1}{2}}\tau_{0}'(x).
\end{split}
\end{align}
Recalling the definition in \eqref{3.4}, we see that the initial line $t=0$ in the $(x,t)$ plane is transformed to the line $\Gamma_0: X+Y=0$ in the $(X,Y)$ plane. This line can be parameterized as $x\mapsto (\overline{X}(x), \overline{Y}(x))=(x,-x)$. Along $\Gamma_0$, we can assign the boundary data $(\overline{h}, \overline{u}, \overline{\alpha}, \overline{\beta}, \overline{p}, \overline{q})$ by setting
\begin{align}\label{3.19}
\left\{
\begin{array}{l}
\dps \overline{h}=(\tau_0(x))^{\fr{\gamma-1}{2}}, \\
\dps \overline{u}=u_0(x),
\end{array}
\right.
\quad
\left\{
\begin{array}{l}
\dps \overline{\alpha}=2\arctan R_0(x), \\
\dps \overline{\beta}=2\arctan S_0(x),
\end{array}
\right.
\quad
\left\{
\begin{array}{l}
\dps \overline{p}=\sqrt{1+R_{0}^2(x)}, \\
\dps \overline{q}=\sqrt{1+S_{0}^2(x)},
\end{array}
\right.
\end{align}
at each point $(x,-x)\in\Gamma_0$.

Let $\Omega_{t^*}$ denote the image of the domain $R\times[0,t^*]$ under the coordinate transformation $(x,t)\rightarrow(X,Y)$. We have the following lemma.
\begin{lem}\label{lem2}
The boundary value problem \eqref{3.14} and \eqref{3.19} admits a unique smooth solution on the domain $\Omega_{t^*}$ up to the boundary $t(X,Y)=t^*$.
\end{lem}
\begin{proof}
We observe that all right-hand side functions in system \eqref{3.14} are locally Lipschitz
continuous. Then, to establish the smooth solution on $\Omega_{t^*}$, it suffices to derive the a priori $L^\infty$-estimates of the solution $(h, u, \alpha, \beta, p, q)(X,Y)$ up to the boundary $t(X,Y)=t^*$.
According to \eqref{2.14}, we know that the variable $\tau(X,Y)$ and then the variable $h(X,Y)$ possess uniform positive upper and lower bounds on $\Omega_{t^*}$ and up to the boundary $t(X,Y)=t^*$. Then there holds
\begin{align}\label{3.20}
\fr{1}{\overline{C}}\leq h(X,Y)\leq \overline{C}.
\end{align}
Moreover, in view of the definitions of $(\alpha, \beta)$ in \eqref{3.1}, we deduce by \eqref{2.14} that
\begin{align}\label{3.21}
\alpha(X,Y),\beta(X,Y)\in[-\pi, 2\arctan M].
\end{align}
Furthermore, from \eqref{3.15a} and \eqref{3.19}, one obtains
\begin{align}\label{3.22}
\begin{split}
p(X,Y)\sin\frac{\alpha(X,Y)}{2}=&\sqrt{1+R_{0}^2(X)}\sin \arctan R_0(X) =R_0(X),\\
q(X,Y)\sin\frac{\beta(X,Y)}{2}=&\sqrt{1+S_{0}^2(Y)}\sin \arctan S_0(Y) =S_0(Y),
\end{split}
\end{align}
which implies that
\begin{align}\label{3.22a}
\begin{split}
u(X,Y)=u_0(-Y)+\int_{-Y}^{X} \fr{1}{2}R_0(x)\ {\rm d}x.
\end{split}
\end{align}
Now, putting \eqref{3.22} and \eqref{3.20} into the equations for $(p,q)$ in \eqref{3.14} yield
\begin{align}\label{3.23}
\begin{split}
p_Y=&-\fr{\gamma+1}{8}h\bigg\{p\cos^2\frac{\alpha}{2}\cdot \Big(q\sin\frac{\beta}{2}\Big) -q\cos\frac{\alpha}{2}\cos\frac{\beta}{2}\cdot \Big(p\sin\frac{\alpha}{2}\Big)\bigg\} \\
=&-\fr{\gamma+1}{8}h\bigg\{p\cos^2\frac{\alpha}{2}\cdot S_0(Y)  -q\cos\frac{\alpha}{2}\cos\frac{\beta}{2}\cdot R_0(X)\bigg\} \\
\leq & \overline{C}M(p+q),
\end{split}
\end{align}
and
\begin{align}\label{3.24}
\begin{split}
q_X=&-\fr{\gamma+1}{8}h\bigg\{q\cos^2\frac{\beta}{2}\cdot \Big(p\sin\frac{\alpha}{2}\Big) -p\cos\frac{\beta}{2}\cos\frac{\alpha}{2}\cdot\Big(q\sin\frac{\beta}{2}\Big)\bigg\} \\
=&-\fr{\gamma+1}{8}h\bigg\{q\cos^2\frac{\beta}{2}\cdot R_0(X) -p\cos\frac{\beta}{2}\cos\frac{\alpha}{2}\cdot S_0(Y)\bigg\} \\
\leq & \overline{C}M (p+q).
\end{split}
\end{align}
Here the facts $p\geq0, q\geq0$ and $|R_0(x)|\leq M, |S_0(x)|\leq M$ are employed. It concludes by \eqref{3.23} and \eqref{3.24} that
\begin{align}\label{3.25}
0\leq p(X,Y), q(X,Y)\leq M\exp\Big(\overline{C}M(X+Y)\Big),
\end{align}
for all $(X,Y)$ on $\Omega_{t^*}$. The proof of the lemma is complete.
\end{proof}

\begin{rem}\label{r2}
From Lemma \ref{lem2}, we see that the solution $(h, u, \alpha, \beta, p, q)(X,Y)$ of problem \eqref{3.14} and \eqref{3.19} is well-defined up to the boundary $t(X,Y)=t^*$. Singularities occur when $\alpha(X,Y)=-\pi$ or $\beta(X,Y)=-\pi$, at which points the Jacobian matrix $\pa(x,t)/\pa(X,Y)$ given in \eqref{3.17} is not invertible.
\end{rem}

\begin{rem}\label{r3}
At the point $(X, Y)\in \Omega_{t^*}$ where $\alpha\neq-\pi$ and $\beta\neq-\pi$, the Jacobian matrix is invertible, having a strictly positive determinant. Moreover, there hold
\begin{align}\label{3.26}
(x,t)=(x(X,Y), t(X,Y)).
\end{align}
Here the function $x(X,Y)$ and $t(X,Y)$ can be obtained by solving one of the equations in \eqref{3.16}. This is because the relations $t_{XY}=t_{YX}$ and $x_{XY}=x_{YX}$ hold by
\begin{align*}
t_{XY}=&\fr{\gamma+1}{16}h^{\fr{2\gamma}{\gamma-1}}pq\sin\fr{\alpha+\beta}{2}=t_{YX},\\[6pt]
x_{XY}=&\fr{\gamma+1}{16}hpq\sin\fr{\alpha-\beta}{2}=x_{YX}.
\end{align*}
\end{rem}

\begin{rem}\label{r4}
There are generally many distinct solutions to system \eqref{3.14a}, \eqref{3.16}, which produce the same solution $(\tau(x,t),u(x,t))$ to \eqref{1.1}-\eqref{1.3}. In fact, suppose that $(\tau, u, \alpha, \beta, p, q, x, t)$ is a particular solution to \eqref{3.14a}, \eqref{3.16}, and $\phi, \psi:\ \mathbb{R}\mapsto \mathbb{R}$ are two $C^2$ bijections satisfying $\phi'>0$ and $\psi'>0$. We define the new variables as follows:
\begin{align}\label{3.27}
\begin{split}
X=\phi(\widetilde{X}),\quad Y=\psi(\widetilde{Y}), \\
(\tilde{\tau}, \tilde{u}, \tilde{\alpha}, \tilde{\beta}, \tilde{x}, \tilde{t})(\widetilde{X}, \widetilde{Y}) =(\tau, u, \alpha, \beta, x, t)(X,Y), \\
\tilde{p}(\widetilde{X}, \widetilde{Y})=p(X,Y)\cdot \phi'(\widetilde{X}),\quad \tilde{q}(\widetilde{X}, \widetilde{Y})=q(X,Y)\cdot \psi'(\widetilde{X}).
\end{split}
\end{align}
One can easily check that the functions $(\tilde{\tau}, \tilde{u},  \tilde{\alpha}, \tilde{\beta}, \tilde{p}, \tilde{q}, \tilde{x}, \tilde{t})(\widetilde{X}, \widetilde{Y})$ also satisfy system \eqref{3.14a}, \eqref{3.16}. Furthermore, the two sets
\begin{align*}
\begin{split}
&\Big\{\Big(t(X,Y), x(X,Y),(\tau(X,Y), u(X,Y))\Big);\ (X,Y)\in \Omega_{t^*}\Big\}, \\
{\rm and}\ &\Big\{\Big(\tilde{t}(\widetilde{X}, \widetilde{Y}), \tilde{x}(\widetilde{X}, \widetilde{Y}), (\tilde{\tau}(\widetilde{X}, \widetilde{Y}), \tilde{u}(\widetilde{X}, \widetilde{Y}))\Big);\ (\widetilde{X}, \widetilde{Y})\in \widetilde{\Omega}_{t^*}\Big\},
\end{split}
\end{align*}
coincide, where $\widetilde{\Omega}_{t^*}=\{(\widetilde{X}, \widetilde{Y}):\ X=\phi(\widetilde{X}),\ Y=\psi(\widetilde{Y}),\ (X,Y)\in\Omega_{t^*}\}$. Specifically,
one may interpret the independent variable transformation in \eqref{3.27} as nothing more than a relabeling of the forward and backward characteristics in the solution $(\tau, u)$.
\end{rem}

\section{The generic structure of solutions}\label{S4}

In this section, we discuss the generic structure of smooth solutions to the semilinear hyperbolic system \eqref{3.14} near the points where $\alpha(X,Y)=-\pi$ or $\beta(X,Y)=-\pi$.

\subsection{Solutions near a given point}\label{S41}

Let a point $(X_0,Y_0)$ be given. We consider the line
$$
\Gamma_\kappa=\Big\{(X,Y):\ X+Y=\kappa\Big\},\qquad \kappa=X_0+Y_0,
$$
and arbitrarily assign the smooth values of $(\alpha, \beta, p, q)$ on $\Gamma_\kappa$:
\begin{align}\label{4.1}
\left\{
\begin{array}{l}
\alpha(s,\kappa-s)=\overline{\alpha}(s), \\
\beta(s,\kappa-s)=\overline{\beta}(s),
\end{array}
\right.
\qquad
\left\{
\begin{array}{l}
p(s,\kappa-s)=\overline{p}(s), \\
q(s,\kappa-s)=\overline{q}(s).
\end{array}
\right.
\end{align}
Let $M_1>0$ be a constant such that
\begin{align}\label{4.2}
|\overline{p}(s)|,|\overline{q}(s)|\leq M_1, \quad \forall\ s\in[X_0-1, X_0+1].
\end{align}
We set $h(X_0,Y_0)=h_0>0$ and
\begin{align}\label{4.3}
\delta=\min\bigg\{1, \fr{1}{(\gamma-1)M_1h_0}\bigg\},
\end{align}
and denote
\begin{align}\label{4.4}
\Omega_\delta=\Big\{(X,Y):\ |X-X_0|\leq \delta,\quad |Y-Y_0|\leq \delta \Big\}.
\end{align}

We have the following lemma.
\begin{lem}\label{lem3}
Let $\overline{h}(s)$ and $\overline{u}(s)$ be smooth solutions of the ODE problems
\begin{align}\label{4.5}
\left\{
\begin{array}{l}
\dps \fr{{\rm d}\overline{h}(s)}{{\rm d}s}=\fr{\gamma-1}{4}\overline{h}^2(s)\overline{p}(s)\sin\fr{\overline{\alpha}(s)}{2} +\fr{\gamma-1}{4}\overline{h}^2(s)\overline{q}(s)\sin\fr{\overline{\beta}(s)}{2}, \\[8pt]
\overline{h}(X_0)=h_0,
\end{array}
\right.
\end{align}
and
\begin{align}\label{4.5a}
\left\{
\begin{array}{l}
\dps \fr{{\rm d}\overline{u}(s)}{{\rm d}s}=\fr{1}{2}\overline{p}(s)\sin\fr{\overline{\alpha}(s)}{2} -\fr{1}{2}\overline{q}(s)\sin\fr{\overline{\beta}(s)}{2}, \\[8pt]
\overline{u}(X_0)=u(X_0, Y_0),
\end{array}
\right.
\end{align}
on the interval $[X_0-\delta, X_0+\delta]$. Then the semilinear hyperbolic system \eqref{3.14} with the boundary value $(\overline{h}(s), \overline{u}(s), \overline{\alpha}(s), \overline{\beta}(s), \overline{p}(s), \overline{q}(s))$ admits a unique smooth solution $(h, u, \alpha, \beta, p, q)(X,Y)$ on the domain $\Omega_\delta$. Moreover, the solution satisfies
\begin{align}\label{4.6}
\fr{1}{4}h_0\leq h(X,Y) \leq 4h_0,\quad \forall\ (X,Y)\in\Omega_\delta.
\end{align}
\end{lem}
\begin{proof}
We first show that the function $\overline{h}(s)$ satisfies
\begin{align}\label{4.7}
\fr{1}{2}h_0\leq \overline{h}(s) \leq 2h_0,\quad \forall\ s\in[X_0-\delta, X_0+\delta].
\end{align}
To prove \eqref{4.7}, one rewrites the equation in \eqref{4.5} as
\begin{align}
\fr{{\rm d}}{{\rm d}s}\Big((\overline{h}(s))^{-1}\Big) =-\fr{\gamma-1}{4}\Big(\overline{p}(s)\sin\fr{\overline{\alpha}(s)}{2} +\overline{q}(s)\sin\fr{\overline{\beta}(s)}{2}\Big), \notag
\end{align}
from which we obtain
\begin{align}\label{4.8}
\overline{h}(s)=h_0\bigg\{1 -h_0\int_{X_0}^s\fr{\gamma-1}{4} \Big(\overline{p}(s)\sin\fr{\overline{\alpha}(s)}{2} +\overline{q}(s)\sin\fr{\overline{\beta}(s)}{2}\Big)\ {\rm d}s \bigg\}^{-1}.
\end{align}
It follows by the choice of $\delta$ in \eqref{4.3} that
\begin{align*}
&h_0\bigg|\int_{X_0}^s\fr{\gamma-1}{4} \Big(\overline{p}(s)\sin\fr{\overline{\alpha}(s)}{2} +\overline{q}(s)\sin\fr{\overline{\beta}(s)}{2}\Big)\ {\rm d}s\bigg| \\
\leq & h_0\cdot \fr{\gamma-1}{4} \cdot 2M_1|s-X_0|\leq \fr{1}{2},
\end{align*}
from which and \eqref{4.8} we arrive at \eqref{4.7}.

Next we show the existence of smooth solutions on the domain $\Omega_\delta$. Similar to the proof of Lemma \ref{lem2}, the crux lies in deriving the a priori estimates of the solution. Recalling \eqref{3.15a} and applying \eqref{4.1} gives for all $(X,Y)\in\Omega_\delta$
\begin{align}\label{4.9}
p(X,Y)\sin\fr{\alpha(X,Y)}{2}&=p(X,\kappa-X)\sin\fr{\alpha(X,\kappa-X)}{2} \notag \\
&=\overline{p}(X)\sin\fr{\overline{\alpha}(X)}{2}.
\end{align}
Moreover, the equation for $h$ in \eqref{3.14} can be rewritten as
$$
\pa_X\Big(\fr{1}{h}\Big)=-\fr{\gamma-1}{4}p\sin\fr{\alpha}{2}.
$$
Integrating the above along variable $x$ from $(\kappa-Y, Y)$ to $(X,Y)$, and utilizing \eqref{4.9}, one obtains
\begin{align*}
h(X,Y)=h(\kappa-Y,Y)\bigg\{1 -h(\kappa-Y,Y)\int_{\kappa-Y}^X\fr{\gamma-1}{4} \overline{p}(s)\sin\fr{\overline{\alpha}(s)}{2}\ {\rm d}s \bigg\}^{-1},
\end{align*}
which together with \eqref{4.7} and \eqref{4.3} leads to
\begin{align}\label{4.10}
h(X,Y)\leq 2h_0\bigg\{1 -2h_0\cdot \fr{\gamma-1}{4}M_1\delta \bigg\}^{-1}\leq 4h_0,
\end{align}
and
\begin{align}\label{4.11}
h(X,Y)\geq \fr{1}{2}h_0\bigg\{1 +2h_0\cdot \fr{\gamma-1}{4}M_1\delta \bigg\}^{-1}\geq \fr{1}{4}h_0,
\end{align}
for all $(X,Y)\in\Omega_\delta$. Combining \eqref{4.10} and \eqref{4.11} yields \eqref{4.6}.

Based on the estimates for $h$ in \eqref{4.6}, we can repeat the estimation arguments for $(p,q)$ presented in \eqref{3.23}-\eqref{3.24} to establish the corresponding estimates of $(p,q)$ on the domain $\Omega_\delta$. The estimates of $(u, \alpha, \beta)$ can be derived by directly integrating their corresponding equations.
The proof of the lemma is complete.
\end{proof}

\subsection{Families of perturbed solutions}\label{S42}

Let $(X_0,Y_0)$ be a given point on the boundary $t(X,Y)=t^*$ such that
\begin{align}\label{4.12}
\alpha(X_0, Y_0)=-\pi,\qquad {\rm or}\qquad \beta(X_0, Y_0)=-\pi,
\end{align}
that is, $(X_0,Y_0)$ is a singular point. Set $\kappa=X_0+Y_0$. Thanks to the conclusion of Lemma \ref{lem2}, we know that
\begin{align}\label{4.13}
\fr{1}{\overline{C}}\leq h(X_0, Y_0)\leq \overline{C},
\end{align}
for some constant $\overline{C}>0$. Hence, by Lemma \ref{lem3}, assigned arbitrarily the smooth values of $(\overline{\alpha}, \overline{\beta}, \overline{p}, \overline{q})(s)$ on $\Gamma_\kappa$, there exists a smooth solution $(h, u, \alpha, \beta, p, q)$ near the point $(X_0,Y_0)$.

\begin{lem}\label{lem4}
Let $(X_0,Y_0)$ be a given singular point on the boundary $t(X,Y)=t^*$. Assume that $(h, u, \alpha, \beta, p, q)$ is a smooth solution of \eqref{3.14} near point $(X_0,Y_0)$.

\noindent {\rm(1)} If $(\alpha, \alpha_X, \alpha_{XX})(X_0, Y_0)=(-\pi, 0, 0)$, then there exists a 3-parameter family of smooth solutions $(h^\theta, u^\theta, \alpha^\theta, \beta^\theta, p^\theta, q^\theta)$ to \eqref{3.14}, depending smoothly on $\theta\in \mathbb{R}^3$, such that the following hold:

{\rm(i)} When $\theta=0\in\mathbb{R}^3$, there has $(h^0, u^0, \alpha^0, \beta^0, p^0, q^0)=(h, u, \alpha, \beta, p, q)$;

{\rm(ii)} At the point $(X_0, Y_0)$, when $\theta=0$ there has
\begin{align}\label{4.14}
{\rm rank}\ D_\theta(\alpha^\theta, \alpha_{X}^\theta, \alpha_{XX}^{\theta})=3,
\end{align}
where
\begin{align*}
D_\theta(\alpha^\theta, \alpha_{X}^\theta, \alpha_{XX}^{\theta})=
\left(
\begin{array}{ccc}
\pa_{\theta_1} \alpha  & \pa_{\theta_2} \alpha  & \pa_{\theta_3} \alpha \\
\pa_{\theta_1} \alpha_X  & \pa_{\theta_2} \alpha_X  & \pa_{\theta_3} \alpha_X \\
\pa_{\theta_1} \alpha_{XX}  & \pa_{\theta_2} \alpha_{XX}  & \pa_{\theta_3} \alpha_{XX}
\end{array}
\right).
\end{align*}

\noindent {\rm(2)} If $(\alpha, \beta, \alpha_{X})(X_0, Y_0)=(-\pi, -\pi, 0)$, then there exists a 3-parameter family of smooth solutions $(h^\theta, u^\theta, \alpha^\theta, \beta^\theta, p^\theta, q^\theta)$ to \eqref{3.14}
that satisfy {\rm (i)} and {\rm (ii)} listed above, with \eqref{4.14} replaced by the relation
\begin{align}\label{4.15}
{\rm rank}\ D_\theta(\alpha^\theta, \beta^\theta, \alpha_{X}^{\theta})=3,
\end{align}
where
\begin{align*}
D_\theta(\alpha^\theta, \beta^\theta, \alpha_{X}^{\theta})=
\left(
\begin{array}{ccc}
\pa_{\theta_1} \alpha  & \pa_{\theta_2} \alpha  & \pa_{\theta_3} \alpha \\
\pa_{\theta_1} \beta  & \pa_{\theta_2} \beta  & \pa_{\theta_3} \beta \\
\pa_{\theta_1} \alpha_{X}  & \pa_{\theta_2} \alpha_{X}  & \pa_{\theta_3} \alpha_{X}
\end{array}
\right).
\end{align*}
\end{lem}
\begin{proof}
The proof is analogous to that of Lemma 4 in Bressan and Chen \cite{BC1}.

Step 1. By smoothly approximating the boundary data, we may assume without loss of generality that $(h, u, \alpha, \beta, p, q)(X,Y)$ is a $C^\infty$ solution defined on $\Omega_\delta$ for system \eqref{3.14} with $(\overline{h}, \overline{u}, \overline{\alpha}, \overline{\beta}, \overline{p}, \overline{q})(s)$ on \(\Gamma_\kappa: \{X+Y=\kappa\}\). We first compute some values at the point $(X_0, Y_0)$. Along $\Gamma_\kappa$, we see that for $s\in[X_0-\delta, X_0+\delta]$
\begin{align}\label{4.16}
\alpha_X-\alpha_Y=\overline{\alpha}'(s),\qquad \beta_X-\beta_Y=\overline{\beta}'(s),\qquad q_X-q_Y=\overline{q}'(s),
\end{align}
from which and the equations in \eqref{3.14}, one has
\begin{align}\label{4.17}
\begin{split}
\alpha_X(X_0,Y_0)=&\overline{\alpha}'(X_0) - \fr{\gamma+1}{4}\overline{h}(X_0)\overline{q}(X_0) \sin\fr{\overline{\alpha}(X_0)}{2}\sin\fr{\overline{\alpha}(X_0)-\overline{\beta}(X_0)}{2}, \\
\beta_Y(X_0,Y_0)=&-\overline{\beta}'(X_0) - \fr{\gamma+1}{4}\overline{h}(X_0)\overline{p}(X_0) \sin\fr{\overline{\beta}(X_0)}{2}\sin\fr{\overline{\beta}(X_0)-\overline{\alpha}(X_0)}{2}, \\
q_Y(X_0,Y_0)=&-\overline{q}'(X_0) -\fr{\gamma+1}{8}\overline{h}(X_0)\overline{p}(X_0) \overline{q}(X_0) \cos\fr{\overline{\beta}(X_0)}{2}\sin\fr{\overline{\alpha}(X_0)-\overline{\beta}(X_0)}{2}.
\end{split}
\end{align}
Differentiating the equation for $\alpha$ in \eqref{4.16} once again, we arrive at
\begin{align}
\overline{\alpha}''(s)=[\alpha_{XX} +\alpha_{YY}-2\alpha_{XY}](s,\kappa-s), \notag
\end{align}
which leads to
\begin{align}\label{4.18}
\alpha_{XX}(X_0,Y_0)=\overline{\alpha}''(X_0)+2\alpha_{XY}(X_0,Y_0)-\alpha_{YY}(X_0,Y_0).
\end{align}
One differentiates the equation for $\alpha$ in \eqref{4.14} with respect to $X$ and $Y$, respectively, to obtain by rearrangement
\begin{align}\label{4.19}
\alpha_{XY}=& -\fr{\gamma+1}{4}\bigg\{q\sin\fr{\alpha}{2}\sin\fr{\alpha-\beta}{2}h_X +h\sin\fr{\alpha}{2}\sin\fr{\alpha-\beta}{2}q_X \notag \\ &+\fr{1}{2}hq\cos\fr{\alpha}{2}\sin\fr{\alpha-\beta}{2}\alpha_X +\fr{1}{2}hq\sin\fr{\alpha}{2}\cos\fr{\alpha-\beta}{2}(\alpha_X-\beta_X)\bigg\} \notag
\\
=&\fr{(\gamma+1)(3-\gamma)}{32}h^2pq\sin^2\fr{\alpha}{2}\sin\fr{\alpha-\beta}{2} -\fr{\gamma+1}{8}hq\sin\fr{2\alpha-\beta}{2}\alpha_X,
\end{align}
and
\begin{align}\label{4.20}
\alpha_{YY}=&-\fr{\gamma+1}{4}\bigg\{q\sin\fr{\alpha}{2}\sin\fr{\alpha-\beta}{2}h_Y +h\sin\fr{\alpha}{2}\sin\fr{\alpha-\beta}{2}q_Y \notag \\ &+\fr{1}{2}hq\cos\fr{\alpha}{2}\sin\fr{\alpha-\beta}{2}\alpha_Y +\fr{1}{2}hq\sin\fr{\alpha}{2}\cos\fr{\alpha-\beta}{2}(\alpha_Y-\beta_Y)\bigg\} \notag
\\
=& -\fr{\gamma^2-1}{16}h^2q^2\sin\fr{\beta}{2}\sin\fr{\alpha}{2}\sin\fr{\alpha-\beta}{2} + \fr{(\gamma+1)^2}{32} h^2q^2\sin\fr{2\alpha-\beta}{2}\sin\fr{\alpha}{2}\sin\fr{\alpha-\beta}{2}   \notag \\
&-\fr{\gamma+1}{4}h\sin\fr{\alpha}{2}\sin\fr{\alpha-\beta}{2}q_Y +\fr{\gamma+1}{8}hq\sin\fr{\alpha}{2}\cos\fr{\alpha-\beta}{2}\beta_Y.
\end{align}
By substituting \eqref{4.17} into \eqref{4.19} and \eqref{4.20}, we deduce
\begin{align}\label{4.21}
\alpha_{XY}(X_0,Y_0)=&\fr{(\gamma+1)(3-\gamma)}{32}(\overline{h}(X_0))^2\overline{p}(X_0) \overline{q}(X_0)\sin^2\fr{\overline{\alpha}(X_0)}{2}\sin\fr{\overline{\alpha}(X_0) -\overline{\beta}(X_0)}{2} \notag \\ &-\fr{\gamma+1}{8}\overline{h}(X_0)\overline{q}(X_0) \sin\fr{2\overline{\alpha}(X_0)-\overline{\beta}(X_0)}{2}\alpha_X(X_0,Y_0) \notag \\
=:&\ \alpha_{12},
\end{align}
and
\begin{align}\label{4.22}
&\alpha_{YY}(X_0,Y_0) \notag \\
=& -\fr{\gamma^2-1}{16}(\overline{h}(X_0))^2(\overline{q}(X_0))^2  \sin\fr{\overline{\beta}(X_0)}{2} \sin\fr{\overline{\alpha}(X_0)}{2}\sin\fr{\overline{\alpha}(X_0)-\overline{\beta}(X_0)}{2} \notag \\
&+\fr{(\gamma+1)^2}{32}(\overline{h}(X_0))^2(\overline{q}(X_0))^2 \sin\fr{2\overline{\alpha}(X_0)-\overline{\beta}(X_0)}{2} \sin\fr{\overline{\alpha}(X_0)}{2}\sin\fr{\overline{\alpha}(X_0)-\overline{\beta}(X_0)}{2}
\notag \\
&-\fr{\gamma+1}{4}\overline{h}(X_0) \sin\fr{\overline{\alpha}(X_0)}{2}\sin\fr{\overline{\alpha}(X_0) -\overline{\beta}(X_0)}{2}q_Y(X_0,Y_0) \notag \\ &+\fr{\gamma+1}{8}\overline{h}(X_0)\overline{q}(X_0) \sin\fr{\overline{\alpha}(X_0)}{2}\cos\fr{\overline{\alpha}(X_0) -\overline{\beta}(X_0)}{2}\beta_Y(X_0,Y_0) \notag \\
=:&\ \alpha_{22}.
\end{align}
One inserts \eqref{4.21} and \eqref{4.22} into \eqref{4.18} to gain
\begin{align}\label{4.23}
\alpha_{XX}(X_0,Y_0)=\overline{\alpha}''(X_0)+2\alpha_{12}-\alpha_{22}.
\end{align}

Step 2. Following Bressan and Chen \cite{BC1}, we construct the perturbations of the boundary data as follows:
\begin{align}\label{4.24}
\left\{
\begin{array}{l}
\overline{\alpha}^{\theta}(s)=\overline{\alpha}(s)+\sum_{i=1}^3\theta_iW_i(s), \\
\overline{\beta}^{\theta}(s)=\overline{\beta}(s)+\sum_{i=1}^3\theta_iZ_i(s),
\end{array}
\right.
\qquad
\left\{
\begin{array}{l}
\overline{p}^{\theta}(s)=\overline{p}(s)+\sum_{i=1}^3\theta_iP_i(s), \\
\overline{q}^{\theta}(s)=\overline{q}(s)+\sum_{i=1}^3\theta_iQ_i(s),
\end{array}
\right.
\end{align}
where $W_i, Z_i, P_i, Q_i$ are suitable $C_{c}^\infty$ functions. At the point $(X_0, Y_0)$, we set
\begin{align}\label{4.25}
\overline{h}^{\theta}(X_0)=h(X_0, Y_0)+\sum_{i=1}^3\theta_iH_i,\qquad \overline{u}^{\theta}(X_0)=u(X_0, Y_0)+\sum_{i=1}^3\theta_iU_i,
\end{align}
for some constants $H_i$ and $U_i$. Let $\tilde{\theta}>0$ be sufficiently small so that for $\|\theta\|\leq \tilde{\theta}$
\begin{align}\label{4.26}
\fr{1}{2}h(X_0, Y_0)\leq \overline{h}^{\theta}(X_0)\leq 2h(X_0, Y_0).
\end{align}
Then we solve the ODE problems analogous to \eqref{4.5}, \eqref{4.5a} to obtain the boundary data for $\overline{h}(s)$ and $\overline{u}(s)$ on $\Gamma_\kappa$ near the point $(X_0,Y_0)$. Moreover, there holds by \eqref{4.26}
\begin{align}\label{4.27}
\fr{1}{4}h(X_0, Y_0)\leq \overline{h}(s)\leq 4h(X_0, Y_0),\quad \forall\ s\in[X_0-\delta', X_0+\delta'],
\end{align}
where $\delta'>0$ is a small constant. According to Lemma \ref{lem3}, for each $\theta\in\mathbb{R}^3$ satisfying $\|\theta\|\leq \tilde{\theta}$, the semilinear hyperbolic system \eqref{3.14} with the boundary data $(\overline{h}(s), \overline{u}(s), \overline{\alpha}(s), \overline{\beta}(s), \overline{p}(s), \overline{q}(s))$ on $\Gamma_\kappa$ admits a unique smooth solution $(h^\theta, u^\theta, \alpha^\theta, \beta^\theta, p^\theta, q^\theta)$ on the domain $\Omega_{\delta''}$, where $\delta''\leq \delta'$ is a small positive constant. In addition, the solution fulfils
\begin{align}\label{4.28}
\fr{1}{8}h(X_0, Y_0)\leq h^\theta(X,Y)\leq 8h(X_0, Y_0),\quad \forall\ (X,Y)\in\Omega_{\delta''},\ \ \forall\ \|\theta\|\leq \tilde{\theta}.
\end{align}

It is observed that the functions $W_i, Z_i, P_i, Q_i$ can be chosen arbitrarily. Consequently,
the following value can be assigned arbitrarily
\begin{align}\label{4.29}
\fr{{\rm d}}{{\rm d}\theta}\fr{{\rm d}^k}{{\rm d}s^k}\overline{\alpha}^\theta\bigg|_{s=X_0, \theta=0},\  \fr{{\rm d}}{{\rm d}\theta}\fr{{\rm d}^k}{{\rm d}s^k}\overline{\beta}^\theta\bigg|_{s=X_0, \theta=0},\ \fr{{\rm d}}{{\rm d}\theta}\fr{{\rm d}^k}{{\rm d}s^k}\overline{p}^\theta\bigg|_{s=X_0, \theta=0},\ \fr{{\rm d}}{{\rm d}\theta}\fr{{\rm d}^k}{{\rm d}s^k}\overline{q}^\theta\bigg|_{s=X_0, \theta=0},
\end{align}
for $k=0,1,\dots$. Furthermore, the quantities
$$
\fr{{\rm d}}{{\rm d}\theta}\fr{{\rm d}^k}{{\rm d}s^k}\overline{h}^\theta\bigg|_{s=X_0, \theta=0},\qquad \fr{{\rm d}}{{\rm d}\theta}\fr{{\rm d}^k}{{\rm d}s^k}\overline{u}^\theta\bigg|_{s=X_0, \theta=0}
$$
can also be chosen arbitrarily, subject to the corresponding relations for the differential equation \eqref{4.5} and its higher-order equations.

Step 3. We now choose the perturbations $(\overline{h}^\theta, \overline{u}^\theta, \overline{\alpha}^\theta, \overline{\beta}^\theta, \overline{p}^\theta, \overline{q}^\theta)$ so that at $s=X_0,\ \theta=0$
\begin{align*}
D_\theta
\left(
\begin{array}{c}
\overline{h}^\theta \\
\overline{u}^\theta \\
\overline{\alpha}^\theta \\
\overline{\beta}^\theta \\
{\overline{\beta}^{\theta}}' \\
{\overline{\alpha}^{\theta}}' \\
{\overline{\alpha}^{\theta}}'' \\
\overline{p}^\theta \\
\overline{q}^\theta \\
{\overline{q}^{\theta}}'
\end{array}
\right)=
\left(
\begin{array}{ccc}
0 & 0 & 0 \\
0 & 0 & 0 \\
1 & 0 & 0 \\
0 & 0 & 0 \\
0 & 0 & 0 \\
0 & 1 & 0 \\
0 & 0 & 1 \\
0 & 0 & 0 \\
0 & 0 & 0 \\
0 & 0 & 0
\end{array}
\right).
\end{align*}
Then we combine \eqref{4.17} and \eqref{4.23} to obtain at the point $(X_0, Y_0)$ and $\theta=0$
\begin{align*}
D_\theta
\left(
\begin{array}{c}
\alpha \\
\alpha_X \\
\alpha_{XX}
\end{array}
\right)
=\left(
\begin{array}{ccc}
1 & 0 & 0 \\
* & 1 & 0 \\
* & * & 1
\end{array}
\right),
\end{align*}
which is the desired \eqref{4.14}.

To arrive at \eqref{4.15}, we can choose the perturbations $(\overline{h}^\theta, \overline{u}^\theta, \overline{\alpha}^\theta, \overline{\beta}^\theta, \overline{p}^\theta, \overline{q}^\theta)$ so that at $s=X_0, \theta=0$
\begin{align*}
D_\theta
\left(
\begin{array}{c}
\overline{h}^\theta \\
\overline{u}^\theta \\
\overline{\alpha}^\theta \\
\overline{\beta}^\theta \\
{\overline{\alpha}^{\theta}}' \\
\overline{p}^\theta \\
\overline{q}^\theta
\end{array}
\right)=
\left(
\begin{array}{ccc}
0 & 0 & 0 \\
0 & 0 & 0 \\
1 & 0 & 0 \\
0 & 1 & 0 \\
0 & 0 & 1 \\
0 & 0 & 0 \\
0 & 0 & 0
\end{array}
\right).
\end{align*}
Then at the point $(X_0, Y_0)$ and $\theta=0$, we have by \eqref{4.17}
\begin{align*}
D_\theta
\left(
\begin{array}{c}
\alpha \\
\beta \\
\alpha_{X}
\end{array}
\right)
=\left(
\begin{array}{ccc}
1 & 0 & 0 \\
0 & 1 & 0 \\
* & * & 1
\end{array}
\right).
\end{align*}
This yields \eqref{4.15} and then the proof of the lemma is complete.
\end{proof}

\subsection{The generic structure of singular points}\label{S43}

Based on Lemma \ref{lem4}, we can show the generic structure of smooth solutions to the semilinear hyperbolic system \eqref{3.14}.

\begin{lem}\label{lem5}
Let $L$ be a large positive number. Consider the compact domain of the form
\begin{align}\label{4.30}
\Omega_L=\Omega_{t^*}\cap\{(X,Y):\ |X|\leq L,\quad |Y|\leq L\}.
\end{align}
Let \(\mathcal{S}\) denote the family of all \(C^2\) solutions to system \eqref{3.14} satisfying \(p>0, q>0\) and $0<h<\infty$ for every \((X,Y)\in \Omega_{t^*}\). Furthermore, let \(\mathcal{S}'\subset\mathcal{S}\) stand for the subfamily consisting of all solutions \((h, u, \alpha, \beta, p, q)\) such that none of the quantities listed below is attained for \((X,Y)\in\Omega_L\):
\begin{align}\label{4.31}
\left\{
\begin{array}{r}
(\alpha, \alpha_X, \alpha_{XX})=(-\pi, 0, 0), \\
(\beta, \beta_Y, \beta_{YY})=(-\pi, 0, 0),
\end{array}
\right.
\end{align}
\begin{align}\label{4.32}
\left\{
\begin{array}{r}
(\alpha, \beta, \alpha_{X})=(-\pi, -\pi, 0), \\
(\alpha, \beta, \beta_{Y})=(-\pi, -\pi, 0).
\end{array}
\right.
\end{align}
Then \(\mathcal{S}'\) is a relatively open and dense subset of \(\mathcal{S}\) with respect to the topology induced by \(C^2(\Omega_L)\).
\end{lem}
\begin{proof}
This lemma can be proved by means of Lemma \ref{lem4} together with Thom's transversality theorem \cite{Thom, BC1}. The proof is identical to that of Lemma 5 in Bressan and Chen \cite{BC1}, so we omit it here.
\end{proof}

By virtue of Lemma \ref{lem5}, we can perform the proof of Theorem 1 in Bressan and Chen \cite{BC1} to obtain the following theorem.
\begin{thm}\label{thm1}
Let $\bar{\tau}_0>0$ be a constant. There exists an open dense set of initial data
\begin{align}\label{4.33}
\mathcal{D}\subset\Big(C^3(\mathbb{R})\cap W^{1,1}(\mathbb{R})\Big)^2,
\end{align}
with
\begin{align}\label{4.34}
\begin{array}{c}
\dps\fr{1}{\overline{c}}\leq \tau_0(x)\leq \overline{c}, \\[6pt]
\dps\min\bigg\{\inf\Big(u_{0}'(x)-(\tau_0(x))^{-\fr{\gamma+1}{2}}\tau_{0}'(x)\Big),\ \inf\Big(u_{0}'(x)+(\tau_0(x))^{-\fr{\gamma+1}{2}}\tau_{0}'(x)\Big)\bigg\}<0,
\end{array}
\end{align}
for some positive constant $\overline{c}$. For initial data $(\tau_0(x), u_0(x))$ satisfying $(\tau_0(x)-\bar{\tau}_0, u_0(x))\in \mathcal{D}$ and \eqref{4.34},
the solution $(\tau(x,t), u(x,t))$ of \eqref{1.1}-\eqref{1.3}, defined on $\mathbb{R}\times[0,t^*]$, is twice continuously differentiable in the complement of finitely many points on the line $t=t^*$. Here $t^*$ is the first time of singularity formation for the solution $(\tau, u)$.
Moreover, for the corresponding solution $(h, u, \alpha, \beta, p, q)$ of \eqref{3.14} with boundary data \eqref{3.19} on the domain $\Omega_{t^*}$, none of the values in \eqref{4.31} and \eqref{4.32} can be attained at any point $(X,Y)$ on the boundary $t(X,Y)=t^*$.
\end{thm}

\section{Structural stability for generic singularities}\label{S5}

In this section, we analyze the asymptotic behavior of the solution in a neighborhood of a singular point for the initial data in $\mathcal{D}$.

\begin{figure}[htbp]
\begin{center}
\qquad
\begin{minipage}[t]{0.45\textwidth}
\includegraphics[scale=0.48]{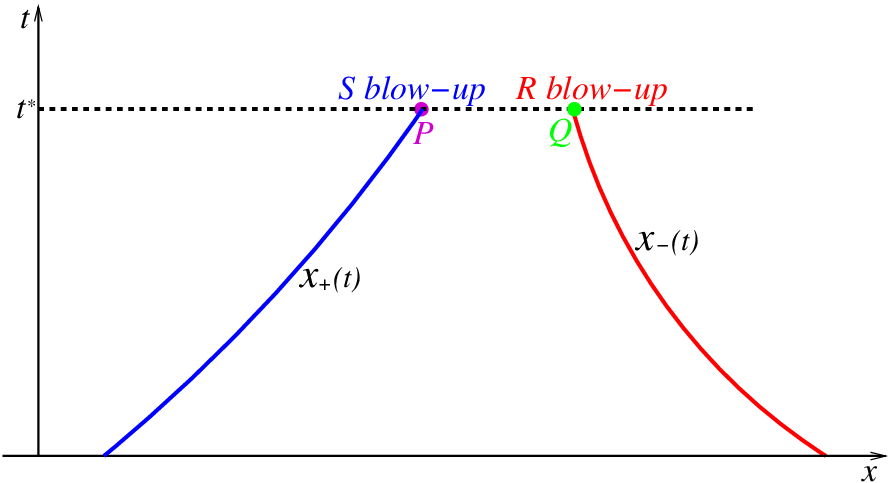}
\end{minipage}\qquad
\begin{minipage}[t]{0.45\textwidth}
\includegraphics[scale=0.42]{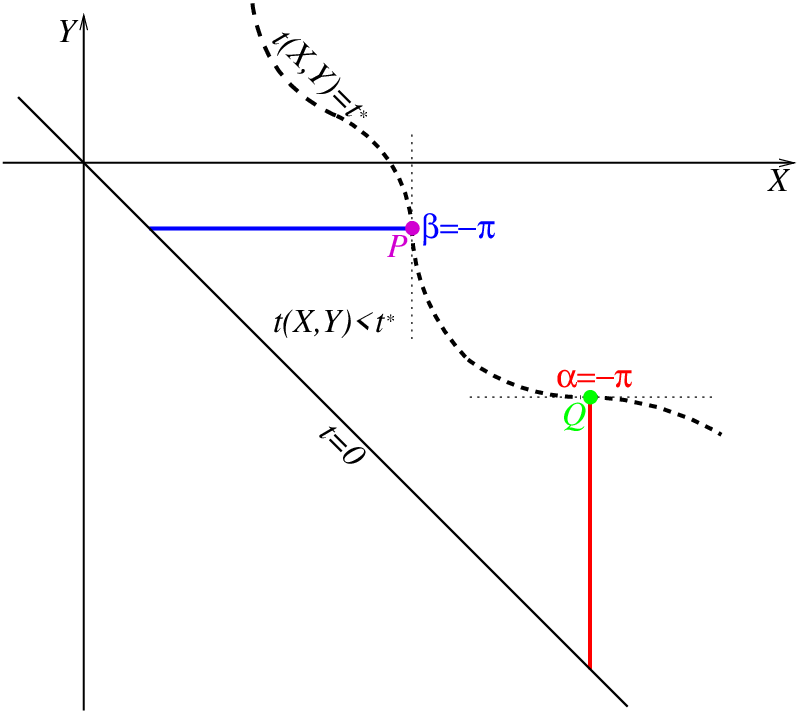}
\end{minipage}
\caption{Schematic diagram of Type-1 singular points. Left: singular points in the $(x,t)$ plane; Right: the corresponding singular points in the $(X,Y)$ plane.}
\label{fig4}
\end{center}
\end{figure}

Let $(\tau_0(x)-\bar{\tau}_0, u_0(x))\in \mathcal{D}$ satisfy \eqref{4.34}. Then the solution $(\tau, u, \alpha, \beta, p, q, x, t)$ of \eqref{3.14a} and \eqref{3.16} is smooth on the domain $\Omega_{t^*}$. Since $\tau$ admits positive lower and upper bounds up to the boundary $t(X,Y)=t^*$, it follows from Lemma \ref{lem3} that the system \eqref{3.14a} possesses a local smooth solution near the boundary $t(X,Y)=t^*$. Moreover, in the solvability domain, $p, q$ remain positive and bounded and $\tau$ has positive lower and upper bounds. For the points where
either $\alpha=-\pi$ or $\beta=-\pi$ on the boundary $t(X,Y)=t^*$, we know by Theorem \ref{thm1} that
\begin{align*}
\begin{split}
\alpha=-\pi\quad {\rm and}\quad \alpha_X=0 \quad &\Longrightarrow \quad \alpha_{XX}\neq 0,\\
\beta=-\pi\quad {\rm and}\quad \beta_Y=0 \quad &\Longrightarrow \quad \beta_{YY}\neq 0,\\
\alpha=-\pi\quad {\rm and}\quad \beta=-\pi \quad &\Longrightarrow \quad \alpha_{X}\neq 0 \quad {\rm and}\quad \beta_{Y}\neq 0.\\
\end{split}
\end{align*}
Thus we distinguish three types of singular points $Q=(X_0,Y_0)$ on the boundary $t(X,Y)=t^*$:
\begin{itemize}

\item[] Type-1: Points where $\alpha=-\pi$, $\alpha_X= 0$, but $\beta\neq-\pi$, $\alpha_{XX}\neq0$ (or else, where $\beta=-\pi$, $\beta_Y= 0$, but $\alpha\neq-\pi$, $\beta_{YY}\neq0$).

\item[] Type-2: Points where $\alpha=-\pi$, but $\alpha_X\neq 0$, $\beta\neq-\pi$ (or else, where $\beta=-\pi$ but $\beta_Y\neq 0$, $\alpha\neq-\pi$).

\item[] Type-3: Points where $\alpha=-\pi$, $\beta= -\pi$, but $\alpha_X\neq 0$, $\beta_Y\neq0$.

\end{itemize}
See Fig. \ref{fig4} for the schematic diagram of Type-1 singular points.

%+++++++++++++++++++++++++++++++++++++++++++++++++++++++++++++++++++++
%+++++++++++++++++++++++++++++++++++++++++++++++++++++++++++++++++++++

Then we have
\begin{thm}\label{thm2}
Consider generic initial data $(\tau_0(x), u_0(x))$ satisfying $(\tau_0(x)-\bar{\tau}_0, u_0(x))\in \mathcal{D}$ and \eqref{4.34}.
Let $(\tau, u, \alpha, \beta, p, q, x, t)$ be the smooth solution of the semilinear hyperbolic system \eqref{3.14a} and \eqref{3.16} on the domain $\Omega_{t^*}$, and $(\tau(x,t), u(x,t))$ be the solution to the Cauchy problem \eqref{1.1}-\eqref{1.3} on the domain $\mathbb{R}\times[0, t^*]$. Let $(X_0, Y_0)$ be a singular point on the boundary $t(X,Y)=t^*$ where $\alpha=-\pi$ and set $(x^*,t^*)=(x(X_0, Y_0), t(X_0, Y_0))$. Then $(X_0, Y_0)$ must be a singular point of Type-1, that is, singularities of the Type-2 and Type-3 cannot occur on $t(X,Y)=t^*$. Moreover, there exists a constant $\zeta\neq0$ such that for $t\leq t^*$
\begin{align}\label{5.1a}
\begin{split}
\tau(x,t)=&\tau(x^*,t^*) +\zeta\cdot [c(\tau(x^*,t^*))(t-t^*)+(x-x^*)]^{\fr{1}{3}} \\ &+O(1)\Big(|t-t^*|+|x-x^*|\Big)^{\fr{4}{9}},
\end{split}
\end{align}
and
\begin{align}\label{5.1aa}
\begin{split}
u(x,t)=&u(x^*,t^*) +c(\tau(x^*, t^*))\zeta\cdot [c(\tau(x^*,t^*))(t-t^*)+(x-x^*)]^{\fr{1}{3}} \\ &+O(1)\Big(|t-t^*|+|x-x^*|\Big)^{\fr{4}{9}}.
\end{split}
\end{align}
In addition, the Riemann invariant $\varpi_+(x,t)$ is smooth at the singular point $(x^*, t^*)$.
\end{thm}

Based on Theorem \ref{thm2}, we can give the following definition.
\begin{defn}\label{def1}
A solution $(\tau(x,t), u(x,t))$ of problem \eqref{1.1}-\eqref{1.3} is said to have only generic singularities for $t\in[0,t^*]$ if it admits a representation of the form
\begin{align}\label{5.1b}
\begin{split}
{\rm Graph}(\tau, u)=\Big\{\big(t(X,Y), x(X,Y),\tau(X,Y), u(X,Y)\big);\ (X,Y)\in \Omega_{t^*}\Big\},
\end{split}
\end{align}
where (i) the functions $(\tau, \alpha, \beta, p, q, x, t)(X,Y)$ and $u(X,Y)$ are smooth; and (ii) the following generic conditions hold:
\begin{align}\label{5.1c}
\begin{split}
{\rm(G1)}&\qquad \alpha=-\pi,\ \ \alpha_X=0\ \Longrightarrow \ \alpha_Y\neq0,\ \ \alpha_{XX}\neq0,\\
{\rm(G2)}&\qquad \beta=-\pi,\ \ \ \beta_Y=0\   \Longrightarrow \ \ \beta_X\neq0,\ \ \beta_{YY}\neq0.
\end{split}
\end{align}
\end{defn}

\subsection{Singular points of Type-1}\label{S51}

Let $Q_1=(X_0,Y_0)$ on the boundary $t(X,Y)=t^*$ be a point of Type-1, where
\begin{align}\label{5.16}
\alpha(Q_1)=-\pi,\qquad \alpha_X(Q_1)=0, \qquad \alpha_{XX}(Q_1)\neq0,\qquad  \beta(Q_1)\neq-\pi.
\end{align}
It suggests by \eqref{3.14a} and \eqref{5.16} that
\begin{align}\label{5.17}
\begin{split}
\alpha_X(Q_1)=& 0,\quad \alpha_{XX}(Q_1)\neq0, \\[4pt]
\alpha_Y(Q_1)=&\Big(\fr{c'}{2c^2}q\sin\fr{\alpha}{2}\sin\fr{\alpha-\beta}{2}\Big)(Q_1)<0, \\[4pt]
\tau_X(Q_1)=&\Big(\fr{1}{2c}p\sin\fr{\alpha}{2}\Big)(Q_1)=:\zeta_{\tau, 1}<0,\\[4pt]
\tau_Y(Q_1)=&\Big(\fr{1}{2c}q\sin\fr{\beta}{2}\Big)(Q_1)=:\zeta_{\tau, 2}\neq0, \\[4pt]
u_X(Q_1)=&\Big(\fr{1}{2}p\sin\fr{\alpha}{2}\Big)(Q_1)=c(\tau(Q_1))\zeta_{\tau, 1},\\[4pt]
u_Y(Q_1)=&\Big(-\fr{1}{2}q\sin\fr{\beta}{2}\Big)(Q_1)=-c(\tau(Q_1))\zeta_{\tau, 2}.
\end{split}
\end{align}
Moreover, we have by \eqref{3.16} and \eqref{5.16}
\begin{align}\label{5.18}
\begin{split}
t_X(Q_1)=& \Big(\fr{1}{2c}p\cos\fr{\alpha}{2}\Big)(Q_1)=0 \\[4pt]
t_Y(Q_1)=&\Big(\fr{1}{2c}q\cos\fr{\beta}{2}\Big)(Q_1):=\zeta_{t,1}\neq0, \\[4pt]
x_X(Q_1)=&\Big(\fr{1}{2}p\cos\fr{\alpha}{2}\Big)(Q_1)=0,\\[4pt]
x_Y(Q_1)=&\Big(-\fr{1}{2}q\cos\fr{\beta}{2}\Big)(Q_1)=-c(\tau(Q_1))\zeta_{t,1}.
\end{split}
\end{align}

To obtain the information on the higher-order derivatives of $t$ at point $Q_1$, we differentiate the equation for $t$ in \eqref{3.16} with respect to $X$ and $Y$ to attain
\begin{align}\label{5.18a}
\begin{split}
t_{XX}=&\Big(\fr{1}{2c}p\cos\fr{\alpha}{2}\Big)_X =\Big(\fr{1}{2c}p\Big)_X\cos\fr{\alpha}{2} -\fr{1}{4c}p\sin\fr{\alpha}{2}\alpha_X, \\[4pt]
t_{XY}=&\Big(\fr{1}{2c}p\cos\fr{\alpha}{2}\Big)_Y =\Big(\fr{1}{2c}p\Big)_Y\cos\fr{\alpha}{2} -\fr{1}{4c}p\sin\fr{\alpha}{2}\alpha_Y,
\end{split}
\end{align}
from which and \eqref{5.16}, one sees that
\begin{align}\label{5.19}
\begin{split}
t_{XX}(Q_1)=& 0, \\[4pt]
t_{XY}(Q_1)=&\Big( -\fr{1}{4c}p\sin\fr{\alpha}{2}\alpha_Y\Big)(Q_1):=\zeta_{t,2}<0.
\end{split}
\end{align}
We further differentiate the first equation in \eqref{5.18a} with respect to $X$ to arrive at
\begin{align}\label{5.20}
t_{XXX}=\Big(\fr{1}{2c}p\Big)_{XX}\cos\fr{\alpha}{2} -\Big(\fr{1}{2c}p\Big)_X\sin\fr{\alpha}{2}\alpha_X -\fr{1}{4c}p\sin\fr{\alpha}{2}\alpha_{XX},
\end{align}
which together with \eqref{5.16} and \eqref{5.17} yields
\begin{align}\label{5.21}
t_{XXX}(Q_1)= \Big(-\fr{1}{4c}p\sin\fr{\alpha}{2}\alpha_{XX}\Big)(Q_1)=:\zeta_{t,3}\neq0.
\end{align}
Analogously, there hold for the variable $x$
\begin{align}\label{5.22a}
\begin{split}
x_{XX}(Q_1)=& 0, \\[4pt]
x_{XY}(Q_1)=&\Big( -\fr{1}{4}p\sin\fr{\alpha}{2}\alpha_Y\Big)(Q_1)=c(\tau(Q_1))\zeta_{t,2},
\end{split}
\end{align}
and then
\begin{align}\label{5.22}
x_{XXX}(Q_1)= \Big(-\fr{1}{4}p\sin\fr{\alpha}{2}\alpha_{XX}\Big)(Q_1)=c(\tau(Q_1))\zeta_{t,3}.
\end{align}

Now in view of the information of $\tau$ and $u$ in \eqref{5.17}, we obtain the local Taylor approximations
\begin{align}\label{5.23}
\begin{split}
\tau(X,Y)=&\tau^* +\zeta_{\tau,1}(X-X_0) +\zeta_{\tau,2}(Y-Y_0) \\
&+O(1)\Big(|X-X_0|^2 +|X-X_0||Y-Y_0|+|Y-Y_0|^2\Big),
\end{split}
\end{align}
and
\begin{align}\label{5.23a}
\begin{split}
u(X,Y)=&u^* +c^*\zeta_{\tau,1}(X-X_0) -c^*\zeta_{\tau,2}(Y-Y_0) \\
&+O(1)\Big(|X-X_0|^2 +|X-X_0||Y-Y_0|+|Y-Y_0|^2\Big).
\end{split}
\end{align}
Moreover, combining \eqref{5.18}, \eqref{5.19} and \eqref{5.22}, we have the Taylor approximations for $t$ and $x$
\begin{align}\label{5.24}
\begin{split}
t(X,Y)=&t^* +\zeta_{t,1}(Y-Y_0) +\zeta_{t,2}(X-X_0)(Y-Y_0) +\fr{1}{6}\zeta_{t,3}(X-X_0)^3\\
&+O(1)\Big(|Y-Y_0|^2 +|X-X_0|^2|Y-Y_0|+|X-X_0|^4\Big),
\end{split}
\end{align}
and
\begin{align}\label{5.24a}
\begin{split}
x(X,Y)=&x^* -c^*\zeta_{t,1}(Y-Y_0) +c^*\zeta_{t,2}(X-X_0)(Y-Y_0)+\fr{c^*}{6}\zeta_{t,3}(X-X_0)^3 \\
&+O(1)\Big(|Y-Y_0|^2 +|X-X_0|^2|Y-Y_0|+|X-X_0|^4\Big).
\end{split}
\end{align}
Here and below, $\tau^*=\tau(X_0, Y_0)$, $u^*=u(X_0, Y_0)$, $x^*=x(X_0, Y_0)$ and $c^*=c(\tau^*)$.
It deduces by \eqref{5.24} and \eqref{5.24a} that
\begin{align}\label{5.26}
\begin{split}
&c^*(t-t^*)-(x-x^*)=2c^*\zeta_{t,1}(Y-Y_0) \\
&\qquad +O(1)\Big(|Y-Y_0|^2 +|X-X_0|^2|Y-Y_0|+|X-X_0|^4\Big),
\end{split}
\end{align}
and
\begin{align}\label{5.27}
\begin{split}
&c^*(t-t^*)+(x-x^*)=\fr{1}{3}c^*\zeta_{t,3}(X-X_0)^3\\
&\qquad +O(1)\Big(|Y-Y_0|^2 +|X-X_0||Y-Y_0|+|X-X_0|^4\Big).
\end{split}
\end{align}
From \eqref{5.26} and \eqref{5.27}, one achieves
\begin{align}\label{5.28}
\begin{split}
Y-Y_0=&\fr{1}{2c^*\zeta_{t,1}}[c^*(t-t^*)-(x-x^*)] \\
&+O(1)\Big(|Y-Y_0|^2 +|X-X_0||Y-Y_0|+|X-X_0|^4\Big),
\end{split}
\end{align}
and
\begin{align}\label{5.29}
\begin{split}
X-X_0=&\Big(\fr{3}{c^*\zeta_{t,3}}\Big)^{\fr{1}{3}}[c^*(t-t^*)+(x-x^*)]^{\fr{1}{3}} \\
&+O(1)\Big(|Y-Y_0|^{\fr{2}{3}} +|X-X_0|^{\fr{1}{3}}|Y-Y_0|^{\fr{1}{3}}+|X-X_0|^{\fr{4}{3}}\Big).
\end{split}
\end{align}
Then we substitute \eqref{5.28} and \eqref{5.29} into \eqref{5.23} to finally obtain
\begin{align}\label{5.30}
\begin{split}
\tau(x,t)=&\tau^* +\zeta_{\tau,1}\Big(\fr{3}{c^*\zeta_{t,3}}\Big)^{\fr{1}{3}}[c^*(t-t^*)+(x-x^*)]^{\fr{1}{3}} \\
&+O(1)\Big(|t-t^*|+|x-x^*|\Big)^{\fr{4}{9}}.
\end{split}
\end{align}
This proves \eqref{5.1a}, with
$$
\zeta=\zeta_{\tau,1}\Big(\fr{3}{c^*\zeta_{t,3}}\Big)^{\fr{1}{3}}\neq0.
$$
Furthermore, it suggests by \eqref{5.23a} that
\begin{align}\label{5.30b}
\begin{split}
u(x,t)=&u^* +c^*\zeta_{\tau,1}\Big(\fr{3}{c^*\zeta_{t,1}}\Big)^{\fr{1}{3}}[c^*(t-t^*)+(x-x^*)]^{\fr{1}{3}} \\
&+O(1)\Big(|t-t^*|+|x-x^*|\Big)^{\fr{4}{9}},
\end{split}
\end{align}
which is the desired expression \eqref{5.1aa}. In addition, we combine \eqref{5.23}, \eqref{5.23a} and \eqref{5.28} to arrive at
\begin{align}\label{5.30c}
\begin{split}
u(x,t)-c^*\tau(x,t)=&(u^*-c^*\tau^*) +\fr{\zeta_{\tau,2}}{\zeta_{t,1}}[(x-x^*)-c^*(t-t^*)] \\
&+O(1)\Big(|t-t^*|+|x-x^*|\Big)^{\fr{5}{3}},
\end{split}
\end{align}
which indicates by \eqref{2.3} that
\begin{align}\label{5.30d}
\begin{split}
\pa_x\varpi_+(x^*, t^*)=&\pa_xu(x^*,t^*)-c^*\pa_x\tau(x^*,t^*)=\fr{\zeta_{\tau,2}}{\zeta_{t,1}}, \\
\pa_t\varpi_+(x^*, t^*)=&\pa_tu(x^*,t^*)-c^*\pa_t\tau(x^*,t^*)=-c^*\fr{\zeta_{\tau,2}}{\zeta_{t,1}}.
\end{split}
\end{align}
Thus the Riemann invariant $\varpi_+$ is smooth up to the singular point $(x^*, t^*)$.

\subsection{Singular points of Type-2}\label{S52}

Assume that a point $Q_2=(X_0,Y_0)$ on $t(X,Y)=t^*$ is a singular point of Type-2, that is, there holds
\begin{align}\label{5.1}
\alpha(Q_2)=-\pi,\qquad \beta(Q_2)\neq-\pi, \qquad \alpha_X(Q_2)\neq0.
\end{align}
We see by \eqref{3.14a}, \eqref{3.16} and \eqref{5.1} that
\begin{align}\label{5.2}
\begin{split}
\alpha_X(Q_2)<& 0 \\[4pt]
\alpha_Y(Q_2)=&\Big(\fr{c'}{2c^2}q\sin\fr{\alpha}{2}\sin\fr{\alpha-\beta}{2}\Big)(Q_2)<0, \\[4pt]
\tau_X(Q_2)=&\Big(\fr{1}{2c}p\sin\fr{\alpha}{2}\Big)(Q_2)=:\sigma_{\tau, 1}<0,\\[4pt]
\tau_Y(Q_2)=&\Big(\fr{1}{2c}q\sin\fr{\beta}{2}\Big)(Q_2)=:\sigma_{\tau, 2}\neq0,
\end{split}
\end{align}
and
\begin{align}\label{5.3}
\begin{split}
t_X(Q_2)=& \Big(\fr{1}{2c}p\cos\fr{\alpha}{2}\Big)(Q_2)=0 \\[4pt]
t_Y(Q_2)=&\Big(\fr{1}{2c}q\cos\fr{\beta}{2}\Big)(Q_2):=\sigma_{t,1}\neq0, \\[4pt]
x_X(Q_2)=&\Big(\fr{1}{2}p\cos\fr{\alpha}{2}\Big)(Q_2)=0,\\[4pt]
x_Y(Q_2)=&\Big(-\fr{1}{2}q\cos\fr{\beta}{2}\Big)(Q_2)=-c(\tau(Q_2))\sigma_{t,1}.
\end{split}
\end{align}
Recalling \eqref{5.18a}, one utilizes \eqref{5.1} and \eqref{5.2} to gain
\begin{align}\label{5.5}
\begin{split}
t_{XX}(Q_2)=& \Big(-\fr{1}{4c}p\sin\fr{\alpha}{2}\alpha_X\Big)(Q_2):=\sigma_{t,2}<0, \\[4pt]
t_{XY}(Q_2)=&\Big( -\fr{1}{4c}p\sin\fr{\alpha}{2}\alpha_Y\Big)(Q_2):=\sigma_{t,3}<0.
\end{split}
\end{align}
Similarly, one has
\begin{align}\label{5.6}
\begin{split}
x_{XX}(Q_2)=& \Big(-\fr{1}{4}p\sin\fr{\alpha}{2}\alpha_X\Big)(Q_1)=c(\tau(Q_1))\sigma_{t,2}, \\[4pt]
x_{XY}(Q_2)=&\Big( -\fr{1}{4}p\sin\fr{\alpha}{2}\alpha_Y\Big)(Q_1)=c(\tau(Q_1))\sigma_{t,3}.
\end{split}
\end{align}
%We combine \eqref{5.2}, \eqref{5.5} and \eqref{5.6} to obtain the local Taylor approximations for $\tau$ and $(t, x)$
%\begin{align}\label{5.7}
%\begin{split}
%\tau(X,Y)=&\tau^* +\beta_{\tau,1}(X-X_0) +\beta_{\tau,2}(Y-Y_0) \\
%&+O(1)\Big(|X-X_0|^2 +|X-X_0||Y-Y_0|+|Y-Y_0|^2\Big),
%\end{split}
%\end{align}
%and
%\begin{align}\label{5.8}
%\begin{split}
%t(X,Y)=&t^* +\beta_{t,1}(Y-Y_0) +\fr{1}{2}\beta_{t,2}(X-X_0)^2 +\beta_{t,3}(X-X_0)(Y-Y_0) \\
%&+O(1)\Big(|Y-Y_0|^2 +|X-X_0|^2|Y-Y_0|+|X-X_0|^3\Big), \\[4pt]
%x(X,Y)=&x^* -c^*\beta_{t,1}(Y-Y_0) +\fr{c^*}{2}\beta_{t,2}(X-X_0)^2 +c^*\beta_{t,3}(X-X_0)(Y-Y_0) \\
%&+O(1)\Big(|Y-Y_0|^2 +|X-X_0|^2|Y-Y_0|+|X-X_0|^3\Big).
%\end{split}
%\end{align}

Next we derive a contradiction. Thanks to \eqref{5.3} and \eqref{5.5}, we observe that the level set where $t(X,Y)=t^*$ is locally the graph of a smooth function $Y=\mathcal{F}(X)$, with $\mathcal{F}(X_0)=Y_0$. Furthermore, the function $Y=\mathcal{F}(X)$ is strictly concave near $X=X_0$ and attains a local maximum $Y_0$ at this point.
The local region below the curve $Y=\mathcal{F}(X)$, denoted by $\Omega_\mathcal{F}$, corresponds to the domain where \(t(X,Y)<t^{*}\). Hence $\alpha(x, Y)>-\pi$ in the region $\Omega_\mathcal{F}$.
Now we consider the level set where $\alpha(X,Y)=-\pi$. In view of the information of $\alpha$ in \eqref{5.2}, one knows that the level set where $\alpha(X,Y)=-\pi$ is locally the graph of a smooth function $Y=\mathcal{G}(X)$, with $\mathcal{G}(X_0)=Y_0$. Moreover, we see that
$$
\mathcal{G}'(X_0)=-\fr{\alpha_X(X_0)}{\alpha_Y(X_0)}\neq0,
$$
which indicates that the curve $Y=\mathcal{F}(X)$ intersects the curve $Y=\mathcal{G}(X)$ at the point $X_0$. Thus there exists some point $(X', \mathcal{G}(X'))\in \Omega_\mathcal{F}$ such that $\alpha(X', G(X'))=-\pi$, which contradicts the fact $\alpha(x, Y)>-\pi$ in the region $\Omega_\mathcal{F}$. Therefore, singularities of the Type-2 cannot occur on the boundary $t(X,Y)=t^*$.

\subsection{Singular points of Type-3}\label{S53}

Assume that a point $Q_3=(X_0,Y_0)$ on $t(X,Y)=t^*$ is a singular point of Type-3, that is, there holds
\begin{align}\label{5.31}
\alpha(Q_3)=-\pi,\qquad \beta(Q_3)=-\pi, \qquad \alpha_{X}(Q_3)\neq0,\qquad  \beta_Y(Q_3)\neq0.
\end{align}
For this case, we recall \eqref{3.14a}, \eqref{3.16} and \eqref{5.16} to find that
\begin{align}\label{5.32}
\begin{split}
\alpha_X(Q_3)\neq&0, \qquad \beta_{Y}(Q_3)\neq0, \\[4pt]
\alpha_Y(Q_3)=&\beta_X(Q_3)=p_Y(Q_3)= q_X(Q_3)=0, \\[4pt]
\tau_X(Q_3)=&\Big(\fr{1}{2c}p\sin\fr{\alpha}{2}\Big)(Q_3)=:\nu_{\tau, 1}<0,\\[4pt]
\tau_Y(Q_3)=&\Big(\fr{1}{2c}q\sin\fr{\beta}{2}\Big)(Q_3)=:\nu_{\tau, 2}<0.
\end{split}
\end{align}
and
\begin{align}\label{5.33}
\begin{split}
t_X(Q_3)=t_Y(Q_3)=x_X(Q_3)=x_Y(Q_3)=0.
\end{split}
\end{align}
We then differentiate the equations for $t$ in \eqref{3.16} with respect to $X$ and $Y$ to obtain the equations of the higher-order derivatives for $t$
\begin{align}\label{5.34}
\begin{split}
t_{XX}=&\Big(\fr{1}{2c}p\cos\fr{\alpha}{2}\Big)_X =\Big(\fr{1}{2c}p\Big)_X\cos\fr{\alpha}{2} -\fr{1}{4c}p\sin\fr{\alpha}{2}\alpha_X, \\[4pt]
t_{XY}=&\Big(\fr{1}{2c}p\cos\fr{\alpha}{2}\Big)_Y =\Big(\fr{1}{2c}p\Big)_Y\cos\fr{\alpha}{2} -\fr{1}{4c}p\sin\fr{\alpha}{2}\alpha_Y, \\[4pt]
t_{YY}=&\Big(\fr{1}{2c}q\cos\fr{\beta}{2}\Big)_Y = \Big(\fr{1}{2c}q\Big)_Y\cos\fr{\beta}{2} -\fr{1}{4c}q\sin\fr{\beta}{2}\beta_Y,
\end{split}
\end{align}
from which and \eqref{5.32} we have
\begin{align}\label{5.35}
\begin{split}
t_{XX}(Q_3)=& \Big(-\fr{1}{4c}p\sin\fr{\alpha}{2}\alpha_X\Big)(Q_3)=:\nu_{t,1}\neq0, \\[4pt]
t_{XY}(Q_3)=&0, \\[4pt]
t_{YY}(Q_3)=& \Big( -\fr{1}{4c}q\sin\fr{\beta}{2}\beta_Y\Big)(Q_3)=:\nu_{t,2}\neq0.
\end{split}
\end{align}
For the variable $x$, one also has
\begin{align}\label{5.36}
\begin{split}
x_{XX}(Q_3)=& \Big(-\fr{1}{4}p\sin\fr{\alpha}{2}\alpha_X\Big)(Q_3)=c(\tau(Q_3))\nu_{t,1}, \\[4pt]
x_{XY}(Q_3)=&0, \\[4pt]
x_{YY}(Q_3)=& \Big( \fr{1}{4}q\sin\fr{\beta}{2}\beta_Y\Big)(Q_3)=-c(\tau(Q_3))\nu_{t,2}.
\end{split}
\end{align}
Hence we acquire the Taylor approximations for $(\tau, t, x)$ by \eqref{5.32}, \eqref{5.35} and \eqref{5.36}
\begin{align}\label{5.37}
\begin{split}
\tau(X,Y)=&\tau^* +\nu_{\tau,1}(X-X_0) +\nu_{\tau,2}(Y-Y_0) \\
& +O(1)\Big(|X-X_0|^2 +|X-X_0||Y-Y_0|+|Y-Y_0|^2\Big),
\end{split}
\end{align}
and
\begin{align}\label{5.38}
\begin{split}
&t(X,Y)=t^* +\fr{1}{2}\nu_{t,1}(X-X_0)^2 +\fr{1}{2}\nu_{t,2}(Y-Y_0)^2 \\
&\quad +O(1)\Big(|X-X_0|^3 +|X-X_0|^2|Y-Y_0|+|X-X_0||Y-Y_0|^2+|Y-Y_0|^3\Big), \\
&x(X,Y)=x^* +\fr{c^*}{2}\nu_{t,1}(X-X_0)^2 -\fr{c^*}{2}\nu_{t,2}(Y-Y_0)^2 \\
&\quad +O(1)\Big(|X-X_0|^3 +|X-X_0|^2|Y-Y_0|+|X-X_0||Y-Y_0|^2+|Y-Y_0|^3\Big).
\end{split}
\end{align}

We now derive a contradiction for this case. If at least one of $\nu_{t,1}<0$ and $\nu_{t,2}<0$ holds, without loss of generality, we assume that $\nu_{t,1}<0$. Then taking $Y=Y_0$ in the expression of $t(X,Y)$ in \eqref{5.38} yields
\begin{align}
\begin{split}
&t(X,Y_0)=t^* +\fr{1}{2}\nu_{t,1}(X-X_0)^2   +O(1)\Big(|X-X_0|^3 \Big), \notag
\end{split}
\end{align}
which implies that there exists a small number $\tilde\delta>0$ such that
\begin{align}\label{5.39}
\begin{split}
&t(X,Y_0)<t^*,  \quad \forall\ X \in U^0(X_0,\tilde\delta),
\end{split}
\end{align}
where $U^0(X_0,\tilde\delta)=\{X:\ 0<|X-X_0|\leq \tilde\delta\}$. Thus $\alpha(X,Y_0)>-\pi$ for $X\in U^0(X_0,\tilde\delta)$.
Recalling \eqref{3.16}, we have
\begin{align}\label{5.40}
\begin{split}
t_X(X,Y_0)=\fr{1}{2c(\tau(X,Y_0))}p(X,Y_0)\cos\fr{\alpha(X,Y_0)}{2},
\end{split}
\end{align}
from which one obtains $t_X(X,Y_0)>0$ for $X\in U^0(X_0,\tilde\delta)$, which contradicts \eqref{5.39}.

If both $\nu_{t,1}>0$ and $\nu_{t,2}>0$ hold, we see by \eqref{5.38} that $t(X, Y)>t^*$ in a punctured neighborhood of $(X_0, Y_0)$. Then starting from the boundary along the line $X=X_0$ toward the point $(X_0,Y_0)$, there exists a point $Y_1<Y_0$ such that $t(X_0,Y_1)=t^*$ and $t(X_0,Y)<t^*$ for $Y<Y_1$.
Since $(X_0,Y_0)$ and $(X_0,Y_1)$ correspond to the same point in the $(x,t)$-plane, $(X_0,Y_1)$ is also a singularity of Type-3. Thus one derives the Taylor approximation for the function $t(X_0, Y)$ at the point $Y_1$ to obtain
\begin{align}\label{5.41}
\begin{split}
&t(X_0,Y)=t^* +\fr{1}{2}\tilde\nu_{t,2}(Y-Y_1)^2   +O(1)\Big(|Y-Y_1|^3 \Big),
\end{split}
\end{align}
where
\begin{align}
\begin{split}
\tilde\nu_{t,2}=t_{YY}(X_0,Y_1)= \Big( -\fr{1}{4c}q\sin\fr{\beta}{2}\beta_Y\Big)(X_0,Y_1)\neq0. \notag
\end{split}
\end{align}
Here we used the fact that $(X_0,Y_1)$ is a Type-3 singular point. Recalling that $t(X_0,Y)<t^*$ for $Y<Y_1$, we find by \eqref{5.41} that $\tilde\nu_{t,2}<0$. By repeating the preceding argument for the case $\nu_{t,1}<0$, we arrive at a contradiction. Hence, singularities of the Type-3 cannot occur on the boundary $t(X,Y)=t^*$.

So far, we have finished the proof of Theorem \ref{thm2}.

\section{Finsler-type norms for smooth solutions}\label{S6}

In this section, we construct two Finsler-type norms on tangent vectors for smooth solutions to measure the transport cost, and then show that they satisfy the Lipschitz property.

\subsection{First order variations}\label{S61}

Let $(\tau, u, R, S)(x,t)$ be a smooth solution to \eqref{2.9}-\eqref{2.9a} for $t\in[0,T]$ with $T<t^*$ satisfying
\begin{align}\label{6.1a}
\int_{\mathbb{R}}|R(x, t)|+|S(x, t)|\ {\rm d}x\leq  \eta,
\end{align}
where $\eta$ is small positive constant to be determined below. We point out that the condition in \eqref{6.1a} can be guaranteed by the initial data. Indeed, one can apply the conservation laws for $(R, S)$
\begin{align}\label{6.1c}
R_t-(cR)_x=0,\qquad S_t+(cS)_x=0,
\end{align}
to obtain that
\begin{align}\label{6.1d}
\int_{\mathbb{R}}|R(x,t)|\ {\rm d}x\leq \int_{\mathbb{R}}|R(x,0)|\ {\rm d}x,\qquad \int_{\mathbb{R}}|S(x,t)|\ {\rm d}x\leq \int_{\mathbb{R}}|S(x,0)|\ {\rm d}x.
\end{align}
Moreover, for sufficiently small $\eta>0$, one can utilize \eqref{6.1a} and the initial assumption $\tau_0(x)-\bar{\tau}_0\in W^{1,1}(\mathbb{R})$ to derive positive upper and lower bounds for $\tau$, that is, there exists a positive constant $\hat{c}$ such that
\begin{align}\label{6.1b}
\fr{1}{\hat{c}}\leq \tau(x,t) \leq \hat{c}.
\end{align}
This argument can be found in, for example, Bressan \cite{Bre}.

We consider a family of perturbed smooth solution $(\tau^\eps, u^\eps, R^\eps, S^\eps)$ to \eqref{2.9}-\eqref{2.9a} as the following form
\begin{align}\label{6.1}
\left\{
\begin{array}{l}
\tau^\eps=\tau+\eps \upsilon+o(\eps),\\
u^\eps=u+\eps \mu+o(\eps),
\end{array}
\right.
\qquad {\rm and}\qquad
\left\{
\begin{array}{l}
R^\eps=R+\eps r +o(\eps),\\
S^\eps=S+\eps s +o(\eps).
\end{array}
\right.
\end{align}
We write $t_{\eps}^*$ for the first time of singularity formation for the solution $(\tau^\eps, R^\eps, S^\eps)$. Let $\eps$ be sufficiently small such that $\inf_{\eps} t_{\eps}^*>T$.
Thus the solution $(\tau, u, R, S)$ and its perturbation $(\tau^\eps, u^\eps, R^\eps, S^\eps)$ are smooth on the domain $\mathbb{R}\times[0,T]$.

From \eqref{2.8a} and \eqref{6.1}, we find that
\begin{align}\label{6.2}
\begin{split}
\dps \tau_{t}^\eps=&\fr{R+S}{2}+\eps\fr{r+s}{2}+o(\eps),\\[4pt]
\dps \tau_{x}^\eps=&\fr{R-S}{2c}+\eps\fr{r-s}{2c}-\eps\fr{c'}{2c^2}(R-S)\upsilon+o(\eps), \\[4pt]
\dps u_{t}^\eps=&\fr{c(R-S)}{2}+\eps\fr{c(r-s)}{2}+\eps\fr{c'}{2}(R-S)\upsilon+o(\eps), \\[4pt]
u_{x}^\eps=&\tau_{t}^\eps.
\end{split}
\end{align}
Given the tangent vectors $r,s$, by \eqref{2.9} and \eqref{6.1}, the perturbations $\upsilon$ and $\mu$ can be uniquely determined by
\begin{align}\label{6.3}
\left\{
\begin{array}{l}
\dps \upsilon_x=-\fr{c'}{2c^2}(R-S)\upsilon +\fr{r-s}{2c}, \\[6pt]
\dps \upsilon(0,t)=0.
\end{array}
\right.
\qquad
\left\{
\begin{array}{l}
\dps \mu_x=\fr{r+s}{2}, \\[6pt]
\dps \mu(0,t)=0.
\end{array}
\right.
\end{align}
Moreover, one also has
\begin{align}\label{6.4}
\begin{split}
\upsilon_t=\fr{r+s}{2},\qquad \mu_t=\fr{c}{2}(r-s)+\fr{c'}{2}(R-S)\upsilon.
\end{split}
\end{align}
Performing direct calculations, we obtain the linear governing closed system for the first order perturbations $\upsilon, r, s$
\begin{align}\label{6.5}
\left\{
\begin{array}{l}
\dps \upsilon_{tt}-c^2(\tau)\upsilon_{xx}=4c c'\tau_x\upsilon_x +2\Big((c')^2\tau_{x}^2 +cc''\tau_{x}^2 +cc'\tau_{xx}\Big)\upsilon, \\[8pt]
\dps r_t-c(\tau)r_x=c'R_x\upsilon +\fr{c''c-(c')^2}{2c^2}(R-S)R\upsilon +\fr{c'}{2c}\Big((2R-S)r -Rs\Big), \\[8pt]
\dps s_t+c(\tau)s_x=-c'S_x\upsilon +\fr{c''c-(c')^2}{2c^2}(S-R)S\upsilon +\fr{c'}{2c}\Big((2S-R)s -Sr\Big).
\end{array}
\right.
\end{align}
Once the solution to system \eqref{6.5} is obtained, the quantity $\mu$ then can be determined.
Recalling the definition $c=\tau^{-\fr{\gamma+1}{2}}$, we see that
\begin{align}\label{6.6}
c'=-\fr{\gamma+1}{2}\tau^{-\fr{\gamma+3}{2}},\quad \fr{c'}{c}=-\fr{\gamma+1}{2\tau},\quad \fr{c''}{c}=\fr{(\gamma+1)(\gamma+3)}{4\tau^2},
\end{align}
which are uniformly bounded for smooth solutions. Here and below, we keep the symbol $c$ in the control systems rather than substituting its explicit expression with respect to $\tau$, since the precise form is not needed for our purposes.

When quantifying the cost of shifting from one solution to another, one naturally considers both vertical and horizontal shifts in the solution space. Note that the tangent flows $\upsilon, \mu$, $r$, and $s$ serve only to characterize vertical shifts of solutions. To describe horizontal shifts of solutions, we further introduce two additional quantities \(w(x,t)\) and \(z(x,t)\) satisfying
\begin{align}\label{6.7}
x_{-}^\eps(t)-x_-(t)=\eps w(x,t)+o(\eps),\qquad x_{+}^\eps(t)-x_+(t)=\eps z(x,t)+o(\eps),
\end{align}
where $x_{\pm}^\eps(t)$ and $x_\pm(t)$ are two forward/backward characteristics starting from initial points $x_{\pm}^\eps(0)$ and $x_\pm(0)$. Obviously, the governing system for $(w,z)$ reads that
\begin{align}\label{6.8}
\left\{
\begin{array}{l}
\dps w_t-cw_x=-c'\Big(\upsilon+\fr{R-S}{2c}w\Big), \\[12pt]
\dps z_t+cz_x=c'\Big(\upsilon+\fr{R-S}{2c}z\Big), \\[8pt]
w(x,0)=w_0(x),\quad z(x,0)=z_0(x).
\end{array}
\right.
\end{align}

Furthermore, to balance the possible increase of energy during
wave interactions, we introduce a pair of interaction potentials $\mathcal{W}^-/\mathcal{W}^+$ for backward/forward directions as follows
\begin{align}\label{6.9}
\begin{split}
\mathcal{W}^-(x,t)=&\eta+\int_{-\infty}^x\fr{S^2}{\sqrt{1+S^2}}(y,t)\ {\rm d}y,\\[4pt] \mathcal{W}^+(x,t)=&\eta+\int_{x}^\infty \fr{R^2}{\sqrt{1+R^2}}(y,t)\ {\rm d}y.
\end{split}
\end{align}
It follows by \eqref{6.1a} and \eqref{6.9} that
\begin{align}\label{6.10}
\begin{split}
\eta\leq \mathcal{W}^-(x,t), \mathcal{W}^+(x,t)\leq 2\eta, \\
\mathcal{W}^\pm(x,t)\leq 2\eta\leq 2\mathcal{W}^\mp(x,t).\
\end{split}
\end{align}
Moreover, we have by \eqref{2.9}
\begin{align}\label{6.11}
\left\{
\begin{array}{l}
\dps (\sqrt{1+R^2})_t-(c\sqrt{1+R^2})_x=-\fr{c'(\tau)}{2c(\tau)}\fr{R-S}{\sqrt{1+R^2}},\\[12pt]
\dps (\sqrt{1+S^2})_t+(c\sqrt{1+S^2})_x=-\fr{c'(\tau)}{2c(\tau)}\fr{S-R}{\sqrt{1+S^2}},
\end{array}
\right.
\end{align}
and then
\begin{align}\label{6.12}
\left\{
\begin{array}{l}
\dps \Big(\fr{S^2}{\sqrt{1+S^2}}\Big)_t +\Big(\fr{cS^2}{\sqrt{1+S^2}}\Big)_x=\fr{c'}{2c}\fr{S^2(S-R)}{(1+S^2)^{3/2}}, \\[12pt]
\dps \Big(\fr{R^2}{\sqrt{1+R^2}}\Big)_t -\Big(\fr{cR^2}{\sqrt{1+R^2}}\Big)_x=\fr{c'}{2c}\fr{R^2(R-S)}{(1+R^2)^{3/2}}.
\end{array}
\right.
\end{align}
Combining \eqref{6.9} and \eqref{6.12} gives
\begin{align*}
\left\{
\begin{array}{l}
\dps \mathcal{W}_t^--c\mathcal{W}_{x}^-=-\fr{2cS^2}{\sqrt{1+S^2}} +\int_{-\infty}^x\fr{c'}{2c}\fr{S^2(S-R)}{(1+S^2)^{3/2}}{\rm d}y, \\[12pt]
\dps \mathcal{W}_t^++c\mathcal{W}_{x}^+=-\fr{2cR^2}{\sqrt{1+R^2}} +\int_{x}^\infty\fr{c'}{2c}\fr{R^2(R-S)}{(1+R^2)^{3/2}}{\rm d}y,
\end{array}
\right.
\end{align*}
from which and \eqref{6.10} one obtains
\begin{align}\label{6.13}
\left\{
\begin{array}{l}
\dps \mathcal{W}_t^--c\mathcal{W}_{x}^-\leq -\fr{2S^2}{\overline{C}\sqrt{1+S^2}} +\overline{C}\eta, \\[12pt]
\dps \mathcal{W}_t^++c\mathcal{W}_{x}^+\leq -\fr{2R^2}{\overline{C}\sqrt{1+R^2}} +\overline{C}\eta.
\end{array}
\right.
\end{align}
Here we retain the notation in Section \ref{S2} that $\overline{C}$ stands for a positive constant arising from the function $c(\tau)$.

\subsection{Finsler norms for tangent vectors}\label{S62}

With the preparation in Section \ref{S61}, we now define the Finsler norm on the space of tangent
vectors $(\upsilon, \mu, r, s)$ and the flow itself $(\tau, u, R, S)$ as follows:
\begin{align}\label{6.15}
\|(\upsilon, \mu, r, s)\|_{(\tau, u, R, S)}:=\inf_{\upsilon, \mu, w, \tilde{r}, z, \tilde{s}}\|(\upsilon, \mu, w, \tilde{r}, z, \tilde{s})\|_{(\tau, u, R, S)},
\end{align}
where the infimum is taken over the set of vertical displacements $\upsilon, \mu, \tilde{r}, \tilde{s}$ and horizontal shifts
$w, z$ which satisfy equations \eqref{6.3}, \eqref{6.4}, \eqref{6.8} and the relations
\begin{align}\label{6.16}
\left\{
\begin{array}{l}
\dps \tilde{r}=r+R_x w-\fr{c'}{4c^2}RS(w-z),\\[12pt]
\dps \tilde{s}=s+S_x z-\fr{c'}{4c^2}RS(w-z).
\end{array}
\right.
\end{align}
Here, to motivate the explicit construction of $\|(\upsilon, \mu, r, s)\|_{(\tau, u, R, S)}$, a reference solution $(R,S)$ and its perturbation $(R^\eps, S^\eps)$ are considered. In first approximation, one has $R^\eps\approx R+\eps r$ ($S^\eps\approx S+\eps s$). On the other hand, the profile $R^\eps$ ($S^\eps$) can be obtained from the graph of $R$ ($S$) through a horizontal shift of $\eps w$ ($\eps z$) succeeded by a vertical shift of $\eps \tilde{r}$ ($\eps \tilde{s}$), provided that
$r=\tilde{r}-R_xw$ ($s=\tilde{s}-S_xz$). The other terms in \eqref{6.16} are the relative shifts.

We define two classes norms of $\|(\upsilon, \mu, w, \tilde{r}, z, \tilde{s})\|_{(\tau, u, R, S)}$ as the following:
\begin{align}\label{6.17a}
\|(\upsilon, \mu, w, \tilde{r}, z, \tilde{s})\|_{(\tau, u, R, S)}^{(i)}:=\sum_{\ell=1}^4\mathcal{I}_\ell,
\end{align}
and
\begin{align}\label{6.17}
\|(\upsilon, \mu, w, \tilde{r}, z, \tilde{s})\|_{(\tau, u, R, S)}^{(ii)}:=\sum_{\ell=1}^6\mathcal{I}_\ell,
\end{align}
where
\begin{align}\label{6.18}
\begin{split}
\mathcal{I}_\ell=\int_{\mathbb{R}}\mathcal{J}_{\ell}^-\mathcal{W}^-\ {\rm d}x +\int_{\mathbb{R}}\mathcal{J}_{\ell}^+\mathcal{W}^+\ {\rm d}x,
\end{split}
\end{align}
and
\begin{align}\label{6.19}
\begin{split}
\mathcal{J}_{1}^-=&|w|\sqrt{1+R^2},\qquad \qquad  \qquad \quad \ \
\mathcal{J}_{1}^+=|z|\sqrt{1+S^2}, \\
\mathcal{J}_{2}^-=&\Big|\upsilon +\fr{Rw}{2c}-\fr{Sz}{2c}\Big|\sqrt{1+R^2}, \qquad \mathcal{J}_{2}^+=\Big|\upsilon +\fr{Rw}{2c}-\fr{Sz}{2c}\Big|\sqrt{1+S^2},  \\
\mathcal{J}_{3}^-=&\Big|\mu +\fr{Rw+Sz}{2}\Big|\sqrt{1+R^2},\qquad
\mathcal{J}_{3}^+=\Big|\mu +\fr{Rw+Sz}{2}\Big|\sqrt{1+S^2},
\\
\mathcal{J}_{4}^-=&|r+R_xw+Rw_x|,\qquad \qquad \quad \ \mathcal{J}_{4}^+=|s+S_xz+Sz_x|, \\
\mathcal{J}_{5}^-=& \fr{|\tilde{r}|}{\sqrt{1+R^2}},\qquad \qquad \qquad \qquad   \ \ \mathcal{J}_{5}^+= \fr{|\tilde{s}|}{\sqrt{1+S^2}}, \\
\mathcal{J}_{6}^-=&\Big|\fr{R\tilde{r}}{\sqrt{1+R^2}} +\sqrt{1+R^2}\Big(w_x+\fr{c'}{4c^2}(w-z)S\Big)\Big|, \\ \mathcal{J}_{6}^+=&\Big|\fr{S\tilde{s}}{\sqrt{1+S^2}} +\sqrt{1+S^2}\Big(z_x+\fr{c'}{4c^2}(w-z)R\Big)\Big|.
\end{split}
\end{align}

We next provide some detailed explanations for $\mathcal{I}_\ell\ (\ell=1,\cdots,6)$.
Roughly speaking, the terms $\mathcal{J}_{\ell}^-\ (\ell=1,\cdots,6)$ mean that
\begin{align*}
\begin{split}
\mathcal{J}_{1}^-\approx &[{\rm change\ in}\ x]\cdot \sqrt{1+R^2}, \\
\mathcal{J}_{2}^-\approx &[{\rm change\ in}\ \tau]\cdot \sqrt{1+R^2}, \\
\mathcal{J}_{3}^-\approx &[{\rm change\ in}\ u]\cdot \sqrt{1+R^2}, \\
\mathcal{J}_{4}^-\approx &[{\rm a\ conserved\ quantity}],\\
\mathcal{J}_{5}^-\approx &[{\rm change\ in}\ \arctan R]\cdot \sqrt{1+R^2},\\
\mathcal{J}_{6}^-\approx &[{\rm change\ in}\ \sqrt{1+R^2}].
\end{split}
\end{align*}
The quantities $\mathcal{J}_{\ell}^+\ (\ell=1,\cdots,6)$ are symmetric for forward waves.

Next we provide a more detailed explanation for each quantity in \eqref{6.17}.

(i) For $\mathcal{I}_1$, the integrand $|w|\sqrt{1+R^2}$ ($|z|\sqrt{1+S^2}$) represents the transportation cost for the base measure of density $\sqrt{1+R^2}$ ($\sqrt{1+S^2}$) moving from $x$ to $x+\eps w(x)$ ($x+\eps z(x)$).

(ii) For $\mathcal{I}_2$, we compute the change in $\tau$ by the relation $x^\eps=x+\eps w+o(\eps)$
\begin{align}\label{6.20}
\begin{split}
\fr{\tau^\eps(x^\eps,t)-\tau(x,t)}{\eps}
&=\upsilon +\tau_xw +o(1) \\
&=\upsilon +\fr{R-S}{2c}w +o(1) \\
&=\Big(\upsilon +\fr{Rw}{2c} -\fr{Sz}{2c}\Big) +\fr{S(z-w)}{2c} +o(1).
\end{split}
\end{align}
Here the term $S(z-w)/2c$ on the right-hand side of \eqref{6.20} is just balanced with
the relative shift term.

(iii) For $\mathcal{I}_3$, one calculates as in \eqref{6.20}
\begin{align}\label{6.20a}
\begin{split}
\fr{u^\eps(x^\eps,t)-u(x,t)}{\eps}
&=\mu +\fr{R+S}{2}w +o(1) \\
&=\Big(\mu +\fr{Rw +Sz}{2} \Big) +\fr{S(w-z)}{2} +o(1).
\end{split}
\end{align}
Here the term $S(w-z)/2$ on the right-hand side of \eqref{6.20a} is also balanced with
the relative shift term.

(iv) For $\mathcal{I}_4$, we shall show the following conservation law
\begin{align}\label{6.21}
\begin{split}
(r+R_x w+Rw_x)_t-[c(r+R_x w+Rw_x)]_x=0,
\end{split}
\end{align}
which will be applied to close the chain of estimates.

(v) For $\mathcal{I}_5$, we can calculate
\begin{align*}
\begin{split}
\arctan R^\eps(x^\eps, t)=&\arctan\Big(R(x^\eps, t) +\eps r(x^\eps, t)+o(\eps)\Big) \\
=&\arctan\Big(R +\eps wR_x +\eps r+o(\eps)\Big)
 \\
=&\arctan R +\eps \fr{r+wR_x}{1+R^2} +o(\eps),
\end{split}
\end{align*}
from which and \eqref{6.16} one obtains
\begin{align*}
\begin{split}
\fr{\arctan R^\eps(x^\eps, t)-\arctan R}{\eps}\sqrt{1+R^2}= \fr{\tilde{r}}{\sqrt{1+R^2}} +\fr{c'}{4c^2}\fr{RS}{\sqrt{1+R^2}}(w-z) +o(1).
\end{split}
\end{align*}
The term $c'RS(w-z)/4c^2\sqrt{1+R^2}$ on the right hand side of the above equality is the relative shift term.

(vi) For $\mathcal{I}_6$, we can also deduce
\begin{align*}
\begin{split}
(R^\eps(x^\eps,t))^2=&R^2(x^\eps,t) +2\eps R(x^\eps,t)r(x^\eps, t)+o(\eps) \\
=& R^2(x,t) +2\eps w(x,t)R(x,t)R_x(x,t)+2\eps R(x,t)r(x,t)+o(\eps),
\end{split}
\end{align*}
from which we find that
\begin{align*}
\begin{split}
&\fr{\sqrt{1+(R^\eps(x^\eps,t))^2}\ {\rm d}x^\eps -\sqrt{1+R^2(x,t)}\ {\rm d}x}{\eps} \\
=& \bigg\{\fr{R(wR_x+r)}{\sqrt{1+R^2}} +\sqrt{1+R^2}w_x +o(1)\bigg\}\ {\rm d}x \\
=& \bigg\{\Big[\fr{R\tilde{r}}{\sqrt{1+R^2}} +\sqrt{1+R^2}\Big(w_x+\fr{c'}{4c^2}(w-z)S \Big)\Big]+ \fr{c'}{4c^2}\fr{S}{\sqrt{1+R^2}}(z-w)+o(1)\bigg\}\ {\rm d}x.
\end{split}
\end{align*}
The term $c'S(z-w)/4c^2\sqrt{1+R^2}$ on the right hand side of the above equality is also the relative shift term.

\begin{rem}\label{rem_12}
When considering $w=z=0$, i.e. without horizontal transportation (or no shift), we see that the metric defined in \eqref{6.17a} can be written as
\begin{align}\label{6.19a}
\begin{split}
&\|(\upsilon, \mu, w, \tilde{r}, z, \tilde{s})\|_{(\tau, u, R, S)}^{(i)} \\
=&\int_{\mathbb{R}}(|\tau^\eps -\tau|+|u^\eps-u|)\Big(\sqrt{1+R^2}\mathcal{W}^- +\sqrt{1+S^2}\mathcal{W}^+\Big)\ {\rm d}x  \\
&+\int_{\mathbb{R}}(|R^\eps -R|+|S^\eps-S|)\ {\rm d}x.
\end{split}
\end{align}
The metric defined in \eqref{6.17} includes more detailed information on the difference between two solutions. Especially,
$\mathcal{J}_{5}^\pm$ terms include the difference of some function of derivative of solution, under transportation. For example,
\begin{equation}
\label{metric_new}
\frac{|\tilde{r}(1+R^2)^{-\frac{1}{2}}|}{\varepsilon}
\approx[\hbox{Variation of }\log(R+\sqrt{1+R^2})]+
o(1).
\end{equation}
When $w=z=0$, i.e. when there is no horizontal transportation (or no shift), \eqref{metric_new} becomes an equation.
\end{rem}

\subsection{Estimates on the norm of tangent vectors}\label{S63}

Denote
\begin{align*}
\begin{split}
\|(\upsilon, \mu, r, s)\|_{(\tau, u, R, S)}^{(i)}=&\inf_{\upsilon, \mu, w, \tilde{r}, z, \tilde{s}}\|(\upsilon, \mu, w, \tilde{r}, z, \tilde{s})\|_{(\tau, u, R, S)}^{(i)},\\
\|(\upsilon, \mu, r, s)\|_{(\tau, u, R, S)}^{(ii)}=&\inf_{\upsilon, \mu, w, \tilde{r}, z, \tilde{s}}\|(\upsilon, \mu, w, \tilde{r}, z, \tilde{s})\|_{(\tau, u, R, S)}^{(ii)}.
\end{split}
\end{align*}
Now, we show that both norms of tangent vectors defined in \eqref{6.17a} and \eqref{6.17} satisfy Lipschitz property.

\begin{lem}\label{lem6}
Given $T>0$, let $(\tau, u, R, S)(x,t)$ be a smooth solution to \eqref{2.9}-\eqref{2.9a} for $t\in[0,T]$ satisfying \eqref{6.1a}. Assume that the first order perturbations $(\upsilon, \mu, r, s)$ satisfy the corresponding system \eqref{6.3}-\eqref{6.5}. Then there exist positive constants $\overline{C}_1$ and $\overline{C}_2$ depending only on the initial data, such that the Lipschitz properties hold
\begin{align}\label{6.22a}
\begin{split}
\|(\upsilon, \mu, r, s)(t)\|_{(\tau, u, R, S)(t)}^{(i)}\leq e^{\overline{C}_1T}\|(\upsilon, \mu, r, s)(0)\|_{(\tau, u, R, S)(0)}^{(i)},
\end{split}
\end{align}
provided that $\eta\leq 1/\overline{C}_1$, and
\begin{align}\label{6.22}
\begin{split}
\|(\upsilon, \mu, r, s)(t)\|_{(\tau, u, R, S)(t)}^{(ii)}\leq e^{\overline{C}_2T}\|(\upsilon, \mu, r, s)(0)\|_{(\tau, u, R, S)(0)}^{(ii)}.
\end{split}
\end{align}
provided that $\eta\leq 1/\overline{C}_2$.
\end{lem}
\begin{proof}
For the purpose of verifying \eqref{6.22a} and \eqref{6.22}, we only need to demonstrate by the definition in \eqref{6.15} that the following Gr\"onwall-type inequalities are valid
\begin{align}\label{6.23a}
\begin{split}
\fr{\rm d}{{\rm d}t}\|(\upsilon, \mu, w, \tilde{r}, z, \tilde{s})(t)\|_{(\tau, u, R, S)(t)}^{(i)}\leq \overline{C}_1\|(\upsilon, \mu, w, \tilde{r}, z, \tilde{s})(t)\|_{(\tau, u, R, S)(t)}^{(i)},
\end{split}
\end{align}
and
\begin{align}\label{6.23}
\begin{split}
\fr{\rm d}{{\rm d}t}\|(\upsilon, \mu, w, \tilde{r}, z, \tilde{s})(t)\|_{(\tau, u, R, S)(t)}^{(ii)}\leq \overline{C}_2\|(\upsilon, \mu, w, \tilde{r}, z, \tilde{s})(t)\|_{(\tau, u, R, S)(t)}^{(ii)},
\end{split}
\end{align}
for any $w,z$ and $\tilde{r}, \tilde{s}$ satisfying \eqref{6.8} and \eqref{6.16}.
For the derivations of \eqref{6.23a} and \eqref{6.23}, we will make repeated use of the following fact:
\begin{align}\label{6.24}
\begin{split}
\fr{{\rm d}}{{\rm d}t}\int_{\mathbb{R}} |f|\mathcal{W}^\pm\ {\rm d}x\leq \int_{\mathbb{R}} |f_t\pm(cf)_x|\mathcal{W}^\pm + |f|(\mathcal{W}^{\pm}_t\pm c\mathcal{W}^{\pm}_x)\ {\rm d}x,
\end{split}
\end{align}
where $f$ is smooth function vanishing at infinity. The proof is divided into several steps.

\textbf{Step 1 (Derivation of $\mathcal{I}_{1}'$).} We combine \eqref{6.8} and \eqref{6.11} to compute
\begin{align*}
\begin{split}
&(w\sqrt{1+R^2})_t-(cw\sqrt{1+R^2})_x \\
=& w[(\sqrt{1+R^2})_t -(c\sqrt{1+R^2})_x] +\sqrt{1+R^2}(w_t-cw_x) \\
=&-\fr{c'}{2c}\fr{R-S}{\sqrt{1+R^2}}w +\sqrt{1+R^2}\Big(-c'(\upsilon+\fr{R-S}{2c}w)\Big) \\
=&-\fr{c'}{2c}\fr{R-S}{\sqrt{1+R^2}}w -c'\Big(v+\fr{Rw}{2c}-\fr{Sz}{2c}\Big)\sqrt{1+R^2} +\fr{c'}{2c}S(w-z)\sqrt{1+R^2} \\
=&-c'\Big(\upsilon+\fr{Rw}{2c}-\fr{Sz}{2c}\Big)\sqrt{1+R^2}  +\fr{c'}{2c}\Big(S-\fr{R-S}{1+R^2}\Big)w\sqrt{1+R^2} -\fr{c'}{2c}Sz\sqrt{1+R^2},
\end{split}
\end{align*}
from which and \eqref{6.13}, \eqref{6.24} we obtain
\begin{align}\label{6.25}
\begin{split}
&\fr{{\rm d}}{{\rm d}t}\int_{\mathbb{R}} \mathcal{J}_{1}^-\mathcal{W}^-\ {\rm d}x=\fr{{\rm d}}{{\rm d}t}\int_{\mathbb{R}}|w|\sqrt{1+R^2}\mathcal{W}^-\ {\rm d}x \\
\leq &\int_{\mathbb{R}}\Big|(w\sqrt{1+R^2})_t -(cw\sqrt{1+R^2})_x\Big|\mathcal{W}^- +|w|\sqrt{1+R^2}(\mathcal{W}^{-}_t-c\mathcal{W}^{-}_x)\ {\rm d}x \\
\leq &\int_{\mathbb{R}}|c'|\Big|\upsilon+\fr{Rw}{2c}-\fr{Sz}{2c}\Big|\sqrt{1+R^2}\mathcal{W}^-\ {\rm d}x +\int_{\mathbb{R}}\fr{|c'|}{c}\Big|S-\fr{R-S}{1+R^2}\Big||w|\sqrt{1+R^2}\mathcal{W}^-\ {\rm d}x \\ &+\int_{\mathbb{R}}\fr{|c'|}{c}|S||z|\sqrt{1+R^2}\mathcal{W}^-\ {\rm d}x +\int_{\mathbb{R}}\Big(-\fr{2S^2}{\overline{C}\sqrt{1+S^2}}+\overline{C}\eta \Big)|w|\sqrt{1+R^2}\ {\rm d}x.
\end{split}
\end{align}
Note that
\begin{align}\label{6.26}
\begin{split}
|S|\leq \sqrt{1+S^2} \leq 1+\fr{S^2}{\sqrt{1+S^2}}.
\end{split}
\end{align}
Putting \eqref{6.26} into \eqref{6.25} and utilizing \eqref{6.10} yields
\begin{align}\label{6.27}
\begin{split}
&\fr{{\rm d}}{{\rm d}t}\int_{\mathbb{R}} \mathcal{J}_{1}^-\mathcal{W}^-\ {\rm d}x \\
\leq &\overline{C}\int_{\mathbb{R}} \mathcal{J}_{2}^-\mathcal{W}^-\ {\rm d}x +\overline{C}\int_{\mathbb{R}} \mathcal{J}_{1}^-\mathcal{W}^-\ {\rm d}x  +\overline{C}\int_{\mathbb{R}}\sqrt{1+R^2}\mathcal{J}_{1}^+\mathcal{W}^+\ {\rm d}x  \\ & -\int_{\mathbb{R}}\fr{S^2}{\overline{C}\sqrt{1+S^2}}(1-\overline{C}\eta)\mathcal{J}_{1}^-\ {\rm d}x \\
\leq & \overline{C}\int_{\mathbb{R}} \mathcal{J}_{2}^-\mathcal{W}^-\ {\rm d}x +\overline{C}\int_{\mathbb{R}} \mathcal{J}_{1}^-\mathcal{W}^-\ {\rm d}x  +\overline{C}\int_{\mathbb{R}}\mathcal{J}_{1}^+\mathcal{W}^+\ {\rm d}x  \\ & +\overline{C}\eta\int_{\mathbb{R}}\fr{R^2}{\sqrt{1+R^2}}\mathcal{J}_{1}^+ \ {\rm d}x -\int_{\mathbb{R}}\fr{S^2}{\overline{C}\sqrt{1+S^2}}(1-\overline{C}\eta)\mathcal{J}_{1}^-\ {\rm d}x.
\end{split}
\end{align}
Similarly, one has
\begin{align}\label{6.28}
\begin{split}
\fr{{\rm d}}{{\rm d}t}\int_{\mathbb{R}} \mathcal{J}_{1}^+\mathcal{W}^+\ {\rm d}x
\leq &\overline{C}\int_{\mathbb{R}} \mathcal{J}_{2}^+\mathcal{W}^+\ {\rm d}x +\overline{C}\int_{\mathbb{R}} \mathcal{J}_{1}^+\mathcal{W}^+\ {\rm d}x  +\overline{C}\int_{\mathbb{R}}\mathcal{J}_{1}^-\mathcal{W}^-\ {\rm d}x  \\ & +\overline{C}\eta\int_{\mathbb{R}}\fr{S^2}{\sqrt{1+S^2}}\mathcal{J}_{1}^- \ {\rm d}x -\int_{\mathbb{R}}\fr{R^2}{\overline{C}\sqrt{1+R^2}}(1-\overline{C}\eta)\mathcal{J}_{1}^+\ {\rm d}x.
\end{split}
\end{align}
Adding \eqref{6.27} and \eqref{6.28} arrives at
\begin{align}\label{6.29}
\begin{split}
\fr{{\rm d}}{{\rm d}t}\mathcal{I}_1
\leq &\overline{C}\mathcal{I}_1 + \overline{C}\mathcal{I}_2 -\fr{1-\overline{C}\eta}{\overline{C}}\int_{\mathbb{R}}\fr{R^2}{\sqrt{1+R^2}}\mathcal{J}_{1}^+\ {\rm d}x \\
&-\fr{1-\overline{C}\eta}{\overline{C}}\int_{\mathbb{R}}\fr{S^2}{\sqrt{1+S^2}}\mathcal{J}_{1}^-\ {\rm d}x.
\end{split}
\end{align}

\textbf{Step 2 (Derivation of $\mathcal{I}_{2}'$).} We directly calculate by \eqref{2.9} and \eqref{6.8}
\begin{align}\label{6.30}
\begin{split}
&\Big(\fr{Rw}{2c}-\fr{Sz}{2c}\Big)_t-c\Big(\fr{Rw}{2c}-\fr{Sz}{2c}\Big)_x \\
=&\fr{w}{2c}(R_t-cR_x)-\fr{z}{2c}(S_t+cS_x) +zS_x \\
&+\fr{R}{2c}(w_t-cw_x) -\fr{S}{2c}(z_t+cz_x) +Sz_x -\fr{c'}{2c^2}(RSw-S^2z) \\
=&\fr{c'}{4c^2}(R-S)Rw -\fr{c'}{4c^2}(S-R)Sz +zS_x \\
&-\fr{c'}{2c}R\Big(\upsilon+\fr{R-S}{2c}w\Big) -\fr{c'}{2c}S\Big(\upsilon+\fr{R-S}{2c}z\Big) +Sz_x-\fr{c'}{2c^2}(RSw-S^2z) \\
=&-\fr{c'}{2c}(R-S)\upsilon -\fr{c'}{c}S\Big(\upsilon+\fr{Rw}{2c}-\fr{Sz}{2c}\Big)+Sz_x+zS_x.
\end{split}
\end{align}
Combining \eqref{6.11} and \eqref{6.30} and applying \eqref{6.3}-\eqref{6.4}, one deduces by
simplifying the result
\begin{align}\label{6.31}
\begin{split}
&\Big[\Big(\upsilon+\fr{Rw}{2c}-\fr{Sz}{2c}\Big)\sqrt{1+R^2}\Big]_t -\Big[c\Big(\upsilon+\fr{Rw}{2c}-\fr{Sz}{2c}\Big)\sqrt{1+R^2}\Big]_x \\
=&\Big(\upsilon+\fr{Rw}{2c}-\fr{Sz}{2c}\Big)\Big((\sqrt{1+R^2})_t-(c\sqrt{1+R^2})_x\Big) +\sqrt{1+R^2}(\upsilon_t-c\upsilon_x) \\ &+\sqrt{1+R^2}\Big[\Big(\fr{Rw}{2c}-\fr{Sz}{2c}\Big)_t-c\Big(\fr{Rw}{2c}-\fr{Sz}{2c}\Big)_x\Big]
 \\
=&-\fr{c'}{c}\Big(S+\fr{R-S}{2(1+R^2)}\Big)\Big(\upsilon+\fr{Rw}{2c}-\fr{Sz}{2c}\Big)\sqrt{1+R^2} \\
&+ (s+zS_x +Sz_x)\sqrt{1+R^2}.
\end{split}
\end{align}
Analogously, there holds
\begin{align}\label{6.32}
\begin{split}
&\Big[\Big(\upsilon+\fr{Rw}{2c}-\fr{Sz}{2c}\Big)\sqrt{1+S^2}\Big]_t +\Big[c\Big(\upsilon+\fr{Rw}{2c}-\fr{Sz}{2c}\Big)\sqrt{1+S^2}\Big]_x \\
=&-\fr{c'}{c}\Big(R+\fr{S-R}{2(1+S^2)}\Big)\Big(\upsilon+\fr{Rw}{2c}-\fr{Sz}{2c}\Big)\sqrt{1+S^2} \\
&+ (r+wR_x +Rw_x)\sqrt{1+S^2}.
\end{split}
\end{align}
By \eqref{6.31}, we employ \eqref{6.24}, \eqref{6.26} and \eqref{6.10} again to obtain
\begin{align}\label{6.33}
\begin{split}
&\fr{{\rm d}}{{\rm d}t}\int_{\mathbb{R}} \mathcal{J}_{2}^-\mathcal{W}^-\ {\rm d}x=\fr{{\rm d}}{{\rm d}t}\int_{\mathbb{R}} \Big|\upsilon+\fr{Rw}{2c}-\fr{Sz}{2c}\Big|\sqrt{1+R^2}\mathcal{W}^-\ {\rm d}x \\
\leq & \overline{C}\int_{\mathbb{R}}  |S|\Big|\upsilon+\fr{Rw}{2c}-\fr{Sz}{2c}\Big|\sqrt{1+R^2}\mathcal{W}^-\ {\rm d}x \\
&+\int_{\mathbb{R}}  \Big|s+zS_x+Sz_x\Big|\sqrt{1+R^2}\mathcal{W}^-\ {\rm d}x \\
&+\int_{\mathbb{R}}  \Big(\overline{C}\eta-\fr{2S^2}{\overline{C}\sqrt{1+S^2}}\Big) \Big|v+\fr{Rw}{2c}-\fr{Sz}{2c}\Big|\sqrt{1+R^2}\ {\rm d}x
\\
\leq & \overline{C}\int_{\mathbb{R}}  \mathcal{J}_{2}^-\mathcal{W}^-\ {\rm d}x  +2\int_{\mathbb{R}}  \mathcal{J}_{4}^+\sqrt{1+R^2}\mathcal{W}^+\ {\rm d}x  -\fr{1-\overline{C}\eta}{\overline{C}}\int_{\mathbb{R}}  \fr{S^2}{\sqrt{1+S^2}} \mathcal{J}_{2}^-\ {\rm d}x.
\end{split}
\end{align}
Performing a similar argument on \eqref{6.32} yields
\begin{align}\label{6.34}
\begin{split}
\fr{{\rm d}}{{\rm d}t}\int_{\mathbb{R}} \mathcal{J}_{2}^+\mathcal{W}^+\ {\rm d}x
\leq & \overline{C}\int_{\mathbb{R}}  \mathcal{J}_{2}^+\mathcal{W}^+\ {\rm d}x  +2\int_{\mathbb{R}}  \mathcal{J}_{4}^-\sqrt{1+S^2}\mathcal{W}^-\ {\rm d}x  \\ &-\fr{1-\overline{C}\eta}{\overline{C}}\int_{\mathbb{R}}  \fr{R^2}{\sqrt{1+R^2}} \mathcal{J}_{2}^+\ {\rm d}x.
\end{split}
\end{align}
By adding \eqref{6.33} and \eqref{6.34} and utilizing \eqref{6.26}, one gains
\begin{align}\label{6.35}
\begin{split}
\fr{{\rm d}}{{\rm d}t}\mathcal{I}_2
\leq & \overline{C}\mathcal{I}_2  +2\int_{\mathbb{R}}  \Big(1+\fr{R^2}{\sqrt{1+R^2}}\Big)\mathcal{J}_{4}^+\mathcal{W}^+\ {\rm d}x  +2\int_{\mathbb{R}}  \Big(1+\fr{S^2}{\sqrt{1+S^2}}\Big)\mathcal{J}_{4}^-\mathcal{W}^-\ {\rm d}x  \\ &-\fr{1-\overline{C}\eta}{\overline{C}}\int_{\mathbb{R}}  \fr{S^2}{\sqrt{1+S^2}} \mathcal{J}_{2}^-\ {\rm d}x -\fr{1-\overline{C}\eta}{\overline{C}}\int_{\mathbb{R}}  \fr{R^2}{\sqrt{1+R^2}} \mathcal{J}_{2}^+\ {\rm d}x \\
\leq &\overline{C}\mathcal{I}_2  +2\mathcal{I}_4  +\int_{\mathbb{R}}  \fr{2R^2}{\sqrt{1+R^2}}\mathcal{J}_{4}^+\mathcal{W}^+\ {\rm d}x  +\int_{\mathbb{R}}  \fr{2S^2}{\sqrt{1+S^2}}\mathcal{J}_{4}^-\mathcal{W}^-\ {\rm d}x  \\ &-\fr{1-\overline{C}\eta}{\overline{C}}\int_{\mathbb{R}}  \fr{S^2}{\sqrt{1+S^2}} \mathcal{J}_{2}^-\ {\rm d}x -\fr{1-\overline{C}\eta}{\overline{C}}\int_{\mathbb{R}}  \fr{R^2}{\sqrt{1+R^2}} \mathcal{J}_{2}^+\ {\rm d}x.
\end{split}
\end{align}

\textbf{Step 3 (Derivation of $\mathcal{I}_{3}'$).}
By a direct calculation, we see by \eqref{2.9} and \eqref{6.8} that
\begin{align}\label{6.35a}
\begin{split}
&\Big(\fr{Rw+Sz}{2}\Big)_t-c\Big(\fr{Rw+Sz}{2}\Big)_x \\
=&\fr{w}{2}(R_t-cR_x)+\fr{z}{2}(S_t+cS_x) -czS_x  +\fr{R}{2}(w_t-cw_x) +\fr{S}{2}(z_t+cz_x) -cSz_x \\
=&\fr{c'}{2}(S-R)\upsilon- c(Sz_x+zS_x),
\end{split}
\end{align}
from which and \eqref{6.3}-\eqref{6.4} one obtains
\begin{align}\label{6.35b}
\begin{split}
&\Big(\mu +\fr{Rw+Sz}{2}\Big)_t-c\Big(\mu+\fr{Rw+Sz}{2}\Big)_x =- c(s+Sz_x+zS_x).
\end{split}
\end{align}
Combing \eqref{6.11} and \eqref{6.35b} yields
\begin{align}\label{6.35c}
\begin{split}
&\Big[\Big(\mu +\fr{Rw+Sz}{2}\Big)\sqrt{1+R^2}\Big]_t -\Big[c\Big(\mu +\fr{Rw+Sz}{2}\Big)\sqrt{1+R^2}\Big]_x \\
=&-\fr{c'}{2c}\fr{R-S}{\sqrt{1+R^2}}\Big(\mu +\fr{Rw+Sz}{2}\Big) -c(s+zS_x +Sz_x)\sqrt{1+R^2}.
\end{split}
\end{align}
We then have
\begin{align}\label{6.35d}
\begin{split}
&\fr{{\rm d}}{{\rm d}t}\int_{\mathbb{R}} \mathcal{J}_{3}^-\mathcal{W}^-\ {\rm d}x=\fr{{\rm d}}{{\rm d}t}\int_{\mathbb{R}} \Big|\mu +\fr{Rw+Sz}{2}\Big|\sqrt{1+R^2}\mathcal{W}^-\ {\rm d}x \\
\leq & \overline{C}\int_{\mathbb{R}}  |S|\Big|\mu +\fr{Rw+Sz}{2}\Big|\sqrt{1+R^2}\mathcal{W}^-\ {\rm d}x \\
&+\int_{\mathbb{R}}  \Big|s+zS_x+Sz_x\Big|\sqrt{1+R^2}\mathcal{W}^-\ {\rm d}x \\
&+\int_{\mathbb{R}}  \Big(\overline{C}\eta-\fr{2S^2}{\overline{C}\sqrt{1+S^2}}\Big) \Big|\mu +\fr{Rw+Sz}{2}\Big|\sqrt{1+R^2}\ {\rm d}x
\\
\leq & \overline{C}\int_{\mathbb{R}}  \mathcal{J}_{3}^-\mathcal{W}^-\ {\rm d}x  +2\int_{\mathbb{R}}  \mathcal{J}_{4}^+\sqrt{1+R^2}\mathcal{W}^+\ {\rm d}x  -\fr{1-\overline{C}\eta}{\overline{C}}\int_{\mathbb{R}}  \fr{S^2}{\sqrt{1+S^2}} \mathcal{J}_{3}^-\ {\rm d}x.
\end{split}
\end{align}
Thus there holds
\begin{align}\label{6.35e}
\begin{split}
\fr{{\rm d}}{{\rm d}t}\mathcal{I}_3
\leq &\overline{C}\mathcal{I}_3  +2\mathcal{I}_4  +\int_{\mathbb{R}}  \fr{2R^2}{\sqrt{1+R^2}}\mathcal{J}_{4}^+\mathcal{W}^+\ {\rm d}x  +\int_{\mathbb{R}}  \fr{2S^2}{\sqrt{1+S^2}}\mathcal{J}_{4}^-\mathcal{W}^-\ {\rm d}x  \\ &-\fr{1-\overline{C}\eta}{\overline{C}}\int_{\mathbb{R}}  \fr{S^2}{\sqrt{1+S^2}} \mathcal{J}_{3}^-\ {\rm d}x -\fr{1-\overline{C}\eta}{\overline{C}}\int_{\mathbb{R}}  \fr{R^2}{\sqrt{1+R^2}} \mathcal{J}_{3}^+\ {\rm d}x.
\end{split}
\end{align}

\textbf{Step 4 (Derivation of $\mathcal{I}_{4}'$).} We first derive the conservation law \eqref{6.21}. From \eqref{6.5}, \eqref{6.8} and \eqref{2.9}, one calculates
\begin{align}\label{6.36}
\begin{split}
&(r+wR_x)_t-(c(r+wR_x))_x \\
=& (r_t-cr_x) -c'\fr{R-S}{2c}r +R_x(w_t-cw_x) +w[R_{xt}-(cR_x)_x] \\
=&\Big(c'R_x\upsilon +\fr{c''c-(c')^2}{2c^2}(R-S)R\upsilon +\fr{c'}{2c}[(2R-S)r-Rs]\Big) \\ &-\fr{c'}{2c}(R-S)r -c'\Big(\upsilon+\fr{R-S}{2c}w\Big)R_x \\
&+w\Big(\fr{c''c-(c')^2}{2c^2}\fr{R-S}{2c}(R-S)R +\fr{c'}{2c}(R_x-S_x)R +\fr{c'}{2c}(R-S)R_x\Big) \\
=&\fr{c'}{2c}R(r-s) +\fr{c'}{2c}(R_x-S_x)Rw +\fr{c''c-(c')^2}{2c^2}(R-S)R\Big(\upsilon+\fr{Rw}{2c}-\fr{Sz}{2c}\Big) \\
&-\fr{c''c-(c')^2}{4c^3}(R-S)RS(w-z).
\end{split}
\end{align}
Moreover, we have by \eqref{6.8} and \eqref{6.3} again
\begin{align*}
\begin{split}
&w_{xt}-(cw_x)_x\\
=&-c''\fr{R-S}{2c}\Big(\upsilon+\fr{R-S}{2c}w\Big)-c'\Big(\upsilon_x+\fr{R-S}{2c}w_x +\fr{R_x-S_x}{2c}w -\fr{R-S}{2c^2}wc'\fr{R-S}{2c}\Big) \\
=&-\fr{c''}{2c}(R-S)\Big(\upsilon+\fr{Rw}{2c}-\fr{Sz}{2c}\Big) +\fr{c''}{4c^2}(R-S)S(w-z) \\
&-c'\Big(-\fr{c'}{2c^2}(R-S)\upsilon +\fr{\tilde{r}-\tilde{s}-wR_x+zS_x}{2c}\Big) \\
&-\fr{c'}{2c}(R-S)w_x -\fr{c'}{2c}(R_x-S_x)w +\fr{(c')^2}{4c^3}(R-S)^2w \\
=& -\fr{c''}{2c}(R-S)\Big(\upsilon+\fr{Rw}{2c}-\fr{Sz}{2c}\Big) +\fr{c''}{4c^2}(R-S)S(w-z) \\
&+\fr{(c')^2}{2c^2}(R-S)\upsilon -\fr{c'}{2c}(\tilde{r}-\tilde{s}) +\fr{c'}{2c}(w-z)S_x -\fr{c'}{2c}(R-S)w_x +\fr{(c')^2}{4c^3}(R-S)^2w,
\end{split}
\end{align*}
from which we obtain
\begin{align}\label{6.37}
\begin{split}
&(Rw_{x})_t-(cRw_x)_x\\
=&R[w_{xt}-(cw_x)_x] +w_x(R_t-cR_x) \\
=&-\fr{c''}{2c}R(R-S)\Big(\upsilon+\fr{Rw}{2c}-\fr{Sz}{2c}\Big) +\fr{c''}{4c^2}(R-S)RS(w-z) \\
&+\fr{(c')^2}{2c^2}R(R-S)\upsilon -\fr{c'}{2c}R(r-s+wR_x-zS_x) +\fr{c'}{2c}R(w-z)S_x \\
&-\fr{c'}{2c}R(R-S)w_x +\fr{(c')^2}{4c^3}R(R-S)^2w +w_x\fr{c'}{2c}(R-S)R \\
=&-\fr{c'}{2c}R(r-s) -\fr{c'}{2c}Rw(R_x-S_x)-\fr{c''c-(c')^2}{2c^2}R(R-S)\Big(\upsilon+\fr{Rw}{2c}-\fr{Sz}{2c}\Big) \\
&-\fr{(c')^2}{4c^3}R(R-S)(Rw-Sz)+\fr{c''}{4c^2}(R-S)RS(w-z)+\fr{(c')^2}{4c^3}R(R-S)^2w \\
=&-\fr{c'}{2c}R(r-s) -\fr{c'}{2c}Rw(R_x-S_x)-\fr{c''c-(c')^2}{2c^2}R(R-S)\Big(\upsilon+\fr{Rw}{2c}-\fr{Sz}{2c}\Big) \\
&+\fr{c'}{2c}(R-S)w_x +\fr{(c')^2}{4c^3}(R-S)^2w.
\end{split}
\end{align}
Adding \eqref{6.36} and \eqref{6.37} gives
\begin{align}\label{6.38}
\begin{split}
(r+R_x w+Rw_x)_t-[c(r+R_x w+Rw_x)]_x=0,
\end{split}
\end{align}
which is the desired conservation law \eqref{6.21}. By symmetry, the following conservation law also holds
\begin{align}\label{6.39}
\begin{split}
(s+S_x z+Sz_x)_t+[c(s+S_x z+Sz_x)]_x=0.
\end{split}
\end{align}

We apply \eqref{6.38}, \eqref{6.24} and \eqref{6.13} to arrive at
\begin{align}\label{6.40}
\begin{split}
&\fr{{\rm d}}{{\rm d}t}\int_{\mathbb{R}} \mathcal{J}_{4}^-\mathcal{W}^-\ {\rm d}x=\fr{{\rm d}}{{\rm d}t}\int_{\mathbb{R}}|r+R_xw+Rw_x|\mathcal{W}^-\ {\rm d}x \\
\leq &\int_{\mathbb{R}}\Big|[(r+R_xw+Rw_x)]_t -[c(r+R_xw+Rw_x)]_x\Big|\mathcal{W}^- \\ &+|r+R_xw+Rw_x|(\mathcal{W}^{-}_t-c\mathcal{W}^{-}_x)\ {\rm d}x \\
= &\int_{\mathbb{R}} |r+R_xw+Rw_x|(\mathcal{W}^{-}_t-c\mathcal{W}^{-}_x) \ {\rm d}x \\
\leq & \int_{\mathbb{R}} \mathcal{J}_{4}^- \Big(\overline{C}\eta-\fr{2S^2}{\overline{C}\sqrt{1+S^2}}\Big)\ {\rm d}x.
\end{split}
\end{align}
Carrying out a similar process on \eqref{6.39} suggests
\begin{align*}
\begin{split}
&\fr{{\rm d}}{{\rm d}t}\int_{\mathbb{R}} \mathcal{J}_{4}^+\mathcal{W}^+\ {\rm d}x
\leq  \int_{\mathbb{R}} \mathcal{J}_{4}^+ \Big(\overline{C}\eta-\fr{2R^2}{\overline{C}\sqrt{1+R^2}}\Big)\ {\rm d}x,
\end{split}
\end{align*}
which together with \eqref{6.40} and \eqref{6.10} gives
\begin{align}\label{6.41}
\begin{split}
\fr{{\rm d}}{{\rm d}t}\mathcal{I}_4
\leq & \int_{\mathbb{R}} \mathcal{J}_{4}^- \Big(\overline{C}\eta-\fr{2S^2}{\overline{C}\sqrt{1+S^2}}\Big)\ {\rm d}x +\int_{\mathbb{R}} \mathcal{J}_{3}^- \Big(\overline{C}\eta-\fr{2S^2}{\overline{C}\sqrt{1+S^2}}\Big)\ {\rm d}x \\
\leq & \overline{C}\mathcal{I}_4 -\int_{\mathbb{R}} \fr{2S^2}{\overline{C}\sqrt{1+S^2}}\mathcal{J}_{4}^-\ {\rm d}x -\int_{\mathbb{R}} \fr{2R^2}{\overline{C}\sqrt{1+R^2}}\mathcal{J}_{4}^+\ {\rm d}x.
\end{split}
\end{align}

\textbf{Step 5 (Derivation of $\mathcal{I}_{5}'$).} We first compute by \eqref{2.9} and \eqref{6.8}
\begin{align*}
\begin{split}
&\Big(\fr{c'}{4c^2}(w-z)RS\Big)_t-\Big(c\cdot\fr{c'}{4c^2}(w-z)RS\Big)_x \\
=& \Big(\fr{c'}{4c^2}\Big)' (\tau_t-c\tau_x)(w-z)RS +\fr{c'}{4c^2}[(w_t-cw_x)-(z_t-cz_x)]RS \\ &+\fr{c'}{4c^2}(w-z)[(RS)_t-(cRS)_x] \\
=&\fr{c''c-2(c')^2}{4c^3}(w-z)RS^2 -\fr{(c')^2}{4c^2}\Big(2\upsilon+\fr{R-S}{2c}(w+z)\Big)RS +\fr{c'}{4c^2}\cdot 2cz_xRS \\
&+\fr{c'}{4c^2}(w-z)\Big(\fr{c'}{2c}(S-R)RS-2cRS_x\Big),
\end{split}
\end{align*}
from which and \eqref{6.36} we obtain
\begin{align}\label{6.42}
\begin{split}
&\tilde{r}_t-(c\tilde{r})_x=\Big[r+wR_x-\fr{c'}{4c^2}(w-z)RS\Big]_t -\Big[c\Big(r+wR_x-\fr{c'}{4c^2}(w-z)RS\Big)\Big]_x \\
=& \fr{c'}{2c}R(\tilde{r}-\tilde{s}-wR_x+zS_x) +\fr{c'}{2c}(R_x-S_x)Rw +\fr{c''c-(c')^2}{2c^2}(R-S)R\Big(\upsilon+\fr{R-S}{2c}w\Big) \\
&-\fr{c''c-2(c')^2}{4c^3}(w-z)RS^2 +\fr{(c')^2}{4c^2}\Big(2\upsilon+\fr{R-S}{2c}(w+z)\Big)RS -\fr{c'}{2c}z_xRS \\
&-\fr{c'}{4c^2}(w-z)\Big(\fr{c'}{2c}(S-R)RS-2cRS_x\Big) \\
=& \fr{c'}{2c}R\tilde{r}-\fr{c'}{2c}R\tilde{s}  +\fr{c''c-(c')^2}{2c^2}(R-S)R\Big(\upsilon+\fr{Rw}{2c}-\fr{Sz}{2c}\Big)  \\
&+\fr{c''c-(c')^2}{2c^2}(R-S)R \cdot\fr{z-w}{2c}S  -\fr{c''c-2(c')^2}{4c^3}(w-z)RS^2 \\
&+\fr{(c')^2}{2c^2}\Big(\upsilon+\fr{Rw}{2c}- \fr{Sz}{2c}\Big)RS+\fr{(c')^2}{2c^2}\Big(\fr{R-S}{4c}(w+z)-\fr{Rw}{2c}+ \fr{Sz}{2c}\Big)RS \\
&-\fr{c'}{2c}z_xRS-\fr{c'}{4c^2}(w-z)\cdot\fr{c'}{2c}(S-R)RS.
\end{split}
\end{align}
Simplifying \eqref{6.42} arrives at
\begin{align}\label{6.43}
\begin{split}
\tilde{r}_t-(c\tilde{r})_x
=& \fr{c'}{2c}R\tilde{r} -\fr{c'}{2c}(R\tilde{s}+ z_xRS) +\fr{c''c-(c')^2}{2c^2}R^2\Big(\upsilon+\fr{Rw}{2c}-\fr{Sz}{2c}\Big) \\ &-\fr{c''c-2(c')^2}{2c^2}RS\Big(\upsilon+\fr{Rw}{2c}-\fr{Sz}{2c}\Big) +\fr{c''c-(c')^2}{4c^3}R^2S(z-w).
\end{split}
\end{align}
One combines \eqref{6.11} and \eqref{6.43} to achieve
\begin{align}\label{6.44}
\begin{split}
&\Big(\fr{\tilde{r}}{\sqrt{1+R^2}}\Big)_t-\Big(\fr{c\tilde{r}}{\sqrt{1+R^2}}\Big)_x \\
=&\fr{c'}{2c}\fr{R+R^2S}{1+R^2}\fr{\tilde{r}}{\sqrt{1+R^2}} -\fr{c'}{2c}\fr{R}{\sqrt{1+R^2}}(s +S_xz+z_xS ) \\ &+\fr{c''c-(c')^2}{2c^2}\fr{R}{\sqrt{1+R^2}}R\Big(\upsilon+\fr{Rw}{2c}-\fr{Sz}{2c}\Big) \\ &-\fr{c''c-2(c')^2}{2c^2}\fr{R}{\sqrt{1+R^2}}S\Big(\upsilon+\fr{Rw}{2c}-\fr{Sz}{2c}\Big) +\fr{2c''c-(c')^2}{8c^3}\fr{R}{\sqrt{1+R^2}}RS(z-w).
\end{split}
\end{align}
Then we have
\begin{align}\label{6.45}
\begin{split}
&\fr{{\rm d}}{{\rm d}t}\int_{\mathbb{R}} \mathcal{J}_{5}^-\mathcal{W}^-\ {\rm d}x=\fr{{\rm d}}{{\rm d}t}\int_{\mathbb{R}}\fr{|\tilde{r}|}{\sqrt{1+R^2}}\mathcal{W}^-\ {\rm d}x \\
\leq &\int_{\mathbb{R}}\Big|\Big(\fr{\tilde{r}}{\sqrt{1+R^2}}\Big)_t -\Big(\fr{c\tilde{r}}{\sqrt{1+R^2}}\Big)_x\Big|\mathcal{W}^- +\fr{|\tilde{r}|}{\sqrt{1+R^2}}(\mathcal{W}^{-}_t-c\mathcal{W}^{-}_x)\ {\rm d}x \\
\leq &\overline{C}\int_{\mathbb{R}} (1+|S|)\fr{|\tilde{r}|}{\sqrt{1+R^2}} \mathcal{W}^-\ {\rm d}x +\overline{C}\int_{\mathbb{R}}  \mathcal{J}_{4}^+\mathcal{W}^-\ {\rm d}x +\overline{C}\int_{\mathbb{R}}  \mathcal{J}_{2}^-\mathcal{W}^-\ {\rm d}x \\ & +\overline{C}\int_{\mathbb{R}}  \mathcal{J}_{2}^+\mathcal{W}^-\ {\rm d}x +\overline{C}\int_{\mathbb{R}}  |R|\mathcal{J}_{1}^+\mathcal{W}^-\ {\rm d}x
+\overline{C}\int_{\mathbb{R}}  |S|\mathcal{J}_{1}^-\mathcal{W}^-\ {\rm d}x \\
&+\int_{\mathbb{R}} \fr{|\tilde{r}|}{\sqrt{1+R^2}}\Big(\overline{C}\eta-\fr{2S^2}{\overline{C}\sqrt{1+S^2}}\Big)\ {\rm d}x.
\end{split}
\end{align}
Note that
\begin{align*}
\begin{split}
1+|S|\leq 2+\fr{S^2}{\sqrt{1+S^2}},\qquad \mathcal{W}^\pm\leq 2\mathcal{W}^\mp.
\end{split}
\end{align*}
One simplifies \eqref{6.45} to obtain
\begin{align}\label{6.46}
\begin{split}
\fr{{\rm d}}{{\rm d}t}\int_{\mathbb{R}} \mathcal{J}_{5}^-\mathcal{W}^-\ {\rm d}x
\leq &\overline{C}\int_{\mathbb{R}} \mathcal{J}_{5}^- \mathcal{W}^-\ {\rm d}x +\overline{C}\int_{\mathbb{R}}  \mathcal{J}_{4}^+\mathcal{W}^+\ {\rm d}x +\overline{C}\mathcal{I}_2 +\overline{C}\mathcal{I}_1 \\
 &+\overline{C}\int_{\mathbb{R}}  \fr{R^2}{\sqrt{1+R^2}}\mathcal{J}_{1}^+\mathcal{W}^+\ {\rm d}x
+\overline{C}\int_{\mathbb{R}}  \fr{S^2}{\sqrt{1+S^2}}\mathcal{J}_{1}^-\mathcal{W}^-\ {\rm d}x \\ &-\fr{1-\overline{C}\eta}{\overline{C}}\int_{\mathbb{R}} \fr{S^2}{\sqrt{1+S^2}} \mathcal{J}_{5}^-\ {\rm d}x.
\end{split}
\end{align}
Hence, after deriving an estimate for $\mathcal{J}_{5}^+$ analogous to that in \eqref{6.46}, we gain the estimate for $\mathcal{I}_5$
\begin{align}\label{6.47}
\begin{split}
\fr{{\rm d}}{{\rm d}t}\mathcal{I}_5
\leq &\overline{C}\mathcal{I}_5 +\overline{C}\mathcal{I}_4 +\overline{C}\mathcal{I}_2 +\overline{C}\mathcal{I}_1 \\
 &+\overline{C}\eta\int_{\mathbb{R}}  \fr{R^2}{\sqrt{1+R^2}}\mathcal{J}_{1}^+\ {\rm d}x
+\overline{C}\eta\int_{\mathbb{R}}  \fr{S^2}{\sqrt{1+S^2}}\mathcal{J}_{1}^-\ {\rm d}x \\ &-\fr{1-\overline{C}\eta}{\overline{C}}\int_{\mathbb{R}} \fr{S^2}{\sqrt{1+S^2}} \mathcal{J}_{5}^-\ {\rm d}x -\fr{1-\overline{C}\eta}{\overline{C}}\int_{\mathbb{R}} \fr{R^2}{\sqrt{1+R^2}} \mathcal{J}_{5}^+\ {\rm d}x.
\end{split}
\end{align}

\textbf{Step 6 (Derivation of $\mathcal{I}_{6}'$).} We first note that
\begin{align}\label{6.48}
\begin{split}
&\fr{R\tilde{r}}{\sqrt{1+R^2}} +\sqrt{1+R^2}\Big(w_x+\fr{c'}{4c^2}(w-z)S\Big) \\
=&\fr{R}{\sqrt{1+R^2}}\Big(\tilde{r}+Rw_x+\fr{c'}{4c^2}(w-z)RS\Big) +\fr{1}{\sqrt{1+R^2}}\Big(w_x+\fr{c'}{4c^2}(w-z)S\Big)\\
=&\fr{R}{\sqrt{1+R^2}}(r+Rw_x+wR_x) +\fr{1}{\sqrt{1+R^2}}\Big(w_x+\fr{c'}{4c^2}(w-z)S\Big).
\end{split}
\end{align}
Recalling \eqref{2.9} and \eqref{6.38} gives
\begin{align}\label{6.49}
\begin{split}
&\Big(\fr{R}{\sqrt{1+R^2}}(r+Rw_x+wR_x)\Big)_t - \Big(\fr{cR}{\sqrt{1+R^2}}(r+Rw_x+wR_x)\Big)_x \\
=& (r+Rw_x+wR_x)\fr{1}{(1+R^2)^{3/2}}(R_t-cR_x) \\
=& \fr{c'}{2c}\fr{R(R-S)}{(1+R^2)^{3/2}}(r+Rw_x+wR_x).
\end{split}
\end{align}
Moreover, we directly calculate
\begin{align*}
\begin{split}
&\Big(w_x+\fr{c'}{4c^2}(w-z)S\Big)_t-\Big(cw_x+c\fr{c'}{4c^2}(w-z)S\Big)_x \\
=& [w_{xt}-(cw_x)_x] +\Big(\fr{c'}{4c^2}\Big)'(\tau_t-c\tau_x)(w-z)S \\
&+\fr{c'}{4c^2}[(w_t-cw_x)-(z_t-cz_x)]S +\fr{c'}{4c^2}(w-z)[S_t-(cS)_x] \\
=& \bigg\{-\fr{c''}{2c}(R-S)\Big(\upsilon+\fr{Rw}{2c}-\fr{Sz}{2c}\Big) +\fr{c''}{4c^2}(R-S)S(w-z) \\ &+\fr{(c')^2}{2c^2}(R-S)\upsilon -\fr{c'}{2c}(\tilde{r}-\tilde{s}) +\fr{c'}{2c}(w-z)S_x -\fr{c'}{2c}(R-S)w_x +\fr{(c')^2}{4c^3}(R-S)^2w \bigg\} \\
&+\fr{c''c-2(c')^2}{4c^3}(w-z)S^2 -\fr{(c')^2}{8c^3}(R-S)S(w-z) +\fr{c'}{2c}Sz_x \\
&-\fr{c'}{2c}(w-z)S_x -\fr{(c')^2}{4c^3}(R-S)S(w-z).
\end{split}
\end{align*}
Simplifying the above equation leads to
\begin{align}\label{6.50}
\begin{split}
&\Big(w_x+\fr{c'}{4c^2}(w-z)S\Big)_t-\Big(cw_x+c\fr{c'}{4c^2}(w-z)S\Big)_x \\
=&-\fr{c'}{2c}(\tilde{r}-\tilde{s}) -\fr{c''c-(c')^2}{2c^2}(R-S)\Big(\upsilon+\fr{Rw}{2c}-\fr{Sz}{2c}\Big) -\fr{c'}{2c}(R-S)w_x +\fr{c'}{2c}Sz_x \\
&+\fr{c''}{4c^2}RS(w-z) +\fr{(c')^2}{8c^3}(S-5R)S(w-z)
\\
=&-\fr{c'}{2c}(\tilde{r}-\tilde{s}) -\fr{c''c-(c')^2}{2c^2}(R-S)\Big(\upsilon+\fr{Rw}{2c}-\fr{Sz}{2c}\Big) \\
&-\fr{c'}{2c}(R-S)\Big(w_x+\fr{c'}{4c^2}(w-z)S\Big) +\fr{c'}{2c}Sz_x+\fr{c''c-2(c')^2}{4c^3}RS(w-z).
\end{split}
\end{align}
One applies \eqref{6.50} and \eqref{2.9} to achieve
\begin{align}\label{6.51}
\begin{split}
&\Big[\fr{1}{\sqrt{1+R^2}}\Big(w_x+\fr{c'}{4c^2}(w-z)S\Big)\Big]_t -\Big[\fr{c}{\sqrt{1+R^2}}\Big(w_x+\fr{c'}{4c^2}(w-z)S\Big)\Big]_x \\
=&-\fr{c'}{2c}\fr{\tilde{r}}{\sqrt{1+R^2}} +\fr{c'}{2c}\fr{1}{\sqrt{1+R^2}}(\tilde{s}+Sz_x) -\fr{c''c-(c')^2}{2c^2}\fr{R-S}{\sqrt{1+R^2}}\Big(\upsilon+\fr{Rw}{2c}-\fr{Sz}{2c}\Big) \\ &-\fr{c'}{2c}\fr{(R-S)+2R^2}{1+R^2}\fr{1}{\sqrt{1+R^2}}\Big(w_x+\fr{c'}{4c^2}(w-z)S\Big) \\ &+\fr{c''c-2(c')^2}{4c^3}\fr{1}{\sqrt{1+R^2}}RS(w-z) \\
=&-\fr{c'}{2c}\fr{\tilde{r}}{\sqrt{1+R^2}} +\fr{c'}{2c}\fr{1}{\sqrt{1+R^2}}(s+Sz_x +S_xz) \\ &-\fr{c''c-(c')^2}{2c^2}\fr{R-S}{\sqrt{1+R^2}}\Big(\upsilon+\fr{Rw}{2c}-\fr{Sz}{2c}\Big) \\ &-\fr{c'}{2c}\fr{(R-S)+2R^2}{1+R^2}\fr{1}{\sqrt{1+R^2}}\Big(w_x+\fr{c'}{4c^2}(w-z)S\Big) \\ &+\fr{2c''c-5(c')^2}{8c^3}\fr{1}{\sqrt{1+R^2}}RS(w-z).
\end{split}
\end{align}
Combining \eqref{6.48}, \eqref{6.49} and \eqref{6.51} yields
\begin{align}\label{6.52}
\begin{split}
&\Big[\fr{R\tilde{r}}{\sqrt{1+R^2}} +\sqrt{1+R^2}\Big(w_x+\fr{c'}{4c^2}(w-z)S\Big)\Big]_t \\ &-\Big[\fr{cR\tilde{r}}{\sqrt{1+R^2}} +c\sqrt{1+R^2}\Big(w_x+\fr{c'}{4c^2}(w-z)S\Big)\Big]_x \\
=& \fr{c'}{c}\fr{R[(R-S)+R^2]}{(1+R^2)^{3/2}}(r+Rw_x+wR_x) -\fr{c'}{2c}\fr{\tilde{r}}{\sqrt{1+R^2}} \\ &+\fr{c'}{2c}\fr{1}{\sqrt{1+R^2}}(s+Sz_x +S_xz)-\fr{c''c-(c')^2}{2c^2}\fr{R-S}{\sqrt{1+R^2}}\Big(\upsilon+\fr{Rw}{2c}-\fr{Sz}{2c}\Big) \\ &-\fr{c'}{2c}\fr{(R-S)+2R^2}{(1+R^2)^{3/2}}\Big[\fr{R\tilde{r}}{\sqrt{1+R^2}} +\sqrt{1+R^2}\Big(w_x+\fr{c'}{4c^2}(w-z)S\Big)\Big] \\ &+\fr{2c''c-5(c')^2}{8c^3}\fr{1}{\sqrt{1+R^2}}RS(w-z).
\end{split}
\end{align}
Note that
$$
\fr{|R[(R-S)+R^2]|}{(1+R^2)^{3/2}}\leq 2+|S|\leq 3+\fr{S^2}{\sqrt{1+S^2}}
$$
Then we have by \eqref{6.24}, \eqref{6.13} and \eqref{6.52}
\begin{align*}
\begin{split}
&\fr{{\rm d}}{{\rm d}t}\int_{\mathbb{R}} \mathcal{J}_{6}^-\mathcal{W}^-\ {\rm d}x
= \fr{{\rm d}}{{\rm d}t}\int_{\mathbb{R}} \Big|\fr{R\tilde{r}}{\sqrt{1+R^2}} +\sqrt{1+R^2}\Big(w_x+\fr{c'}{4c^2}(w-z)S\Big)\Big|\mathcal{W}^-\ {\rm d}x \\
\leq & \overline{C}\int_{\mathbb{R}} \Big(1+\fr{S^2}{\sqrt{1+S^2}}\Big)|r+Rw_x+wR_x| \mathcal{W}^-\ {\rm d}x \\
&+\overline{C}\int_{\mathbb{R}} \fr{|\tilde{r}|}{\sqrt{1+R^2}} \mathcal{W}^-\ {\rm d}x +\overline{C}\int_{\mathbb{R}} \fr{1}{\sqrt{1+R^2}}|s+Sz_x +S_xz|\mathcal{W}^-\ {\rm d}x \\ & +\overline{C}\int_{\mathbb{R}} (\sqrt{1+R^2}+\sqrt{1+S^2})\Big|\upsilon+\fr{Rw}{2c}-\fr{Sz}{2c}\Big|\mathcal{W}^-\ {\rm d}x \\ &+\overline{C}\int_{\mathbb{R}}\Big(1+\fr{S^2}{\sqrt{1+S^2}}\Big)\Big|\fr{R\tilde{r}}{\sqrt{1+R^2}} +\sqrt{1+R^2}\Big(w_x+\fr{c'}{4c^2}(w-z)S\Big)\Big| \mathcal{W}^-\ {\rm d}x \\
&+\overline{C}\int_{\mathbb{R}}  [|S||w|\sqrt{1+R^2}+|z|\sqrt{1+S^2}] \mathcal{W}^-\ {\rm d}x \\
&+\int_{\mathbb{R}}\Big|\fr{R\tilde{r}}{\sqrt{1+R^2}} +\sqrt{1+R^2}\Big(w_x+\fr{c'}{4c^2}(w-z)S\Big)\Big| \Big(\overline{C}\eta-\fr{2S^2}{\overline{C}\sqrt{1+S^2}}\Big)\ {\rm d}x,
\end{split}
\end{align*}
which together with \eqref{6.10} obtains
\begin{align}\label{6.53}
\begin{split}
&\fr{{\rm d}}{{\rm d}t}\int_{\mathbb{R}} \mathcal{J}_{6}^-\mathcal{W}^-\ {\rm d}x \\
\leq & \overline{C}\int_{\mathbb{R}} \Big(1+\fr{S^2}{\sqrt{1+S^2}}\Big) \mathcal{J}_{4}^- \mathcal{W}^-\ {\rm d}x  +\overline{C}\int_{\mathbb{R}} \mathcal{J}_{5}^- \mathcal{W}^-\ {\rm d}x +\overline{C}\int_{\mathbb{R}} \mathcal{J}_{4}^+\mathcal{W}^+\ {\rm d}x \\ & +\overline{C}\int_{\mathbb{R}}\mathcal{J}_{2}^-\mathcal{W}^- +\mathcal{J}_{2}^+\mathcal{W}^+\ {\rm d}x +\overline{C}\int_{\mathbb{R}}\Big(1+\fr{S^2}{\sqrt{1+S^2}}\Big)\mathcal{J}_{6}^-\mathcal{W}^-\ {\rm d}x \\
&+\overline{C}\int_{\mathbb{R}}  \Big[\Big(1+\fr{S^2}{\sqrt{1+S^2}}\Big)\mathcal{J}_{1}^-+\mathcal{J}_{1}^+\Big] \mathcal{W}^-\ {\rm d}x  +\int_{\mathbb{R}}\mathcal{J}_{6}^- \Big(\overline{C}\eta-\fr{2S^2}{\overline{C}\sqrt{1+S^2}}\Big)\ {\rm d}x \\
\leq & \overline{C}\int_{\mathbb{R}}\mathcal{J}_{6}^-\mathcal{W}^-\ {\rm d}x +\overline{C}\int_{\mathbb{R}} \mathcal{J}_{5}^- \mathcal{W}^-\ {\rm d}x + \overline{C} \mathcal{I}_4 +\overline{C} \mathcal{I}_2 +\overline{C} \mathcal{I}_1 \\
&+ \overline{C}\eta\int_{\mathbb{R}} \fr{S^2}{\sqrt{1+S^2}} \mathcal{J}_{4}^-\ {\rm d}x  + \overline{C}\eta\int_{\mathbb{R}} \fr{S^2}{\sqrt{1+S^2}} \mathcal{J}_{1}^-\ {\rm d}x -\fr{1-\overline{C}\eta}{\overline{C}}\int_{\mathbb{R}} \fr{S^2}{\sqrt{1+S^2}}\mathcal{J}_{6}^- \ {\rm d}x.
\end{split}
\end{align}
We perform a similar argument for $\mathcal{J}_{6}^+$ as in \eqref{6.53} to gain
\begin{align}\label{6.54}
\begin{split}
\fr{{\rm d}}{{\rm d}t}\mathcal{I}_6
\leq & \overline{C}\mathcal{I}_6 +\overline{C}\mathcal{I}_5 + \overline{C} \mathcal{I}_4 +\overline{C} \mathcal{I}_2 +\overline{C} \mathcal{I}_1+ \overline{C}\eta\int_{\mathbb{R}} \fr{S^2}{\sqrt{1+S^2}} \mathcal{J}_{4}^-\ {\rm d}x \\
&+ \overline{C}\eta\int_{\mathbb{R}} \fr{R^2}{\sqrt{1+R^2}} \mathcal{J}_{4}^+\ {\rm d}x+ \overline{C}\eta\int_{\mathbb{R}} \fr{S^2}{\sqrt{1+S^2}} \mathcal{J}_{1}^-\ {\rm d}x + \overline{C}\eta\int_{\mathbb{R}} \fr{R^2}{\sqrt{1+R^2}} \mathcal{J}_{1}^+\ {\rm d}x \\ & -\fr{1-\overline{C}\eta}{\overline{C}}\int_{\mathbb{R}} \fr{S^2}{\sqrt{1+S^2}}\mathcal{J}_{6}^- \ {\rm d}x -\fr{1-\overline{C}\eta}{\overline{C}}\int_{\mathbb{R}} \fr{R^2}{\sqrt{1+R^2}}\mathcal{J}_{6}^+ \ {\rm d}x.
\end{split}
\end{align}

\textbf{Step 7 (Integration).}
We first combine \eqref{6.17a}, \eqref{6.29}, \eqref{6.35}, \eqref{6.35e} and \eqref{6.41} to obtain
\begin{align}\label{6.55a}
\begin{split}
&\fr{{\rm d}}{{\rm d}t}\|(\upsilon, \mu, w, \tilde{r}, z, \tilde{s})(t)\|_{(\tau, u, R, S)(t)}^{(i)} =\sum_{i=\ell}^4\fr{{\rm d}}{{\rm d}t}\mathcal{I}_\ell \\
\leq & \overline{C}_1\sum_{\ell=1}^4\mathcal{I}_\ell -\fr{1-\overline{C}_1\eta}{\overline{C}_1}\int_{\mathbb{R}}\fr{R^2}{\sqrt{1+R^2}} \Big(\mathcal{J}_{1}^+  +\mathcal{J}_{4}^+ \Big)\ {\rm d}x \\ &-\fr{1-\overline{C}_1\eta}{\overline{C}_1}\int_{\mathbb{R}}\fr{S^2}{\sqrt{1+S^2}} \Big(\mathcal{J}_{1}^-  +\mathcal{J}_{4}^- \Big)\ {\rm d}x \\
\leq &\overline{C}_1\sum_{\ell=1}^4\mathcal{I}_\ell =\overline{C}_1\|(\upsilon, \mu, w, \tilde{r}, z, \tilde{s})(t)\|_{(\tau, u, R, S)(t)}^{(i)},
\end{split}
\end{align}
for some positive constant $\overline{C}_1$, provided that $\overline{C}_1\eta\leq 1$. This finishes the proof of \eqref{6.23a}.

Next, combining \eqref{6.17}, \eqref{6.29}, \eqref{6.35}, \eqref{6.35e}, \eqref{6.41}, \eqref{6.47} and \eqref{6.54} yields
\begin{align}\label{6.55}
\begin{split}
&\fr{{\rm d}}{{\rm d}t}\|(\upsilon, \mu, w, \tilde{r}, z, \tilde{s})(t)\|_{(\tau, u, R, S)(t)}^{(ii)} =\sum_{\ell=1}^6\fr{{\rm d}}{{\rm d}t}\mathcal{I}_\ell \\
\leq & \overline{C}_2\sum_{i=1}^6\mathcal{I}_\ell -\fr{1-\overline{C}_2\eta}{\overline{C}_2}\int_{\mathbb{R}}\fr{R^2}{\sqrt{1+R^2}} \Big(\mathcal{J}_{1}^+ +\mathcal{J}_{4}^+ +\mathcal{J}_{6}^+\Big)\ {\rm d}x \\ &-\fr{1-\overline{C}_2\eta}{\overline{C}_2}\int_{\mathbb{R}}\fr{S^2}{\sqrt{1+S^2}} \Big(\mathcal{J}_{1}^- +\mathcal{J}_{4}^- +\mathcal{J}_{6}^-\Big)\ {\rm d}x \\
\leq &\overline{C}_2\sum_{\ell=1}^6\mathcal{I}_\ell =\overline{C}_2\|(\upsilon, \mu, w, \tilde{r}, z, \tilde{s})(t)\|_{(\tau, u, R, S)(t)}^{(ii)},
\end{split}
\end{align}
for some positive constant $\overline{C}_2$, provided that $\overline{C}_2\eta\leq 1$. This proves \eqref{6.23}. The proof of the lemma is complete.
\end{proof}

\section{Metric for the general case}\label{S7}

In this section, we first extend the Lipschitz metrics for smooth solutions in Section \ref{S6} to piecewise smooth solutions with only generic singularities, and subsequently carry out the extension to the general case.

\subsection{Tangent vectors in transformed coordinates}\label{S71}

Let $(\tau(x,t), u(x,t))$ be a reference solution of \eqref{1.1} and $(\tau^\eps(x,t), u^\eps(x,t))$ be a family of perturbed solutions satisfying $t_{\eps}^*\geq t^*$. This can be guaranteed that when we consider unidirectional perturbations. In fact, it is sufficient for the initial data of the perturbation to satisfy $R^\eps(x,0)\geq R(x,0)$ and $S^\eps(x,0)\geq S(x,0)$, where $R^\eps(x,0)=u_{x}^\eps(x,0)+c(\tau^\eps(x,0))\tau_{x}^\eps(x,0)$ and $S^\eps(x,0)=u_{x}^\eps(x,0)-c(\tau^\eps(x,0))\tau_{x}^\eps(x,0)$.

In the $(X,Y)$ plane, let $(\tau, u, \alpha, \beta, p, q, x, t)$ and $(\tau^\eps, u^\eps, \alpha^\eps, \beta^\eps, p^\eps, q^\eps, x^\eps, t^\eps)$ be the corresponding smooth solutions of \eqref{3.14a} and \eqref{3.16} defined on $\Omega_{t^*}$ and $\Omega_{t_{\eps}^*}$. Obviously, there holds $\Omega_{t^*}\subseteq \Omega_{t_{\eps}^*}$. For each $\eps$, the curves given by $X=Const.$ and $Y=Const.$ correspond, respectively, to the backward and forward characteristics of the solutions $(\tau^\eps(x,t), u^\eps(x,t))$. At time
$t=0$ we choose the parameterizations as
\begin{align}\label{7.1}
\begin{split}
X^\eps(0, x+\eps w(x,0))=x,\qquad Y^\eps(0, x+\eps z(x,0))=-x,
\end{split}
\end{align}
where $w(x,0)=w_0(x)$ and $z(x,0)=z_0(x)$ are the shifts in \eqref{6.8}. For $\varsigma\leq t^*$, consider the curve in $(X,Y)$ plane
\begin{align}\label{7.2}
\begin{split}
\Gamma_{\varsigma} =&\{(X,Y),\ t(X,Y)=\varsigma\} \\
=& \{(X, Y(\varsigma, X));\ X\in \mathbb{R}\} =\{(X(\varsigma, Y), Y);\ Y\in \mathbb{R}\},
\end{split}
\end{align}
and denote the corresponding perturbed curve by
\begin{align}\label{7.3}
\begin{split}
\Gamma_{\varsigma}^\eps =&\{(X,Y),\ t^\eps (X,Y)=\varsigma\} \\
=& \{(X, Y^\eps(\varsigma, X));\ X\in \mathbb{R}\} =\{(X^\eps(\varsigma, Y), Y);\ Y\in \mathbb{R}\}.
\end{split}
\end{align}
Assume that the perturbed solutions take the form
\begin{align}\label{7.4}
\begin{split}
&(\tau^\eps, u^\eps, \alpha^\eps, \beta^\eps, p^\eps, q^\eps, x^\eps, t^\eps)\\
=& (\tau, u, \alpha, \beta, p, q, x, t)+\eps (\Xi, \mathcal{U}, \mathcal{A}, \mathcal{B}, \mathcal{P}, \mathcal{Q}, \mathcal{X}, \mathcal{T}) +o(\eps).
\end{split}
\end{align}
It is noted by \eqref{6.1b} that the coefficients of system \eqref{3.14a} and \eqref{3.16} are smooth under the condition \eqref{6.1a} for sufficiently small $\eta$. Then one can find that the first
order perturbations $(\Xi, \mathcal{A}, \mathcal{B}, \mathcal{P}, \mathcal{Q}, \mathcal{X}, \mathcal{T})$ satisfy a linearized system and are well defined for $(X,Y)\in \Omega_{t^*}$.

We next use the perturbation variables $(\Xi, \mathcal{U}, \mathcal{A}, \mathcal{B}, \mathcal{P}, \mathcal{Q}, \mathcal{X}, \mathcal{T})$ to express the quantities $(\upsilon, \mu, \tilde{r}, w, \tilde{s}, z)$ defined in Section \ref{S6}. From the definitions
$$
t^\eps(X, Y^\eps(\varsigma, X))=\varsigma=t^\eps(X^\eps(\varsigma, Y), Y),
$$
we see by the implicit function theorem and \eqref{3.16} that, at $\eps=0$
\begin{align*}
\begin{split}
\fr{\pa X^\eps}{\pa \eps}\Big|_{\eps=0}=&-\fr{\pa t^\eps}{\pa \eps}\cdot (t_X)^{-1}=-\mathcal{T}\fr{2c}{p\cos\fr{\alpha}{2}},\\
\fr{\pa Y^\eps}{\pa \eps}\Big|_{\eps=0}=&-\fr{\pa t^\eps}{\pa \eps}\cdot (t_Y)^{-1}=-\mathcal{T}\fr{2c}{q\cos\fr{\beta}{2}}.
\end{split}
\end{align*}

(1) The change in $x$ can be calculated by
\begin{align}\label{7.5}
\begin{split}
w=&\lim_{\eps\rightarrow0}\fr{x^\eps(X,Y^\eps(\varsigma, X))-x(X,Y(\varsigma, X))}{\eps} \\
=&\mathcal{X}(X, Y(\varsigma, X))+x_Y\cdot\fr{\pa Y^\eps}{\pa \eps}\Big|_{\eps=0} =(\mathcal{X}+c\mathcal{T})(X, Y(\varsigma, X)),
\end{split}
\end{align}
and similarly,
\begin{align}\label{7.6}
\begin{split}
z=&\lim_{\eps\rightarrow0}\fr{x^\eps(X^\eps(\varsigma, Y),Y)-x(X(\varsigma, Y),Y)}{\eps} =(\mathcal{X}-c\mathcal{T})(X(\varsigma, Y), Y).
\end{split}
\end{align}

(2) For the change in $\tau$, one obtains by \eqref{3.14a}
\begin{align*}
\begin{split}
\upsilon +\tau_xw&=\lim_{\eps\rightarrow0} \fr{\tau^\eps(X,Y^\eps(\varsigma, X))-\tau(X,Y(\varsigma, X))}{\eps} \\
&=\Xi(X,Y(\varsigma, X))+\tau_Y\cdot \fr{\pa Y^\eps}{\pa \eps}\Big|_{\eps=0} =\Big(\Xi-\mathcal{T}\tan\fr{\beta}{2}\Big)(X,Y(\varsigma, X)),
\end{split}
\end{align*}
from which we see that
\begin{align}\label{7.7}
\begin{split}
\upsilon +\fr{Rw}{2c}-\fr{Sz}{2c}=\upsilon +\tau_x w+\fr{S(w-z)}{2c} =\Xi.
\end{split}
\end{align}

(3) For the change in $u$, a similar procedure as above yields
\begin{align*}
\begin{split}
\mu +u_xw&=\lim_{\eps\rightarrow0} \fr{u^\eps(X,Y^\eps(\varsigma, X))-u(X,Y(\varsigma, X))}{\eps} \\
&=\Big(U-c\mathcal{T}\tan\fr{\alpha}{2}\Big)(X,Y(\varsigma, X)),
\end{split}
\end{align*}
subsequently
\begin{align}\label{7.7a}
\begin{split}
\mu+\fr{Rw+Sz}{2} =\mu+u_xw -\fr{S(w-z)}{2}= \mathcal{U}.
\end{split}
\end{align}

(4) To derive $(\tilde{r}, \tilde{s})$, we first have
\begin{align*}
\begin{split}
r+wR_x=&\fr{{\rm d}}{{\rm d}\eps}\tan\fr{\alpha^\eps(X,Y^\eps(\varsigma, X))}{2}\Big|_{\eps=0} =\fr{1}{2}\Big(\mathcal{A}-\mathcal{T}\fr{2c}{q\cos\fr{\beta}{2}}\alpha_Y\Big)\sec^2\fr{\alpha}{2}, \\
s+zS_x=&\fr{{\rm d}}{{\rm d}\eps}\tan\fr{\beta^\eps(X^\eps(\varsigma, Y),Y)}{2}\Big|_{\eps=0} =\fr{1}{2}\Big(\mathcal{B}-\mathcal{T}\fr{2c}{p\cos\fr{\alpha}{2}}\beta_X\Big)\sec^2\fr{\beta}{2},
\end{split}
\end{align*}
which together with \eqref{6.16} and \eqref{7.5}-\eqref{7.6} gives
\begin{align}\label{7.8}
\begin{split}
\tilde{r}=&\fr{1}{2} \Big(\mathcal{A}-\mathcal{T}\fr{2c}{q\cos\fr{\beta}{2}}\alpha_Y\Big)\sec^2\fr{\alpha}{2} -\fr{c'}{2c}\mathcal{T}\tan\fr{\alpha}{2}\tan\fr{\beta}{2}, \\
\tilde{s}=&\fr{1}{2} \Big(\mathcal{B}-\mathcal{T}\fr{2c}{p\cos\fr{\alpha}{2}}\beta_X\Big)\sec^2\fr{\beta}{2} -\fr{c'}{2c}\mathcal{T}\tan\fr{\alpha}{2}\tan\fr{\beta}{2}.
\end{split}
\end{align}

(5) We notice that
\begin{align*}
\begin{split}
r+R_xw+Rw_x=&\tilde{r} +Rw_x+\fr{c'}{4c^2}RS(w-z),\\
s+S_xz+Sz_x=&\tilde{s} +Sz_x+\fr{c'}{4c^2}RS(w-z).
\end{split}
\end{align*}
The terms $w_x$ and $z_x$ correspond to the change of base measure with density 1.
Note that along $\Gamma_\varsigma$
\begin{align}\label{7.9}
\sqrt{1+R^2}{\rm d} x=p{\rm d} X,\qquad \sqrt{1+S^2}{\rm d} x=-q{\rm d} Y.
\end{align}
Then one can compute the change in base measure with density 1 as
\begin{align}\label{7.10}
\begin{split}
&\lim_{\eps\rightarrow0}\fr{p^\eps(X,Y^\eps(\varsigma, X))\cos\fr{\alpha^\eps(X,Y^\eps(\varsigma, X))}{2}-p(X,Y(\varsigma, X))\cos\fr{\alpha(X,Y(\varsigma, X))}{2}}{\eps} \\
=&\Big(P+p_Y\fr{\pa Y^\eps}{\pa \eps}\Big|_{\eps=0}\Big)\cos\fr{\alpha}{2} -\fr{1}{2}p\sin\fr{\alpha}{2}\Big(\mathcal{A}+\alpha_Y\fr{\pa Y^\eps}{\pa \eps}\Big|_{\eps=0}\Big) \\
=&\Big(P-\mathcal{T}\fr{2c}{q\cos\fr{\beta}{2}}p_Y\Big)\cos\fr{\alpha}{2} -\fr{1}{2}p\sin\fr{\alpha}{2}\Big(\mathcal{A}-\mathcal{T}\fr{2c}{q\cos\fr{\beta}{2}}\alpha_Y\Big).
\end{split}
\end{align}
Analogously, there holds
\begin{align}\label{7.11}
\begin{split}
&\lim_{\eps\rightarrow0}\fr{q^\eps(X^\eps(\varsigma, Y),Y)\cos\fr{\beta^\eps(X^\eps(\varsigma, Y),Y)}{2}-q(X(\varsigma, Y),Y)\cos\fr{\beta(X(\varsigma, Y),Y)}{2}}{\eps} \\
=&\Big(Q-\mathcal{T}\fr{2c}{p\cos\fr{\alpha}{2}}q_X\Big)\cos\fr{\beta}{2} -\fr{1}{2}q\sin\fr{\beta}{2}\Big(\mathcal{B}-\mathcal{T}\fr{2c}{p\cos\fr{\alpha}{2}}\beta_X\Big).
\end{split}
\end{align}

Combining \eqref{7.5}-\eqref{7.11}, the weighted norms \eqref{6.17a} and \eqref{6.17} can be written as a line integral over the line $\Gamma_\varsigma$ as follows:
\begin{align}\label{7.12}
\|(\upsilon, \mu, w, \tilde{r}, z, \tilde{s})\|_{(\tau, u, R, S)}^{(i)} =\sum_{\ell=1}^4\widetilde{\mathcal{I}}_\ell,
\end{align}
and
\begin{align}\label{7.13}
\|(\upsilon, \mu, w, \tilde{r}, z, \tilde{s})\|_{(\tau, u, R, S)}^{(ii)} =\sum_{\ell=1}^6\widetilde{\mathcal{I}}_\ell,
\end{align}
where
\begin{align}\label{7.14}
\begin{split}
\widetilde{\mathcal{I}}_\ell=\int_{\Gamma_\varsigma} \Big\{|\widetilde{\mathcal{J}}_{\ell}^-|\mathcal{W}^-\ {\rm d}X +|\widetilde{\mathcal{J}}_{\ell}^+|\mathcal{W}^+\ {\rm d}Y\Big\}.
\end{split}
\end{align}
The expressions of $\widetilde{\mathcal{J}}_{i}^-$ and $\widetilde{\mathcal{J}}_{i}^+$ are
\begin{align}\label{7.15}
\begin{split}
\widetilde{\mathcal{J}}_{1}^-=&p(\mathcal{X}+c\mathcal{T}),\qquad
\widetilde{\mathcal{J}}_{1}^+=q(\mathcal{X}-c\mathcal{T}), \\
\widetilde{\mathcal{J}}_{2}^-=&p\Xi, \qquad \widetilde{\mathcal{J}}_{2}^+=q\Xi, \qquad
\widetilde{\mathcal{J}}_{3}^-=p\mathcal{U}, \qquad \widetilde{\mathcal{J}}_{3}^+=q\mathcal{U}, \\
\widetilde{\mathcal{J}}_{4}^-=& p\cos\fr{\alpha}{2}\Big[\fr{1}{2} \Big(\mathcal{A}-\mathcal{T}\fr{2c}{q\cos\fr{\beta}{2}}\alpha_Y\Big)\sec^2\fr{\alpha}{2} -\fr{c'}{2c}\mathcal{T}\tan\fr{\alpha}{2}\tan\fr{\beta}{2}\Big] \\ &+\tan\fr{\alpha}{2}\Big[\Big(P-\mathcal{T}\fr{2c}{q\cos\fr{\beta}{2}}p_Y\Big)\cos\fr{\alpha}{2} -\fr{1}{2}p\sin\fr{\alpha}{2} \Big(\mathcal{A}-\mathcal{T}\fr{2c}{q\cos\fr{\beta}{2}}\alpha_Y\Big)\Big] \\
&+\fr{c'}{4c^2}\tan\fr{\alpha}{2}\tan\fr{\beta}{2}\cdot 2c\mathcal{T}\cdot p\cos\fr{\alpha}{2},  \\
\widetilde{\mathcal{J}}_{4}^+=& q\cos\fr{\beta}{2}\Big[\fr{1}{2} \Big(\mathcal{B}-\mathcal{T}\fr{2c}{q\cos\fr{\alpha}{2}}\beta_X\Big)\sec^2\fr{\beta}{2} -\fr{c'}{2c}\mathcal{T}\tan\fr{\beta}{2}\tan\fr{\alpha}{2}\Big] \\ &+\tan\fr{\beta}{2}\Big[\Big(Q-\mathcal{T}\fr{2c}{q\cos\fr{\alpha}{2}}q_X\Big)\cos\fr{\beta}{2} -\fr{1}{2}q\sin\fr{\beta}{2} \Big(\mathcal{B}-\mathcal{T}\fr{2c}{q\cos\fr{\alpha}{2}}\beta_X\Big)\Big] \\ &+\fr{c'}{4c^2}\tan\fr{\alpha}{2}\tan\fr{\beta}{2}\cdot 2c\mathcal{T}\cdot q\cos\fr{\beta}{2},
\end{split}
\end{align}
and
\begin{align}\label{7.16}
\begin{split}
\widetilde{\mathcal{J}}_{5}^-=& p\cos^2\fr{\alpha}{2}\Big\{\fr{1}{2} \Big(\mathcal{A}-\mathcal{T}\fr{2c}{q\cos\fr{\beta}{2}}\alpha_Y\Big)\sec^2\fr{\alpha}{2} -\fr{c'}{2c}\mathcal{T}\tan\fr{\alpha}{2}\tan\fr{\beta}{2}\Big\}, \\
\widetilde{\mathcal{J}}_{5}^+=& q\cos^2\fr{\beta}{2}\Big\{\fr{1}{2} \Big(\mathcal{B}-\mathcal{T}\fr{2c}{p\cos\fr{\alpha}{2}}\beta_X\Big)\sec^2\fr{\beta}{2} -\fr{c'}{2c}\mathcal{T}\tan\fr{\alpha}{2}\tan\fr{\beta}{2}\Big\}, \\
\widetilde{\mathcal{J}}_{6}^-=&\fr{1}{2}p\sin\alpha\Big[\fr{1}{2} \Big(\mathcal{A}-\mathcal{T}\fr{2c}{q\cos\fr{\beta}{2}}\alpha_Y\Big)\sec^2\fr{\alpha}{2} -\fr{c'}{2c}\mathcal{T}\tan\fr{\alpha}{2}\tan\fr{\beta}{2}\Big] \\ &+\fr{1}{\cos\fr{\alpha}{2}}\Big[\Big(P-\mathcal{T}\fr{2c}{q\cos\fr{\beta}{2}}p_Y\Big)\cos\fr{\alpha}{2} -\fr{1}{2}p\sin\fr{\alpha}{2}\Big(\mathcal{A}-\mathcal{T}\fr{2c}{q\cos\fr{\beta}{2}}\alpha_Y\Big) \Big] \\
&+\fr{c'}{4c^2}\cdot 2c\mathcal{T}\cdot p\tan\fr{\beta}{2}, \\
\widetilde{\mathcal{J}}_{6}^+=&\fr{1}{2}q\sin\beta\Big[\fr{1}{2} \Big(\mathcal{B}-\mathcal{T}\fr{2c}{p\cos\fr{\alpha}{2}}\beta_X\Big)\sec^2\fr{\beta}{2} -\fr{c'}{2c}\mathcal{T}\tan\fr{\alpha}{2}\tan\fr{\beta}{2}\Big] \\ &+\fr{1}{\cos\fr{\beta}{2}}\Big[\Big(Q-\mathcal{T}\fr{2c}{p\cos\fr{\alpha}{2}}q_X\Big)\cos\fr{\beta}{2} -\fr{1}{2}q\sin\fr{\beta}{2}\Big(\mathcal{B}-\mathcal{T}\fr{2c}{p\cos\fr{\alpha}{2}}\beta_X\Big) \Big] \\
&+\fr{c'}{4c^2}\cdot 2c\mathcal{T}\cdot q\tan\fr{\alpha}{2}.
\end{split}
\end{align}
Making use of \eqref{3.14}, we can simplify as the expressions in \eqref{7.15} and \eqref{7.16} as
\begin{align}\label{7.17}
\begin{split}
\widetilde{\mathcal{J}}_{1}^-=&p(\mathcal{X}+c\mathcal{T}),\qquad  \qquad  \qquad  \qquad
\widetilde{\mathcal{J}}_{1}^+=q(\mathcal{X}-c\mathcal{T}), \\
\widetilde{\mathcal{J}}_{2}^-=&p\Xi, \qquad \widetilde{\mathcal{J}}_{2}^+=q\Xi, \qquad \qquad \quad
\widetilde{\mathcal{J}}_{3}^-=p\mathcal{U}, \qquad \widetilde{\mathcal{J}}_{3}^+=q\mathcal{U}, \\
\widetilde{\mathcal{J}}_{4}^-=& P\sin\fr{\alpha}{2} +\fr{1}{2}p\cos\fr{\alpha}{2}\mathcal{A}, \qquad \quad \ \
\widetilde{\mathcal{J}}_{4}^+=  Q\sin\fr{\beta}{2} +\fr{1}{2}q\cos\fr{\beta}{2}\mathcal{B}, \\
\widetilde{\mathcal{J}}_{5}^-=& \fr{1}{2}p\mathcal{A}-\fr{c'}{2c}\mathcal{T}p\sin^2\fr{\alpha}{2}, \qquad \qquad
\widetilde{\mathcal{J}}_{5}^+= \fr{1}{2}q\mathcal{B}-\fr{c'}{2c}\mathcal{T}q\sin^2\fr{\beta}{2}, \\
\widetilde{\mathcal{J}}_{6}^-=& P +\fr{c'}{4c}p\mathcal{T}\sin\alpha,  \qquad \qquad \qquad
\widetilde{\mathcal{J}}_{6}^+= Q +\fr{c'}{4c}q\mathcal{T}\sin\beta.
\end{split}
\end{align}
It is obvious from \eqref{7.17} that the integrands $\widetilde{\mathcal{J}}_{\ell}^-$ and $\widetilde{\mathcal{J}}_{\ell}^+$ are smooth for $\ell=1,\cdots,6$. Note that the smoothness of these terms permits us to pass to the limit for the metrics defined later.

\subsection{Length of piecewise regular paths}\label{S72}

Thanks to Theorems \ref{thm1} and \ref{thm2}, we know that, for generic initial data $(\tau_0(x), u_0(x))$ satisfying $(\tau_0(x)-\bar{\tau}_0, u_0(x))\in \mathcal{D}$ and \eqref{4.34}, the solution $(\tau(x,t), u(x,t))$ of \eqref{1.1} is piecewise smooth and has only generic singularities on the line $t=t^*$. In what follows, we give the definition of piecewise regular paths.

\begin{defn}\label{def2}
We say a path of solution $\Upsilon^t: \vartheta\mapsto (\tau^\vartheta, u^\vartheta)$ $(\vartheta\in[0,1])$ is piecewise regular for $t\in[0,T]$, if the following hold:

{\rm (i)} For each $\vartheta\in[0,1]$, the solution $(\tau^\vartheta, u^\vartheta)$ of \eqref{1.1} is smooth for $t\in[0,T)$.

{\rm (ii)} There exist finitely many values $\vartheta_i\ (i=1,\cdots, N)$ such that for $\vartheta\in[0,1]\setminus\{\vartheta_1,\cdots, \vartheta_N\}$, the solution $(\tau^\vartheta, u^\vartheta)$ has only generic singularities, if they exist, on the line $t=T$.
\end{defn}

To achieve our objective, we present below the result concerning the density for the set of piecewise regular paths, which follows directly from Theorem \ref{thm1}.
\begin{lem}\label{lem7}
Let $T>0$. Assume that $\vartheta\mapsto(\tau^\vartheta, u^\vartheta, \alpha^\vartheta, \beta^\vartheta, p^\vartheta, q^\vartheta, x^\vartheta, t^\vartheta)$ is a smooth path of solutions to \eqref{3.14a}, \eqref{3.16}. Then there exists a sequence of paths of solutions $\vartheta\mapsto(\tau_{n}^\vartheta, u_{n}^\vartheta, \alpha_{n}^\vartheta, \beta_{n}^\vartheta, p_{n}^\vartheta, q_{n}^\vartheta, x_{n}^\vartheta, t_{n}^\vartheta)$ such that

{\rm (i)} For each $n\geq1$, the path of the corresponding solution of \eqref{1.1} $\vartheta\mapsto (\tau_{n}^\vartheta, u_{n}^\vartheta, )$ is piecewise regular for $t\in[0,T]$ in the sense of Definition \ref{def2}.

{\rm (ii)} For the bounded domain $\Omega_L=\Omega_T\cap \{(X,Y):\ |X|\leq L, |Y|\leq L\}$, where $L$ is any positive number, the functions $(\tau_{n}^\vartheta, u_{n}^\vartheta, \alpha_{n}^\vartheta, \beta_{n}^\vartheta, p_{n}^\vartheta, q_{n}^\vartheta, x_{n}^\vartheta, t_{n}^\vartheta)$ converge to $(\tau^\vartheta, u^\vartheta, \alpha^\vartheta, \beta^\vartheta, p^\vartheta, q^\vartheta, x^\vartheta, t^\vartheta)$ uniformly in $C^k([0,1]\times\Omega_L)$, for every $k\geq1$.
\end{lem}
\begin{proof}
The proof is similar to that of Theorem 2 in Bressan and Chen \cite{BC1} (also see Theorem 4 in \cite{BC2}), we omit it here.
\end{proof}

By virtue of the density conclusion given in Lemma \ref{lem7}, constructing a Lipschitz metric reduces to show that the weighted length of a regular path satisfies the same estimates as the smooth paths in Section \ref{S6}. We now define the length of a piecewise regular path $\Upsilon^t: \vartheta\mapsto (\tau^\vartheta, u^\vartheta)$.

\begin{defn}\label{def3}
Let $T>0$ be given and $\Upsilon^t: \vartheta\mapsto (\tau^\vartheta, u^\vartheta)\ (\vartheta\in[0,1])$ be a piecewise regular path for $t\in[0, T]$ in the sense of Definition \ref{def2}. The two classis of lengths for $\Upsilon^t$, denoted by $\|\Upsilon^t\|^{(i)}$ and $\|\Upsilon^t\|^{(ii)}$, are defined as
\begin{align}\label{7.19}
\begin{split}
\|\Upsilon^t\|^{(i)}=\inf_{\Upsilon^t}\int_{0}^1 \Big\{\sum_{\ell=1}^4 \int_{\Gamma_{t}^\vartheta}\Big(|(\widetilde{\mathcal{J}}_{\ell}^-)^\vartheta|(\mathcal{W}^-)^\vartheta\ {\rm d}X +|(\widetilde{\mathcal{J}}_{\ell}^+)^\vartheta|(\mathcal{W}^+)^\vartheta\ {\rm d}Y\Big)\Big\}\ {\rm d}\vartheta,
\end{split}
\end{align}
and
\begin{align}\label{7.20}
\begin{split}
\|\Upsilon^t\|^{(ii)}=\inf_{\Upsilon^t}\int_{0}^1 \Big\{\sum_{\ell=1}^6 \int_{\Gamma_{t}^\vartheta}\Big(|(\widetilde{\mathcal{J}}_{\ell}^-)^\vartheta|(\mathcal{W}^-)^\vartheta\ {\rm d}X +|(\widetilde{\mathcal{J}}_{\ell}^+)^\vartheta|(\mathcal{W}^+)^\vartheta\ {\rm d}Y\Big)\Big\}\ {\rm d}\vartheta,
\end{split}
\end{align}
where the infimum is taken over all piecewise smooth relabelings of the $X$-$Y$ coordinates, and $\Gamma_{\varsigma}^\vartheta:=\{(X,Y);\ t^\vartheta(X,Y)=\varsigma\}$.
\end{defn}
\begin{rem}\label{r5}
According to Remark \ref{r4} in Section \ref{S3}, there exist infinitely many paths of solutions of \eqref{3.14a}, \eqref{3.16} yielding the same path of solutions to \eqref{2.9}-\eqref{2.9a}. The values of the  integrands in \eqref{7.19} and \eqref{7.20} depend on the choice of relabeling for the $(X,Y)$ coordinates. Accordingly, we take the infimum over all of these relabelings.
\end{rem}

Then we have
\begin{thm}\label{thm4}
Let $T>0$ be given and $\Upsilon^t: \vartheta\mapsto (\tau^\vartheta, u^\vartheta)\ (\vartheta\in[0,1])$ be a piecewise regular path for $t\in[0, T]$ in the sense of Definition \ref{def2}. For each $\theta\in[0,1]$, we further assume that the solution $(\tau^\vartheta(x,t), u^\vartheta(x,t))$ satisfies
\begin{align}\label{7.21}
\begin{split}
\int_{\mathbb{R}}|R^\vartheta(x, t)|+|S^\vartheta(x, t)|\ {\rm d}x\leq  \eta,
\end{split}
\end{align}
for some constant $\eta>0$, where $R^\vartheta=u_{x}^\vartheta+c(\tau^\vartheta)\tau_{x}^\vartheta$ and $S^\vartheta=u_{x}^\vartheta-c(\tau^\vartheta)\tau_{x}^\vartheta$.
Then the lengths satisfy
\begin{align}\label{7.22}
\begin{split}
\|\Upsilon^t\|^{(i)}\leq e^{\overline{C}_1T} \|\Upsilon^0\|^{(i)},
\end{split}
\end{align}
provided that $\eta\leq 1/\overline{C}_1$, and
\begin{align}\label{7.23}
\begin{split}
\|\Upsilon^t\|^{(ii)}\leq e^{\overline{C}_2T} \|\Upsilon^0\|^{(ii)},
\end{split}
\end{align}
provided that $\eta\leq 1/\overline{C}_2$. Here the constants $\overline{C}_{1}$ and $\overline{C}_{2}$ are given in \eqref{6.22a} and \eqref{6.22}, respectively.
\end{thm}
\begin{proof}
To show \eqref{7.22} and \eqref{7.23}, it is sufficient to establish the following inequalities
\begin{align}\label{7.24}
\begin{split}
&\|(\upsilon^\vartheta, \mu^\vartheta, r^\vartheta, s^\vartheta)(t)\|_{(\tau^\vartheta, u^\vartheta, R^\vartheta, S^\vartheta)(t)}^{(\iota)} \\
\leq& e^{\overline{C}T}\|(\upsilon^\vartheta, \mu^\vartheta, r^\vartheta, s^\vartheta)(0)\|_{(\tau^\vartheta, u^\vartheta, R^\vartheta, S^\vartheta)(0)}^{(\iota)},\quad (\iota=i, ii)
\end{split}
\end{align}
for $\vartheta\in[0,1]\setminus\{\vartheta_1,\cdots, \vartheta_N\}$. In fact, by means of Definition \ref{def3}, we see that for any fixed $\eps>0$, there exists a relabeling of variables $(X,Y)$ such that there holds at time $t=0$
\begin{align}\label{7.25}
\begin{split}
\int_{0}^1 \Big\{\sum_{\ell} \int_{\Gamma_{0}^\vartheta}\Big(|(\widetilde{\mathcal{J}}_{\ell}^-)^\vartheta|(\mathcal{W}^-)^\vartheta\ {\rm d}X +|(\widetilde{\mathcal{J}}_{\ell}^+)^\vartheta|(\mathcal{W}^+)^\vartheta\ {\rm d}Y\Big)\Big\}\ {\rm d}\vartheta \leq \|\Upsilon^0\|^{(\iota)}+\eps.
\end{split}
\end{align}
One integrates \eqref{7.24} with respect to $\vartheta$ and utilizes \eqref{7.25} to obtain
$$
\|\Upsilon^t\|^{(\iota)}\leq e^{\overline{C}T}\Big(\|\Upsilon^0\|^{(\iota)}+\eps\Big),
$$
which together with the arbitrariness of $\eps$ leads directly to \eqref{7.22} and \eqref{7.23}.

Now, let $\Upsilon^t: \vartheta\mapsto (\tau^\vartheta, u^\vartheta)$ be a given piecewise regular path fulfilling \eqref{7.21} on $[0,T]$. According to Definition \ref{def2}, for every $\vartheta\in[0,1]\setminus\{\vartheta_1,\cdots, \vartheta_N\}$, the solution $(\tau^\vartheta, u^\vartheta )$ is smooth for $t\in[0, T)$ and has only generic regularities, if they exist, on $t=T$. By Lemma \ref{lem2}, we know that the solution $(\tau^\vartheta, u^\vartheta, \alpha^\vartheta, \beta^\vartheta, p^\vartheta, q^\vartheta, x^\vartheta, t^\vartheta)$ of \eqref{3.14a}, \eqref{3.16} is smooth on $\Omega_T$ in the $(X,Y)$ plane. Thus the tangent vector is well-defined for all $\vartheta\in[0,1]$ and $t\in[0,T]$. In what follows, we aim at establishing \eqref{7.24}.

If the solution $(\tau^\vartheta, u^\vartheta)$ is smooth up to $t=T$, then one can obtain the estimate \eqref{7.24} directly from \eqref{6.22a} and \eqref{6.22}.

For the case that $(\tau^\vartheta, u^\vartheta)$ is smooth in $[0,T)$ and has generic singularities on $t=T$, we next demonstrate that the onset of these generic singularities does not affect the conclusion. Indeed, we first observe that the map
\begin{align}\label{7.26a}
\begin{split}
t\mapsto\sum_{\ell} \int_{\Gamma_{t}^\vartheta}\Big(|(\widetilde{\mathcal{J}}_{\ell}^-)^\vartheta|(\mathcal{W}^-)^\vartheta\ {\rm d}X +|(\widetilde{\mathcal{J}}_{\ell}^+)^\vartheta|(\mathcal{W}^+)^\vartheta\ {\rm d}Y\Big),
\end{split}
\end{align}
is continuous up to $t=T$. On the other hand, the estimate \eqref{7.24} is valid for $t<T$, because the solution is smooth in this regime. Owing to the continuity of the map in \eqref{7.26a}, we can take $t\rightarrow T^-$ to obtain that \eqref{7.24} holds up to $t=T$. The proof of the theorem is complete.
\end{proof}

\subsection{Construction of the geodesic distance}\label{S73}

We are going to construct a geodesic distance $d^*(\cdot,\cdot)$ on the space $(C^1\cap W^{1,1})^2$ for the compressible Euler equations.

Let $\bar{\tau}_0, \eta$ and $\bar{c}$ be three positive constants. Denote the set
\begin{align}\label{7.26}
\mathcal{D}_{\eta}:=\left\{(\tau_0(x), u_0(x))\left|
\begin{array}{ll}
\dps (\tau_0-\bar{\tau}_0, u_0)\in \Big(C^1(\mathbb{R})\cap W^{1,1}(\mathbb{R})\Big)^2, \\[6pt]
\dps \fr{1}{\bar{c}}\leq \tau_0(x)\leq \bar{c},\\[10pt]
\dps \int_{\mathbb{R}}[|R_0(x)|+|S_0(x)|] {\rm d}x\leq \eta,
\end{array}
\right.
\right\},
\end{align}
where $R_0(x)=u_{0}'(x)+c(\tau_0(x))\tau_{0}'(x)$ and $S_0(x)=u_{0}'(x)-c(\tau_0(x))\tau_{0}'(x)$. According to Theorem \ref{thm1}, there exists an open dense subset $\mathcal{D}$ of $(C^3(\mathbb{R})\cap W^{1,1}(\mathbb{R})^2$ such that, for the initial data $(\tau_0-\bar{\tau}_0, u_0)\in \mathcal{D}$, the solution $(\tau(x,t), u(x,t)$ of \eqref{1.1}-\eqref{1.3} is smooth on $t\in[0,t^*)$ and has only generic singularities at $t=t^*$. Here $t^*$ is the first time of singularity formation for the solution $(\tau, u)$. We next define a geodesic distance by taking the infimum over all piecewise regular paths connecting two solutions of \eqref{1.1}-\eqref{1.3} on the set $\mathcal{D}^\infty:=C_{0}^\infty\cap \mathcal{D}$. Then this distance can be extended by continuity from $\mathcal{D}^\infty$ to a larger space.

\begin{defn}\label{def4}
Let $(\tau(x,t), u(x,t))$ and $(\hat{\tau}(x,t), \hat{u}(x,t))$ be the two solutions of \eqref{1.1}-\eqref{1.3} with initial data in $\mathcal{D}^\infty\cap \mathcal{D}_\eta$. Denote by $t^*$ and $\hat{t}^*$ the times at which their first singularities form, respectively. For any time $t\leq t_d:=\min\{t^*, \hat{t}^*\}$, we define the geodesic distance $d^*((\tau, u),(\hat{\tau}, \hat{u}))$ as follows:
\begin{align}\label{7.27}
&d_{*}^{(\iota)}((\tau, u),(\hat{\tau}, \hat{u})) \notag \\
:=&\inf\left\{ \|\Upsilon^t\|^{(\iota)}
\left|
\begin{array}{ll}
\dps \Upsilon^t\ {\rm is\ a\ piecewise\ regular\ path\ for}\ t\in[0,t_d]\ {\rm satisfying} \\[6pt]
\dps \Upsilon^t(0)=(\tau, u),\quad  \Upsilon^t(1)=(\hat{\tau}, \hat{u}), \\[6pt]
\dps \int_{\mathbb{R}}\big(|R^\vartheta|+|S^\vartheta|\big)\ {\rm d}x\leq \eta,\quad \forall\ \vartheta\in[0,1],
\end{array}
\right.
\right\},
\end{align}
for $\iota=i, ii$, where $\eta$ is a small positive constant ensuring positive upper and lower bounds for $\tau^\vartheta$.
\end{defn}

Since the concatenation of two piecewise regular paths remains a piecewise regular path, one concludes that $d_{*}^{(\iota)}(\cdot, \cdot)$ defined in \eqref{7.27} is indeed a metric. Based on this definition of metric, we can define the distance on the space $(C^1\cap W^{1,1})^2$.
\begin{defn}\label{def5}
Let $(\tau_0, u_0)$ and $(\hat{\tau}_0, \hat{u}_0)$ in $\mathcal{D}_\eta$ be two initial data, and $(\tau(x,t), u(x,t))$ and $(\hat{\tau}(x,t), \hat{u}(x,t))$ be the corresponding solutions of \eqref{1.1}-\eqref{1.3}. Denote by $t^*$ and $\hat{t}^*$ the times at which their first singularities form, respectively. For any time $t\leq t_d:=\min\{t^*, \hat{t}^*\}$, we define $d_{*}^{(\iota)}((\tau, u),(\hat{\tau}, \hat{u}))$ as follows:
\begin{align}\label{7.28}
d_{*}^{(\iota)}((\tau, u),(\hat{\tau}, \hat{u})) :=\lim_{n\rightarrow\infty}d_{*}^{(\iota)}((\tau^{(n)}, u^{(n)}),(\hat{\tau}^{(n)}, \hat{u}^{(n)})),
\end{align}
where $(\tau^{(n)}, u^{(n)})$ and $(\hat{\tau}^{(n)}, \hat{u}^{(n)})$ are two arbitrary sequences of solutions to \eqref{1.1}-\eqref{1.3} with the initial data $(\tau_{0}^{(n)}, u_{0}^{(n)})$ and $(\hat{\tau}_{0}^{(n)}, \hat{u}_{0}^{(n)})$ in $\mathcal{D}^\infty\cap \mathcal{D}_\eta$ satisfying
\begin{align}\label{7.29}
\|(\tau_{0}^{(n)}-\tau_0, \hat{\tau}_{0}^{(n)}-\hat{\tau}_0)\|_{C^1\cap W^{1,1}}\rightarrow0,\quad
\|(u_{0}^{(n)}-u_0, \hat{u}_{0}^{(n)}-\hat{u}_0)\|_{C^1\cap W^{1,1}}\rightarrow0.
\end{align}
Here $\|\cdot\|_{C^1\cap W^{1,1}}:= \|\cdot\|_{C^1}+\|\cdot\|_{W^{1,1}}$.
\end{defn}

Regarding to $d_{*}^{(\iota)}((\tau, u), (\hat{\tau}, \hat{u}))$, we have the following Lipschitz continuous property.
\begin{thm}\label{thm5}
With the setting from Definition \ref{def5}, for any time $t\leq t_d$,
the two solutions $(\tau(x,t), u(x,t))$ and $(\hat{\tau}(x,t), \hat{u}(x,t))$ of \eqref{1.1}-\eqref{1.3} satisfy
\begin{align}\label{7.30}
d_{*}^{(\iota)}((\tau, u),(\hat{\tau}, \hat{u}))\leq \widehat{C} d_{*}^{(\iota)}((\tau_0, u_{0}),(\hat{\tau}_0, \hat{u}_{0})),
\end{align}
for some positive constant $\widehat{C}$ depending only on $\bar{c}$ and $t_d$, provided that $\eta$ is sufficiently small.
\end{thm}
\begin{proof}
We first know that the solution with initial data in $\mathcal{D}^\infty\cap \mathcal{D}_\eta$ is Lipschitz continuous. Thus the limit in \eqref{7.28} is independent of the selection of sequences, which ensures that the definition of the metric $d_{*}^{(\iota)}(\cdot,\cdot)$ is well-defined.
Moreover, due to the fact that $\mathcal{D}^\infty\cap \mathcal{D}_\eta$ is a dense set in the space $(C^1\cap W^{1,1})^2$, we then can directly extend the Lipschitz metric to the general case by Theorem \ref{thm4}. The proof of the theorem is complete.
\end{proof}

\subsection{Comparison with other metrics}\label{S74}

Finally, we compare the distance $d_{*}^{(\iota)}(\cdot,\cdot)$ other types of metrics.

\begin{prop}\label{prop2}
Let $(\tau_0, u_0)$ and $(\hat{\tau}_0, \hat{u}_0)$ be any two initial data
in $\mathcal{D}_\eta$, and $(\tau, u)$ and $(\hat{\tau}, \hat{u})$ be the corresponding solutions of \eqref{1.1}-\eqref{1.3}. Let $t^*$ and $\hat{t}^*$ be the corresponding times at which their first singularities form. For small enough $\eta>0$ depending only on the initial data, there exists some positive constant $\widehat{C}$ such that there holds for any time $t\leq t_d:=\min\{t^*, \hat{t}^*\}$
\begin{align}\label{7.31a}
\begin{split}
&d_{*}^{(\iota)}((\tau, u), (\hat{\tau}, \hat{u})) \\
\leq & \widehat{C}\Big(\|\tau_0-\hat{\tau}_0\|_{C^1} + \|\tau_0-\hat{\tau}_0\|_{W^{1,1}} + \|u_0-\hat{u}_0\|_{C^1}+  \|u_0-\hat{u}_0\|_{W^{1,1}} \Big), \quad (\iota=i, ii).
\end{split}
\end{align}
\end{prop}
\begin{proof}
According to the Lipschitz continuous dependence with respect to the Finsler norm (i.e. Theorem \ref{thm5}), it suffices to establish the inequality \eqref{7.31a} at $t=0$, that is
\begin{align}\label{7.31}
\begin{split}
&d_{*}^{(\iota)}((\tau_0, u_0), (\hat{\tau}_0, \hat{u}_0)) \\
\leq & \widehat{C}\Big(\|\tau_0-\hat{\tau}_0\|_{C^1} + \|\tau_0-\hat{\tau}_0\|_{W^{1,1}} + \|u_0-\hat{u}_0\|_{C^1}+  \|u_0-\hat{u}_0\|_{W^{1,1}} \Big),
\end{split}
\end{align}
for some positive constant $\widehat{C}$. To prove \eqref{7.31}, from the definition in \eqref{7.27}, we only need to construct one path $\Upsilon^0$ connecting $(\tau_0, u_0)$ and $(\hat{\tau}_0, \hat{u}_0)$ such that
\begin{align}\label{7.32}
\begin{split}
\|\Upsilon^0\|^{(\iota)}
\leq & \widehat{C}\Big(\|u_0-\hat{u}_0\|_{C^1}+  \|u_0-\hat{u}_0\|_{W^{1,1}} +  \|\tau_0-\hat{\tau}_0\|_{C^1} + \|\tau_0-\hat{\tau}_0\|_{W^{1,1}}\Big),
\end{split}
\end{align}

We choose the path $\Upsilon^0$:
\begin{align}\label{7.33}
\begin{split}
u_{0}^\vartheta=(1-\vartheta)u_0 +\vartheta\hat{u}_0,\qquad \tau_{0}^\vartheta =\Big((1-\vartheta)\tau_{0}^{-\fr{\gamma-1}{2}} +\vartheta\hat{\tau}_{0}^{-\fr{\gamma-1}{2}}  \Big)^{-\fr{2}{\gamma-1}},
\end{split}
\end{align}
from which and the properties of $(\tau_0, u_0)$ and $(\hat{\tau}_0, \hat{u}_0)$,  one sees that
\begin{align}\label{7.34}
\begin{split}
&(\tau_{0}^\vartheta|_{\vartheta=0}, u_{0}^\vartheta|_{\vartheta=0})=(\tau_0, u_0), \qquad  (\tau_{0}^\vartheta|_{\vartheta=1}, u_{0}^\vartheta|_{\vartheta=1})=(\hat{\tau}_0, \hat{u}_0), \\ &
R_{0}^\vartheta =(1-\vartheta)R_0+\vartheta \hat{R}_0,\qquad \quad \ S_{0}^\vartheta =(1-\vartheta)S_0+\vartheta \hat{S}_0, \\
& \fr{1}{\bar{c}}\leq \tau_{0}^\vartheta\leq \bar{c}, \qquad \qquad \qquad \qquad   \int_{\mathbb{R}} (|R_{0}^\vartheta| +|S_{0}^\vartheta|)\ {\rm d}x\leq \eta.
\end{split}
\end{align}
Now we choose $w^\vartheta=z^\vartheta=0$ to simplify the terms $\mathcal{J}_{\ell}^\pm$ in \eqref{6.19} to
\begin{align}\label{7.35}
\begin{split}
& \mathcal{J}_{1}^-=0, \quad \mathcal{J}_{1}^+=0, \quad
\mathcal{J}_{2}^-=|\upsilon^\vartheta|\sqrt{1+(R^\vartheta)^2}, \quad \mathcal{J}_{2}^+=|\upsilon^\vartheta|\sqrt{1+(S^\vartheta)^2},  \\
& \mathcal{J}_{3}^-=|\mu^\vartheta |\sqrt{1+(R^\vartheta)^2},\quad \mathcal{J}_{3}^+=|\mu^\vartheta |\sqrt{1+(S^\vartheta)^2}, \quad
\mathcal{J}_{4}^-=|r^\vartheta|, \quad \mathcal{J}_{4}^+=|s^\vartheta|, \\
&\mathcal{J}_{5}^-= \fr{|r^\vartheta|}{\sqrt{1+(R^\vartheta)^2}},\quad  \mathcal{J}_{5}^+= \fr{|s^\vartheta|}{\sqrt{1+(S^\vartheta)^2}}, \\
&\mathcal{J}_{6}^-=\Big|\fr{R^\vartheta r^\vartheta}{\sqrt{1+(R^\vartheta)^2}}\Big|, \quad \mathcal{J}_{6}^+=\Big|\fr{S^\vartheta s^\vartheta}{\sqrt{1+(S^\vartheta)^2}} \Big|.
\end{split}
\end{align}
Then the norms in \eqref{6.17a} and \eqref{6.17} reduce, respectively, to
\begin{align}\label{7.36}
\begin{split}
&\|(\upsilon^\vartheta, \mu^\vartheta, w^\vartheta, \tilde{r}, z^\vartheta, \tilde{s}^\vartheta)\|_{(\tau^\vartheta, u^\vartheta, R^\vartheta, S^\vartheta)}^{(i)}=\sum_{\ell=1}^4\mathcal{I}_{\ell}^\vartheta \\ =&\int_{\mathbb{R}}\Big(|\upsilon^\vartheta|\sqrt{1+(R^\vartheta)^2}(\mathcal{W}^-)^\vartheta +|\upsilon^\vartheta|\sqrt{1+(S^\vartheta)^2}(\mathcal{W}^+)^\vartheta\Big)\ {\rm d}x  \\ &+\int_{\mathbb{R}}\Big(|\mu^\vartheta|\sqrt{1+(R^\vartheta)^2}(\mathcal{W}^-)^\vartheta +|\mu^\vartheta|\sqrt{1+(S^\vartheta)^2}(\mathcal{W}^+)^\vartheta\Big)\ {\rm d}x \\ &+\int_{\mathbb{R}}\Big(|r^\vartheta|(\mathcal{W}^-)^\vartheta +|s^\vartheta|(\mathcal{W}^+)^\vartheta\Big)\ {\rm d}x,
\end{split}
\end{align}
and
\begin{align}\label{7.37}
\begin{split}
&\|(\upsilon^\vartheta, \mu^\vartheta, w^\vartheta, \tilde{r}^\vartheta, z^\vartheta, \tilde{s}^\vartheta)\|_{(\tau^\vartheta, u^\vartheta, R^\vartheta, S^\vartheta)}^{(ii)}=\sum_{\ell=1}^6\mathcal{I}_\ell \\
=& \|(\upsilon^\vartheta, \mu^\vartheta, w^\vartheta, \tilde{r}, z^\vartheta, \tilde{s}^\vartheta)\|_{(\tau^\vartheta, u^\vartheta, R^\vartheta, S^\vartheta)}^{(i)} \\ &+\int_{\mathbb{R}}\Big(\fr{|r^\vartheta|}{\sqrt{1+(R^\vartheta)^2}}(\mathcal{W}^-)^\vartheta +\fr{|s^\vartheta|}{\sqrt{1+(S^\vartheta)^2}}(\mathcal{W}^+)^\vartheta\Big)\ {\rm d}x \\ &+\int_{\mathbb{R}}\Big(\fr{|R^\vartheta r^\vartheta|}{\sqrt{1+(R^\vartheta)^2}}(\mathcal{W}^-)^\vartheta +\fr{|S^\vartheta s^\vartheta|}{\sqrt{1+(S^\vartheta)^2}}(\mathcal{W}^+)^\vartheta\Big)\ {\rm d}x.
\end{split}
\end{align}

We next calculate $\upsilon^\vartheta, \mu^\vartheta$ and $(r^\vartheta, s^\vartheta)$ by \eqref{7.33} and \eqref{7.34}. First, we obtain
\begin{align}\label{7.38}
\begin{split}
r^\vartheta =\fr{{\rm d}}{{\rm d}\vartheta} R_{0}^\vartheta =\hat{R}_{0}-R_0, \qquad
s^\vartheta =\fr{{\rm d}}{{\rm d}\vartheta} S_{0}^\vartheta =\hat{S}_{0}-S_0.
\end{split}
\end{align}
Moreover, one has
\begin{align*}
\begin{split}
\upsilon^\vartheta =&\fr{{\rm d}}{{\rm d}\vartheta} \tau_{0}^\vartheta = -\fr{2}{\gamma-1}\Big((1-\vartheta)\tau_{0}^{-\fr{\gamma-1}{2}} +\vartheta\hat{\tau}_{0}^{-\fr{\gamma-1}{2}}  \Big)^{-\fr{\gamma+1}{\gamma-1}} \Big(\hat{\tau}_{0}^{-\fr{\gamma-1}{2}} -\tau_{0}^{-\fr{\gamma-1}{2}}\Big),
\end{split}
\end{align*}
which implies that
\begin{align}\label{7.39}
\begin{split}
\fr{1}{\overline{C}}|\hat{\tau}_{0}-\tau_0|\leq |\upsilon^\vartheta| \leq \overline{C}|\hat{\tau}_{0}-\tau_0|,
\end{split}
\end{align}
for some constant $\overline{C}>0$ depending only on $\bar{c}$. Finally, it holds
\begin{align}\label{7.40}
\begin{split}
\mu^\vartheta =&\fr{{\rm d}}{{\rm d}\vartheta} u_{0}^\vartheta = \hat{u}_{0}-u_0.
\end{split}
\end{align}
Putting \eqref{7.38}-\eqref{7.40} into \eqref{7.36} and \eqref{7.37} obtains for $\iota=i, ii$
\begin{align*}
\begin{split}
&\|(\upsilon^\vartheta, \mu^\vartheta, w^\vartheta, \tilde{r}, z^\vartheta, \tilde{s}^\vartheta)\|_{(\tau^\vartheta, u^\vartheta, R^\vartheta, S^\vartheta)}^{(\iota)}  \\ \leq &\int_{\mathbb{R}} \overline{C}|\tau_0-\hat{\tau}_0|\Big(1+|R_{0}^\vartheta| +|S_{0}^\vartheta|\Big)  \cdot 2\eta\ {\rm d}x  +\int_{\mathbb{R}} |u_0-\hat{u}_0|\Big(1+|R_{0}^\vartheta| +|S_{0}^\vartheta|\Big)  \cdot 2\eta\ {\rm d}x \\ &+2\int_{\mathbb{R}}\Big(|R_0-\hat{R}_{0}| +|S_0-\hat{S}_{0}|\Big)\cdot 2\eta\ {\rm d}x \\
\leq & \overline{C}\eta \Big(1+\|(\tau_0, u_0)\|_{C^1} +\|(\hat{\tau}_0, \hat{u}_0)\|_{C^1}\Big) \int_{\mathbb{R}}\big(|\tau_0-\hat{\tau}_0| +|R_0-\hat{R}_{0}| +|S_0-\hat{S}_{0}|\big)\ {\rm d}x \\
\leq &\widehat{C}\Big(\|\tau_0-\hat{\tau}_0\|_{W^{1,1}} +\|u_0-\hat{u}_0\|_{W^{1,1}}\Big),
\end{split}
\end{align*}
which leads to \eqref{7.32}. This finishes the proof of the proposition.
\end{proof}

\begin{prop}\label{prop3}
Let $(\tau_0, u_0)$ and $(\hat{\tau}_0, \hat{u}_0)$ be any two initial data
in $\mathcal{D}_\eta$, and $(\tau, u)$ and $(\hat{\tau}, \hat{u})$ be the corresponding solutions of \eqref{1.1}-\eqref{1.3}. Denote by $t^*$ and $\hat{t}^*$ the corresponding times at which their first singularities form. Let $\eta>0$ be a sufficiently small constant depending only on the initial data.
For any time $t\leq t_d:=\min\{t^*, \hat{t}^*\}$, there holds
\begin{align}\label{7.41}
\begin{split}
\|\tau-\hat{\tau}\|_{L^1} + \|u-\hat{u}\|_{L^1}
\leq \widehat{C} d_{*}^{(\iota)}((\tau, u), (\hat{\tau}, \hat{u})),\quad (\iota=i, ii),
\end{split}
\end{align}
for some positive constant $\widehat{C}$.
\end{prop}
\begin{proof}
We directly calculate
\begin{align}\label{7.42}
\begin{split}
\|\tau-\hat{\tau}\|_{L^1} =&\int_{R}\Big|\int_{0}^1\fr{{\rm d}\tau^\vartheta}{{\rm d}\vartheta}\ {\rm d}\vartheta\Big|\ {\rm d}x \\
\leq & \inf_{\Upsilon^t}\int_{0}^1\int_{\mathbb{R}}\Big|\fr{{\rm d}\tau^\vartheta}{{\rm d}\vartheta}\Big|\ {\rm d}x\ {\rm d}\vartheta
= \inf_{\Upsilon^t}\int_{0}^1\int_{\mathbb{R}}|\upsilon^\vartheta|\ {\rm d}x\ {\rm d}\vartheta,
\end{split}
\end{align}
and
\begin{align}\label{7.43}
\begin{split}
\|u-\hat{u}\|_{L^1}=\int_{R}\Big|\int_{0}^1\fr{{\rm d}u^\vartheta}{{\rm d}\vartheta}\ {\rm d}\vartheta\Big|\ {\rm d}x \leq \inf_{\Upsilon^t}\int_{0}^1\int_{\mathbb{R}}|\mu^\vartheta|\ {\rm d}x\ {\rm d}\vartheta.
\end{split}
\end{align}
Moreover, one has
\begin{align}\label{7.44}
\begin{split}
|\upsilon^\vartheta|\leq & \Big|\upsilon^\vartheta +\fr{R^\vartheta w^\vartheta}{2c^\vartheta}-\fr{S^\vartheta z^\vartheta}{2c^\vartheta}\Big| +\fr{|R^\vartheta|\cdot|w^\vartheta|}{2c^\vartheta} +\fr{|S^\vartheta|\cdot|z^\vartheta|}{2c^\vartheta} \\
\leq & \Big|\upsilon^\vartheta +\fr{R^\vartheta w^\vartheta}{2c^\vartheta}-\fr{S^\vartheta z^\vartheta}{2c^\vartheta}\Big| +\fr{|w^\vartheta|\sqrt{1+(R^\vartheta)^2}}{2c^\vartheta} +\fr{|z^\vartheta|\sqrt{1+(S^\vartheta)^2}}{2c^\vartheta},
\end{split}
\end{align}
and
\begin{align}\label{7.45}
\begin{split}
|\mu^\vartheta|\leq  \Big|\mu^\vartheta +\fr{R^\vartheta w^\vartheta +S^\vartheta z^\vartheta}{2}\Big| +\fr{|w^\vartheta|\sqrt{1+(R^\vartheta)^2}}{2} +\fr{|z^\vartheta|\sqrt{1+(S^\vartheta)^2}}{2}.
\end{split}
\end{align}
Combining \eqref{7.42}-\eqref{7.42} and the definition of $d_{*}^{(\iota)}(\cdot, \cdot)$ arrives at the desired estimate \eqref{7.41}. The proof of the proposition is complete.
\end{proof}

\section*{Data Availability Statement}
No data was used for the research described in the paper.

\section*{Statements and Declarations}

The authors declare that they have no conflict of interest.

\section*{Acknowledgements}

G. Chen was partially supported by National Science Foundation (DMS-2306258, DMS-2605028), Y. Hu was partially supported by National Natural Science Foundation of China (12171130) and Natural Science Foundation of Zhejiang province of China (LMS25A010014), and Y. Shen was
partially supported by National Science Foundation (DMS-2206218).

\end{document}